\pdfoutput=1
\documentclass[onefignum,onetabnum]{siamart190516}
\usepackage[utf8]{inputenc}
\usepackage{graphicx}
\usepackage{float}
\usepackage{placeins}
\usepackage{braket,amsfonts}
\usepackage{array}

\usepackage{pgfplots}
\pgfplotsset{compat=1.18}

\usepackage{diagbox}
\usepackage{amssymb}
\usepackage{amsmath,bm}
\usepackage{multirow}
\usepackage{makecell}
\usepackage{graphicx}

\usepackage[compatibility=false]{caption}
\usepackage{subcaption}
\usepackage[titletoc, title]{appendix}
\usepackage{booktabs} 
\usepackage{lipsum}
\usepackage{amsfonts}
\usepackage{epstopdf}
\usepackage{algorithmic}
\usepackage{amsopn}
\usepackage{mathrsfs}  
\usepackage{bm}
\newcommand{\vertiii}[1]{{\left\vert\kern-0.25ex\left\vert\kern-0.25ex\left\vert #1 
		\right\vert\kern-0.25ex\right\vert\kern-0.25ex\right\vert}}
\newsiamthm{claim}{Claim}
\newsiamremark{remark}{Remark}
\newsiamremark{assumption}{Assumption}
\newsiamremark{hypothesis}{Hypothesis}
\crefname{hypothesis}{Hypothesis}{Hypotheses}

\usepackage{algorithmic}

\usepackage{graphicx,epstopdf}

\Crefname{ALC@unique}{Line}{Lines}

\usepackage{amsopn}

\usepackage{xspace}
\usepackage{bold-extra}
\usepackage[most]{tcolorbox}

\colorlet{texcscolor}{blue!50!black}
\colorlet{texemcolor}{red!70!black}
\colorlet{texpreamble}{red!70!black}
\colorlet{codebackground}{black!25!white!25}

\usepackage{xcolor}
\graphicspath{{./}{numerics_rev/}{new_numerics/}{stability_tests/}}

\DeclareTotalTCBox{\code}{ v O{} }
{
	fontupper=\ttfamily\color{black},
	nobeforeafter,
	tcbox raise base,
	colback=codebackground,colframe=white,
	top=0pt,bottom=0pt,left=0mm,right=0mm,
	leftrule=0pt,rightrule=0pt,toprule=0mm,bottomrule=0mm,
	boxsep=0.5mm,
	#2}{#1}

\title{An Efficient Long-Time Stable Second-Order Ensemble Method for the Parameterized Navier--Stokes Equations}

\author{
	Daozhi Han
	\thanks{
		Department of Mathematics,
		The State University of New York at Buffalo,
		Buffalo, NY 14260, USA
		(\email{daozhiha@buffalo.edu}).
	}
	\and
	Jiaxuan Liu
	\thanks{
		School of Mathematics and Statistics,
		Central South University,
		Changsha, Hunan 410083, China
		(\email{jiaxuanliu@csu.edu.cn}).
	}
	\and
	Xiaoming Wang
	\thanks{
		School of Mathematics Science,
		Eastern Institute of Technology,
		Ningbo, Zhejiang 315200, China
		(\email{wxm@eitech.edu.cn}).
	}
	\and
	Zhoushun Zheng
	\thanks{
		School of Mathematics and Statistics,
		Central South University,
		Changsha, Hunan 410083, China
		(\email{zszheng@csu.edu.cn}).
	}
}

\headers
{A LONG-TIME STABLE ENSEMBLE METHOD FOR NAVIER--STOKES}
{D. Han, J. Liu, X. Wang, and Z. Zheng}

\ifpdf
\hypersetup{ pdftitle={An Efficient Long-Time Stable Second-Order Ensemble Method for the Parameterized Navier-Stokes Equations} }
\fi

\usepackage{hyperref}

\usepackage{float}
\makeatletter
\AtBeginDocument{%
  \def\refstepcounter@noarg#1{%
    \cref@old@refstepcounter{#1}%
    \cref@constructprefix{#1}{\cref@result}%
    \@ifundefined{cref@#1@alias}
      {\def\@tempa{#1}}
      {\def\@tempa{\csname cref@#1@alias\endcsname}}%
    \protected@edef\cref@currentlabel{%
      [\@tempa][\arabic{#1}][\cref@result]%
      \csname p@#1\endcsname\csname the#1\endcsname}}%
  \def\refstepcounter@optarg[#1]#2{%
    \refstepcounter@noarg{#2}%
    \@ifundefined{cref@#1@alias}
      {\def\@tempa{#1}}
      {\def\@tempa{\csname cref@#1@alias\endcsname}}%
    \protected@edef\cref@currentlabel{%
      \expandafter\cref@override@label@type
      \cref@currentlabel\@nil{\@tempa}}}%
}
\makeatother

\begin{document}
	\maketitle
	\begin{abstract}

Ensemble simulations for uncertainty quantification, ensemble
forecasting, and the computation of long-time statistics run over long
time horizons, where stability must hold uniformly in time and without
restrictive time-step conditions. We propose and analyze a second-order
ensemble method for parameterized incompressible Navier--Stokes
equations that combines three properties: all members share one
time-independent Brinkman (shifted Stokes) operator, the scheme is dissipative
uniformly in time for every time-step size, and its admissibility
condition on the shared viscosity is sharp. The method uses the
second-order backward differentiation formula (BDF2), a mean-reverting
scalar auxiliary variable (mr-SAV) reformulation of convection, and a
viscosity splitting that treats a shared viscosity $\bar\nu$ implicitly
and the deviations $\nu_j-\bar\nu$ explicitly. For constant sample
viscosities satisfying $0<\nu_j<\frac43\bar\nu$, a filtered energy
estimate yields uniform-in-time stability for every time-step size and
uniformly bounded forcing, with an absorbing radius independent of the
initial data and the time step. For a given $\bar\nu$, the upper bound
is sharp. The unconditional long-time stability extends to spatially
varying random viscosities, under a signed pointwise condition that
keeps $\nu_j$ within $(\bar\nu/3,\frac43\bar\nu)$, bounded away from the
endpoints uniformly in space, and on the full range
$0<\nu_j<\frac43\bar\nu$ for proportional spatial profiles. Under the
full-range condition for constant viscosities, and for sufficiently
regular solutions in two space dimensions, we prove second-order temporal
error estimates for the velocity, the auxiliary variable, and the discrete
time primitive of the pressure error. Numerical experiments confirm the
long-time stability theory, for large time steps and in the sharpness of
the viscosity condition. They also show that the method is accurate,
member by member and in its ensemble statistics, that mean reversion
keeps long-time computations robust, and that the shared operator halves
the computing time compared with member-by-member computation.

\end{abstract}

\begin{keywords}
	Navier--Stokes equations, ensemble method, mean-reverting scalar auxiliary variable, BDF2, long-time stability, error analysis, parameterized flows
	\end{keywords}

\begin{AMS}
	35Q30, 65C05, 65M12, 65M15, 65M60, 76D05
\end{AMS}

\section{Introduction}
\label{sec:introduction}

Many problems in fluid dynamics require the solution of the
incompressible Navier--Stokes (NS) equations for a collection of initial
conditions, body forces, boundary data, or physical parameters, that is,
the computation of an ensemble of flows. Such ensembles arise in
uncertainty quantification by sampling-based and other nonintrusive
methods, such as Monte Carlo sampling, stochastic collocation, and
nonintrusive polynomial chaos \cite{XiuHesthaven2005,Hosder2006UQ}, as
well as in sensitivity analysis, ensemble forecasting and data
assimilation, and parameter studies in design and control. If an
ensemble of $J$ members is computed member by member, a linear system
with its own coefficient matrix must be solved for each member at every
time step, and for large ensembles the cost quickly becomes
prohibitive. The ensemble time-stepping approach of Jiang
and Layton \cite{JiangLayton2014} reduces this cost by discretizing all
members with one common implicit operator: the convective term is split
into a part involving the ensemble-mean velocity, treated implicitly, and
a fluctuation part, treated explicitly, so that at every time step the $J$
members share a single coefficient matrix and differ only in their
right-hand sides. One factorization or preconditioner per time step then
serves the whole ensemble, and the $J$ right-hand sides can be processed
by block solvers. This idea has since been developed in many directions,
including higher-order time stepping \cite{Jiang2015,GJW-CMAM2019},
ensemble eddy-viscosity regularizations for higher Reynolds numbers
\cite{JiangLayton2015}, natural convection and magnetohydrodynamics
\cite{Fiordilino2018,MohebujjamanRebholz2017}, reduced-order modeling
\cite{GunzburgerJiangSchneier2017}, and parameterized problems
\cite{GunzburgerJiangWang2019IMA,LuoWang2018}.

Two issues motivate the present work. First, implicit diffusion gives
member-dependent matrices when the sample viscosities differ. Viscosity
splitting $\nu_j=\bar\nu+\nu_j'$, with the ensemble mean $\bar\nu$
treated implicitly and the deviation $\nu_j'$ explicitly, restores a
shared matrix \cite{GunzburgerJiangWang2019IMA,GJW-CMAM2019}. However,
the explicitly treated velocity fluctuation introduces a CFL-type
constraint depending on the mesh size and the velocity spread, which in
existing analyses interacts with the explicit viscosity term. For the
second-order scheme of \cite{GJW-CMAM2019}, the sufficient stability
conditions are $\max_j|\nu_j'|/\bar\nu\le\sqrt\mu/3$, with a parameter
$0\le\mu<1$, together with a CFL-type condition whose bound is at most
$1-\sqrt\mu$; as the ratio approaches $1/3$, the permitted time step
therefore shrinks for a fixed mesh and fixed nonzero velocity
fluctuations. In the numerical tests of \cite{GJW-CMAM2019}, the
threshold $1/3$ appeared to be sharp. Such restrictions are particularly
inconvenient in long-time computations, where the velocity spread may
grow.

Second, long-time and statistical computations call for stability
estimates that are uniform in time and dissipative, that is, of
absorbing-ball type with a radius independent of the initial data.
Together with compactness and continuity properties, such estimates are
basic ingredients in the analysis of global attractors and invariant
measures and in the approximation of long-time statistics
\cite{Wang2010a,Wang2016}. For a single NS flow, they are available for
semi-implicit schemes under time-step restrictions
\cite{GTWWW2012,Wang2012}. For ensembles, the first-order
eddy-viscosity scheme EEV2 of Jiang and Layton \cite{JiangLayton2015} is
unconditionally energy stable for a sufficiently large eddy-viscosity
parameter, with a uniform-in-time bound for homogeneous boundary data;
however, its shared operator changes in time with the ensemble, and for
finite element velocities that are not exactly divergence-free a
time-step condition remains. The ensemble scheme of Jiang and Yang
\cite{JiangYang2021SISC}, based on a stabilized scalar auxiliary variable
(SAV) for the kinetic energy, keeps one coefficient matrix fixed in time
and is stable without a time-step restriction under a
viscosity-fluctuation condition. With nontrivial forcing, however, its
modified energy is not dissipative: the bounds grow with the accumulated
forcing, and long-time stability does not follow. Related SAV ensemble
methods exist for Stokes--Darcy flow \cite{JiangYang2021CMAME}.

In this paper we propose a second-order ensemble scheme whose central
feature is a combination of three properties that existing ensemble
schemes do not provide together. First, all members share one shifted
Stokes (Brinkman) operator that is fixed in time, so that one
factorization or preconditioner serves all members and all time steps,
and, for a deterministic shared viscosity, also ensembles of different
sizes. Second, the scheme is dissipative uniformly in time for every
time-step size: its bounds are of absorbing-ball type, with a radius
independent of the initial data and of the time step. Third, its only
condition, $0<\nu_j<\frac43\bar\nu$ for the sample viscosities relative
to the implicitly treated viscosity $\bar\nu$, is sharp. This condition
restricts the choice of $\bar\nu$, not the ensemble, since any constant
$\bar\nu>\frac34\max_j\nu_j$ is admissible. A wide spread
$\max_j\nu_j/\min_j\nu_j$ thus costs larger constants and weaker
damping of stiff modes, but not stability (Remark~\ref{rem:nubar-choice}). The key
ingredient is the mean-reverting scalar auxiliary variable (mr-SAV)
reformulation of convection, also used in \cite{HanWang-M3AS} for
geophysical models in vorticity form. The SAV approach was introduced for
gradient flows \cite{ShenXuYang2018,ShenXuYang2019} and has been adapted
to the NS equations in several forms
\cite{LinYangDong2019,LiShen2020,LiShenLiu2022,HuangShen2021}. Here a
scalar $q_j$ with $q_j(0)=1$ multiplies the convective term of the $j$th
member and solves $q_j'+\gamma(q_j-1)=b(\bm u_j,\bm u_j,\bm u_j)$, where
$b$ is the convective trilinear form and $\gamma>0$ a relaxation
parameter. For divergence-free velocities the right-hand side vanishes,
so $q_j\equiv1$ and the NS equations are recovered; at the discrete level
the relaxation term supplies the dissipation $\gamma|q_j|^2$ that controls
$q_j$ uniformly in time. Mean reversion also pulls $q_j$ back toward its
exact value $1$: without it, the deviation $q_j-1$, which scales the
computed convection, accumulates over long times
(Section~\ref{subsubsec:kolmogorov-mean-reversion}); mean reversion
prevents this at no extra cost and without altering the continuous
problem. Combined with
BDF2 and the viscosity splitting, this gives the mr-SAV--BDF2 ensemble
scheme \eqref{scheme:sav-bdf2-ensemble}. It treats $\bar\nu$ implicitly
and extrapolates the viscosity deviation and both velocity arguments of
convection, so the coupling between the velocity and $q_j^{n+1}$ is
linear, and each step reduces to two solves per member with the shared
operator.

Our main results, Theorem~\ref{thm-long-time-stability} and
Corollary~\ref{cor:absorbing-radius}, prove
uniform-in-time stability for every time-step size $k>0$ and uniformly
bounded forcing under $-\bar\nu<\nu_j'<\bar\nu/3$, with physical and
ensemble bounds and an absorbing radius independent of the initial data
and of the time step. The proof combines a modified BDF2 $G$-identity, in
the spirit of \cite{Dahlquist1976,HairerWanner1996}, with a filtered
velocity test and a Lyapunov energy for the resulting stable filter. The
SAV cancellation of the convective terms survives the filtered test
because of the skew-symmetry of the trilinear form for divergence-free
velocities. For a given $\bar\nu$ the upper bound is sharp: already for
the Stokes equations, $\nu_j'>\bar\nu/3$ amplifies high-frequency modes
for every $k>0$, and in a disk the full scheme then has exponentially
growing solutions (Lemma~\ref{prop:sharp-threshold}). Thus the threshold
$1/3$ identified numerically in \cite{GJW-CMAM2019} cannot be improved
for positive deviations. The unconditional long-time stability extends
to spatially varying random viscosities
(Theorem~\ref{thm:random-viscosity-stability}), under a signed pointwise
condition that keeps $\nu_j(x)$ within $(\bar\nu(x)/3,\frac43\bar\nu(x))$,
bounded away from the endpoints uniformly in $x$, and on the full range
for proportional sample fields. The analysis is
carried out in the time-discrete, space-continuous setting with
homogeneous Dirichlet conditions; it extends to spatial discretizations
that produce exactly divergence-free velocities or use a skew-symmetric
convective form. Numerical experiments confirm the long-time stability
theory: when the viscosity condition holds, the computed solutions remain
bounded for all tested time steps, up to $k=10$; when it is violated,
they blow up once the time step exceeds the predicted critical value
(Section~\ref{sec:large-steps}). They also show that the method is
accurate, member by member and in its ensemble statistics, and that mean
reversion keeps long-time computations robust
(Sections~\ref{subsec:convergence-tests}--\ref{sec:Kolmogorovflow}).
Computing all members with the shared operator is about twice as fast as
computing them one by one with the same scheme, and about eleven times as
fast as with a standard linearly implicit scheme
(Section~\ref{exp:sec:sampling}).

The finite-time error analysis in Appendix~\ref{app:error-analysis}, for
constant viscosities in two space dimensions, uses the full-range
condition: a gradient energy for the stable filter controls the
nonlinear error terms and yields second-order temporal convergence
without a time-step restriction; its extension to spatially varying
viscosities is only outlined (Remark~\ref{rem:spatial-error}). The error
estimates concern the temporal discretization only.

The rest of the article is organized as follows.
Section~\ref{sec:sav-bdf2-ensemble-stability} introduces the problem, the
mr-SAV reformulation, the scheme and its decoupled implementation, and
proves unique solvability. Section~\ref{sec:long-time-stability} proves
the uniform-in-time estimates and the sharpness of the threshold $1/3$,
and extends the estimates to spatially varying viscosities.
Section~\ref{numerics} presents
the numerical experiments and Section~\ref{sec:conclusions} concludes.
Appendix~\ref{app:error-analysis} treats the backward Euler startup step
and contains the finite-time error analysis.

\section{The ensemble method}
\label{sec:sav-bdf2-ensemble-stability}

\subsection{The parameterized Navier--Stokes problem}
\label{subsec:model-functional-setting}

Let $\Omega\subset\mathbb R^d$, $d=2,3$, be a bounded connected Lipschitz
domain. For each ensemble member $j=1,\ldots,J$, we seek the velocity
$\bm u_j:\Omega\times[0,\infty)\to\mathbb R^d$ and the pressure
$p_j:\Omega\times (0,\infty)\to\mathbb R$ satisfying
\begin{subequations}
\label{eq:parameterized-nse}
\begin{align}
    &\partial_t\bm u_j
    +(\bm u_j\cdot\nabla)\bm u_j
    -\nabla\cdot(\nu_j\nabla\bm u_j)
    +\nabla p_j
    =
    \bm F_j,
    \quad
    \nabla\cdot\bm u_j=0
    \;\;\text{in }\Omega\times(0,\infty),
    \label{eq:parameterized-nse-system}
    \\
    &\qquad\qquad\qquad\quad\bm u_j
    =
    \bm 0
    \;\;\text{on }\partial\Omega\times(0,\infty),
    \qquad\qquad\qquad
    \bm u_j(\cdot,0)=\bm u_j^0
    \quad\text{in }\Omega .
    \label{eq:parameterized-nse-data}
\end{align}
\end{subequations}
Here $\nu_j>0$ denotes the constant viscosity associated with the
$j$th ensemble member (spatially varying viscosities are treated in
Section~\ref{subsec:random-viscosity}), $\bm F_j$ is the prescribed body force, and
$\bm u_j^0$ is the corresponding initial velocity.

We use the velocity and pressure spaces
\begin{equation}
\label{def:velocity-pressure-spaces}
  \bm X:=H_0^1(\Omega)^d,
  \qquad
  Q:=L_0^2(\Omega)
  :=
  \left\{
    r\in L^2(\Omega):(r,1)=0
  \right\},
\end{equation}
and define the solenoidal subspace
\begin{equation}
\label{def:divergence-free-space}
  \bm V
  :=
  \left\{
    \bm v\in\bm X:
    (\nabla\cdot\bm v,r)=0
    \quad\forall r\in Q
  \right\}.
\end{equation}
The initial velocities are assumed to satisfy
\begin{equation}
\label{assumption:compatible-initial-data}
  \bm u_j^0\in\bm V,
  \qquad j=1,\ldots,J.
\end{equation}

The $L^2(\Omega)$ inner product and norm are denoted by
$(\cdot,\cdot)$ and $\|\cdot\|$, respectively. For
$\bm f\in\bm X'$, define
\begin{equation}
\label{def:dual-norm}
  \|\bm f\|_{-1}
  :=
  \sup_{\bm v\in\bm X\setminus\{\bm0\}}
  \frac{\langle\bm f,\bm v\rangle}
       {\|\nabla\bm v\|},
\end{equation}
where $\langle\cdot,\cdot\rangle$ denotes the duality pairing between
$\bm X'$ and $\bm X$, so that
\begin{equation}
\label{eq:duality-estimate}
  |\langle\bm f,\bm v\rangle|
  \le
  \|\bm f\|_{-1}\|\nabla\bm v\|,
  \qquad \bm v\in\bm X.
\end{equation}
We also write $\|\bm f\|_{\bm X'}:=\|\bm f\|_{-1}$.

Since $\bm X=H_0^1(\Omega)^d$, Poincar\'e's inequality gives
\begin{equation}
\label{eq:poincare}
  \|\nabla\bm v\|^2
  \ge c_0\|\bm v\|^2,
  \qquad
  \bm v\in\bm X,
\end{equation}
for some constant $c_0>0$ depending only on $\Omega$. Moreover, the
pair $(\bm X,Q)$ satisfies the continuous inf--sup condition
\begin{equation}
\label{eq:continuous-infsup}
  \inf_{r\in Q\setminus\{0\}}
  \sup_{\bm v\in\bm X\setminus\{\bm0\}}
  \frac{(r,\nabla\cdot\bm v)}
       {\|r\|\,\|\nabla\bm v\|}
  \ge \beta_0>0.
\end{equation}
For a conforming pair $(\bm X_h,Q_h)\subset(\bm X,Q)$, \eqref{eq:poincare}
holds on $\bm X_h$ with the same $c_0$; uniformly inf--sup stable pairs
satisfy, by definition, the discrete inf--sup condition with an
$h$-independent constant.

For $\bm w,\bm u,\bm v\in\bm X$, define
\begin{equation}
\label{def:trilinear-form}
  b(\bm w,\bm u,\bm v)
  :=
  \int_\Omega
  (\bm w\cdot\nabla)\bm u\cdot\bm v\,dx.
\end{equation}
For $d=2,3$, H\"older's inequality, the Sobolev embeddings, and
Poincar\'e's inequality yield
\begin{equation}
\label{eq:trilinear-continuity}
  |b(\bm w,\bm u,\bm v)|
  \le
  C_b
  \|\nabla\bm w\|
  \|\nabla\bm u\|
  \|\nabla\bm v\|,
  \qquad
  \bm w,\bm u,\bm v\in\bm X,
\end{equation}
where $C_b>0$ depends only on $\Omega$ and $d$. If $\bm w\in\bm V$,
integration by parts gives
\begin{equation}
\label{eq:skew-general}
  b(\bm w,\bm u,\bm v)=-b(\bm w,\bm v,\bm u),
  \qquad
  b(\bm w,\bm v,\bm v)=0,
  \qquad \bm u,\bm v\in\bm X;
\end{equation}
in particular, $b(\bm u,\bm u,\bm u)=0$ for $\bm u\in\bm V$.

To split the viscous term, we fix a constant $\bar\nu>0$, the implicitly
treated viscosity shared by all members, and set
\begin{equation}
\label{def:viscosity-splitting}
  \nu_j':=\nu_j-\bar\nu,
  \qquad
  \eta_j:=-\frac{\nu_j'}{\bar\nu}=1-\frac{\nu_j}{\bar\nu},
  \qquad
  \rho:=\max_{1\le j\le J}|\nu_j'|.
\end{equation}
For the main results we assume
\begin{equation}
\label{assumption:viscosity-fluctuation}
 0<\nu_j<\frac43\,\bar\nu,
 \qquad j=1,\ldots,J,
\end{equation}
equivalently $-\bar\nu<\nu_j'<\bar\nu/3$, or $\eta_j\in(-\frac13,1)$.
The lower bound is the positivity of the viscosities. The upper bound,
$\bar\nu>\frac34\max_j\nu_j$, restricts the choice of $\bar\nu$, not the
ensemble, and is sharp for a given $\bar\nu$
(Lemma~\ref{prop:sharp-threshold}). For the ensemble mean
$\bar\nu=\frac1J\sum_{i=1}^J\nu_i$, used in most experiments of
Section~\ref{numerics} but not required by the analysis, it reads
$\max_j\nu_j<\frac43\cdot\frac1J\sum_{i=1}^J\nu_i$; for the midrange
$\frac12(\max_j\nu_j+\min_j\nu_j)$, it reads $\max_j\nu_j<2\min_j\nu_j$.
The constants below depend on the viscosities only through $\bar\nu$ and
the ratios $\nu_j/\bar\nu$, and deteriorate as a ratio approaches $0$ or
$\frac43$; bounds stated uniformly over a family of ensembles assume that
the ratios stay in a fixed compact subinterval of $(0,\frac43)$. If
$\bar\nu$ is computed from the samples, all results hold for each
realization, and groupwise if the ensemble is split into groups with
their own shared viscosities (Remark~\ref{rem:nubar-choice}(iv)).

Let $k>0$ denote the time-step size and $t^n:=nk$. The scheme uses a
discrete forcing sequence $\{\bm F_j^{n+1}\}_{n\ge0}\subset\bm X'$, which
we assume to be uniformly bounded:
\begin{equation}
\label{assumption:discrete-force}
  M_F
  :=
  \max_{1\le j\le J}
  \sup_{n\ge0}
  \|\bm F_j^{n+1}\|_{-1}
  <\infty.
\end{equation}
Whenever a bound is asserted uniformly over a family of time-step sizes,
the same $M_F$ is assumed for all of them. If
$\bm F_j\in L^\infty(0,\infty;\bm X')$, the time-cell average
\begin{equation}
\label{def:discrete-force-average}
  \bm F_j^{n+1}
  :=
  \frac1k
  \int_{t^n}^{t^{n+1}}\bm F_j(s)\,ds
\end{equation}
satisfies \eqref{assumption:discrete-force} with
$M_F\le\max_j\|\bm F_j\|_{L^\infty(0,\infty;\bm X')}$ for every $k$.
When the temporal regularity of $\bm F_j$ permits pointwise evaluation,
one may instead set $\bm F_j^{n+1}:=\bm F_j(t^{n+1})$.

\subsection{\texorpdfstring{Mean-reverting scalar auxiliary variable reformulation}{Mean-reverting scalar auxiliary variable reformulation}}
\label{subsec:fsav-reformulation}

Let $\gamma>0$ be prescribed. For each ensemble member
$j=1,\ldots,J$, introduce a scalar auxiliary variable
$q_j:[0,\infty)\to\mathbb R$ satisfying $q_j(0)=1$. The corresponding
mean-reverting scalar auxiliary variable (mr-SAV) formulation is given by
\begin{subequations}
\label{eq:continuous-fsav-system}
\begin{align}
  \partial_t\bm u_j
  +q_j(\bm u_j\cdot\nabla)\bm u_j
  -\nabla\cdot(\nu_j\nabla\bm u_j)
  +\nabla p_j
  &=
  \bm F_j,
  \qquad
  \nabla\cdot\bm u_j=0,
  \label{eq:continuous-fsav-system1}
  \\
  q_j'
  +\gamma q_j
  -b(\bm u_j,\bm u_j,\bm u_j)
  &=
  \gamma,
  \label{eq:continuous-fsav-system2}
\end{align}
\end{subequations}
with the boundary and initial conditions of
\eqref{eq:parameterized-nse}. The relaxation term $\gamma(1-q_j)$ drives
$q_j$ toward its exact value $1$, whence the name; at the discrete level
its linear part produces the dissipation $\gamma|q_j|^2$ that controls the
auxiliary variable uniformly in time
(Section~\ref{sec:long-time-stability}). For a sufficiently regular
solution with $\bm u_j\in\bm V$, \eqref{eq:skew-general}
reduces \eqref{eq:continuous-fsav-system2} to $q_j'+\gamma(q_j-1)=0$, so
$q_j\equiv1$, and \eqref{eq:continuous-fsav-system} is equivalent to
\eqref{eq:parameterized-nse} within this solution class.

\subsection{The ensemble scheme}
\label{subsec:coupled-scheme}

Let $H$ be a real Hilbert space with inner product $(\cdot,\cdot)$ and
norm $\|\cdot\|$. For an $H$-valued sequence $\{w^n\}_{n\ge0}$, define the
BDF2 difference quotient and second-order extrapolation by
\begin{equation}
\label{def:bdf2-and-extrapolation}
  \delta_t w^{n+1}
  :=
  \frac{3w^{n+1}-4w^n+w^{n-1}}{2k},
  \qquad
  w^{*,n+1}
  :=
  2w^n-w^{n-1},
  \qquad n\ge1.
\end{equation}

Introduce the two-level vector $W^n:=(w^n,w^{n-1})^T$, $n\ge1$,
and define
\begin{equation}
\label{def:G-matrix}
  G
  :=
  \frac14
  \begin{pmatrix}
    5&-2\\
    -2&1
  \end{pmatrix},
  \qquad
  \|W\|_G^2
  :=
  (W,GW)_{H\times H}.
\end{equation}
Since $G$ is symmetric positive definite, $\|\cdot\|_G$ defines a norm
on $H\times H$. The choices $H=L^2(\Omega)^d$ and $H=\mathbb R$
correspond, respectively, to vector- and scalar-valued sequences, and
$H=L^2(\Omega)^{d\times d}$ to sequences of gradients. The energy
estimates rest on the following identity for the BDF2 difference
quotient and on the equivalence of $\|\cdot\|_G$ with the norm of
$H\times H$.

\begin{lemma}
\label{lem:bdf2-G-identity}
For every $H$-valued sequence $\{w^n\}_{n\ge0}$ and $n\ge1$,
\begin{equation}
\label{eq:bdf2-G-identity-preliminary}
  k(\delta_t w^{n+1},w^{n+1})
  =
  \|W^{n+1}\|_G^2
  -
  \|W^n\|_G^2
  +
  \frac14
  \|w^{n+1}-2w^n+w^{n-1}\|^2.
\end{equation}
Moreover, for every $W=(w_1,w_2)^T\in H\times H$,
\begin{equation}
\label{eq:G-norm-equivalence-preliminary}
  c_l\|W\|_G^2
  \le
  \|W\|^2
  \le
  c_u\|W\|_G^2,
\end{equation}
where $\|W\|^2:=\|w_1\|^2+\|w_2\|^2$, $c_l:=4/(3+2\sqrt2)=12-8\sqrt2$, and
$c_u:=4/(3-2\sqrt2)=12+8\sqrt2$. For $H=\mathbb R$ we write $|\cdot|_G$
instead of $\|\cdot\|_G$.
\end{lemma}
The identity follows by expansion, and
\eqref{eq:G-norm-equivalence-preliminary} from the eigenvalues
$(3\pm2\sqrt2)/4$ of $G$.

For each ensemble member $j=1,\ldots,J$, prescribe
$\bm u_j^0\in\bm V$ and $q_j^0=1$, and let a first level
$(\bm u_j^1,q_j^1)\in\bm V\times\mathbb R$ be given.
For $n\ge1$, given
$(\bm u_j^{n-1},\bm u_j^n,q_j^{n-1},q_j^n)$, determine
$(\bm u_j^{n+1},p_j^{n+1},q_j^{n+1})
\in\bm X\times Q\times\mathbb R$ such that
\begin{subequations}
\label{scheme:sav-bdf2-ensemble}
\begin{align}
  &(\delta_t\bm u_j^{n+1},\bm v)
  +q_j^{n+1}
   b(\bm u_j^{*,n+1},\bm u_j^{*,n+1},\bm v)
  +\bar\nu(\nabla\bm u_j^{n+1},\nabla\bm v)
  \nonumber\\
  &\qquad
  +\nu_j'(\nabla\bm u_j^{*,n+1},\nabla\bm v)
  -(p_j^{n+1},\nabla\cdot\bm v)
  =\langle\bm F_j^{n+1},\bm v\rangle,
  \label{scheme:sav-bdf2-momentum}
  \\
  &(\nabla\cdot\bm u_j^{n+1},r)=0,
  \label{scheme:sav-bdf2-continuity}
  \\
  &\delta_tq_j^{n+1}+\gamma q_j^{n+1}
  -b(\bm u_j^{*,n+1},\bm u_j^{*,n+1},\bm u_j^{n+1})
  =\gamma,
  \label{scheme:sav-bdf2-q}
\end{align}
\end{subequations}
where the first two equations hold for all
$(\bm v,r)\in\bm X\times Q$.

The stability results of Section~\ref{sec:long-time-stability} hold for
arbitrary starting levels $(\bm u_j^0,q_j^0)$ and $(\bm u_j^1,q_j^1)$ in
$\bm V\times\mathbb R$; the error
analysis computes the first level by the mr-SAV backward Euler step of
Appendix~\ref{subsec:be-startup}.

Only $\bar\nu$ is treated implicitly; the viscosity deviation and both
velocity arguments of convection are extrapolated. Hence, for fixed $k$
and a fixed spatial discretization, the implicit velocity--pressure
block, and thus the coefficient matrix of the Brinkman subproblems below,
is the same for all members and all BDF2 steps. Moreover, testing
\eqref{scheme:sav-bdf2-momentum} with $\bm u_j^{n+1}$ and multiplying
\eqref{scheme:sav-bdf2-q} by $q_j^{n+1}$ produces the term
$q_j^{n+1}b(\bm u_j^{*,n+1},\bm u_j^{*,n+1},\bm u_j^{n+1})$ with opposite
signs. This exact cancellation, the basis of the energy analysis, does
not require $\bm u_j^{*,n+1}$ to be divergence-free; the filtered test
of Section~\ref{subsec:asymmetric-stability} also needs
$b(\bm u_j^{*,n+1},\bm u_j^{*,n+1},\bm u_j^{*,n+1})=0$, which holds as
$\bm u_j^{*,n+1}\in\bm V$ (Proposition~\ref{prop:stepwise-solvability}).

\subsection{Efficient implementation}
\label{subsec:decoupling-solvability}

For each member, a step of the scheme reduces to one scalar
equation and two linear Brinkman subproblems. Define
\begin{equation}
\label{def:brinkman-bilinear-form}
  a_k(\bm w,\bm v)
  :=\frac{3}{2k}(\bm w,\bm v)
  +\bar\nu(\nabla\bm w,\nabla\bm v),
\end{equation}
and, for each $j$ and $n\ge1$, introduce
\begin{align}
  \langle\bm f_{1j}^{n+1},\bm v\rangle
  &:=\langle\bm F_j^{n+1},\bm v\rangle
  -\nu_j'(\nabla\bm u_j^{*,n+1},\nabla\bm v)
  +\left(\frac{4\bm u_j^n-\bm u_j^{n-1}}{2k},\bm v\right),
  \label{def:brinkman-rhs-one}
  \\
  \langle\bm f_{2j}^{n+1},\bm v\rangle
  &:=-b(\bm u_j^{*,n+1},\bm u_j^{*,n+1},\bm v).
  \label{def:brinkman-rhs-two}
\end{align}
Both functionals belong to $\bm X'$, the second by
\eqref{eq:trilinear-continuity}. For $i=1,2$, let
$(\bm u_{ij}^{n+1},p_{ij}^{n+1})\in\bm X\times Q$
satisfy
\begin{subequations}
\label{eq:brinkman-subproblems-weak}
\begin{align}
  a_k(\bm u_{ij}^{n+1},\bm v)
  -(p_{ij}^{n+1},\nabla\cdot\bm v)
  &=\langle\bm f_{ij}^{n+1},\bm v\rangle,
  &&\bm v\in\bm X,
  \label{eq:brinkman-momentum-weak}
  \\
  (\nabla\cdot\bm u_{ij}^{n+1},r)&=0,
  &&r\in Q.
  \label{eq:brinkman-continuity-weak}
\end{align}
\end{subequations}
These subproblems differ only in their right-hand sides. Since
\eqref{scheme:sav-bdf2-momentum}--\eqref{scheme:sav-bdf2-continuity} is
the Brinkman problem with right-hand side
$\bm f_{1j}^{n+1}+q_j^{n+1}\bm f_{2j}^{n+1}$, linearity gives
\begin{equation}
\label{eq:affine-decomposition}
  \bm u_j^{n+1}
  =\bm u_{1j}^{n+1}+q_j^{n+1}\bm u_{2j}^{n+1},
  \qquad
  p_j^{n+1}=p_{1j}^{n+1}+q_j^{n+1}p_{2j}^{n+1}.
\end{equation}
Substitution into \eqref{scheme:sav-bdf2-q} yields
\begin{equation}
\label{eq:scalar-q-update-equation}
  D_j^{n+1}q_j^{n+1}
  =\gamma+\frac{4q_j^n-q_j^{n-1}}{2k}
  +b(\bm u_j^{*,n+1},\bm u_j^{*,n+1},\bm u_{1j}^{n+1}),
\end{equation}
where
\begin{equation}
\label{def:scalar-denominator}
  D_j^{n+1}
  :=\frac{3}{2k}+\gamma
  -b(\bm u_j^{*,n+1},\bm u_j^{*,n+1},\bm u_{2j}^{n+1}).
\end{equation}

The next proposition shows that each step is uniquely solvable and
that $D_j^{n+1}\ge3/(2k)+\gamma$, so that
\eqref{eq:scalar-q-update-equation} determines $q_j^{n+1}$.

\begin{proposition}
\label{prop:stepwise-solvability}
Let $n\ge1$, $k>0$, $\bar\nu>0$, $\gamma>0$, $\nu_j'\in\mathbb R$, and
$\bm F_j^{n+1}\in\bm X'$. For any $\bm u_j^{n-1},\bm u_j^n\in\bm X$ and
$q_j^{n-1},q_j^n\in\mathbb R$, the system
\eqref{scheme:sav-bdf2-ensemble} admits a unique solution
$(\bm u_j^{n+1},p_j^{n+1},q_j^{n+1})
\in\bm X\times Q\times\mathbb R$, and
\begin{equation}
\label{eq:positive-scalar-denominator}
  D_j^{n+1}
  =\frac{3}{2k}+\gamma
  +\frac{3}{2k}\|\bm u_{2j}^{n+1}\|^2
  +\bar\nu\|\nabla\bm u_{2j}^{n+1}\|^2
  \ge\frac{3}{2k}+\gamma.
\end{equation}
\end{proposition}

\begin{proof}

Since $a_k(\bm v,\bm v)=\frac3{2k}\|\bm v\|^2+\bar\nu\|\nabla\bm v\|^2$,
the form $a_k$ is bounded and coercive on $\bm X$; with
\eqref{eq:continuous-infsup}, mixed variational theory gives unique
solvability of \eqref{eq:brinkman-subproblems-weak} for $i=1,2$. Taking
$\bm v=\bm u_{2j}^{n+1}$ in \eqref{eq:brinkman-momentum-weak} with $i=2$,
the pressure term vanishes by \eqref{eq:brinkman-continuity-weak} with
$r=p_{2j}^{n+1}$, so that
$-b(\bm u_j^{*,n+1},\bm u_j^{*,n+1},\bm u_{2j}^{n+1})
=a_k(\bm u_{2j}^{n+1},\bm u_{2j}^{n+1})$; inserted into
\eqref{def:scalar-denominator}, this proves
\eqref{eq:positive-scalar-denominator}. Hence
\eqref{eq:scalar-q-update-equation} determines $q_j^{n+1}$ uniquely,
\eqref{eq:affine-decomposition} gives the velocity and pressure, and
direct substitution shows that they solve
\eqref{scheme:sav-bdf2-ensemble}. Conversely, by uniqueness for the
Brinkman problem, every solution of \eqref{scheme:sav-bdf2-ensemble} has
the form \eqref{eq:affine-decomposition}, and its scalar component
satisfies \eqref{eq:scalar-q-update-equation}; this proves uniqueness.

Finally, $\bm u_{1j}^{n+1},\bm u_{2j}^{n+1}\in\bm V$ by
\eqref{eq:brinkman-continuity-weak}, so every computed level lies in
$\bm V$ whatever the history; with starting levels in $\bm V$, also
$\bm u_j^{*,n+1}\in\bm V$ for all $n\ge1$.
\end{proof}

\section{Long-time stability}
\label{sec:long-time-stability}

We prove uniform-in-time stability with a time-step-independent
absorbing radius for constant sample viscosities, show that the
threshold is sharp, and extend the analysis to spatially varying
viscosities.

For $n\ge1$, collect the velocity and auxiliary variable at
two consecutive time levels into
\begin{equation}
\label{eq:stability-two-level-vectors}
  V_j^n:=
  \begin{pmatrix}
    \bm u_j^n\\
    \bm u_j^{n-1}
  \end{pmatrix},
  \qquad
  Q_j^n:=
  \begin{pmatrix}
    q_j^n\\
    q_j^{n-1}
  \end{pmatrix},
\end{equation}
and define their combined BDF2 $G$-energy by
\begin{equation}
\label{eq:basic-G-energy}
  \mathcal G_j^n
  :=\|V_j^n\|_G^2+|Q_j^n|_G^2.
\end{equation}
\subsection{\texorpdfstring{Uniform-in-time stability}{Uniform-in-time stability}}
\label{subsec:asymmetric-stability}

For constant sample viscosities, assume the sharp condition
\eqref{assumption:viscosity-fluctuation} introduced in
Section~\ref{subsec:model-functional-setting}.
For each member, with $\eta_j\in(-\frac13,1)$ from
\eqref{def:viscosity-splitting}, set
\begin{equation}
\label{eq:full-range-multiplier}
 \bm w_j^{n+1}:=\bm u_j^{n+1}-\eta_j\bm u_j^{*,n+1}.
\end{equation}
Since $\bm u_j^{*,n+1}\in\bm V$, the identity
$b(\bm u_j^{*,n+1},\bm u_j^{*,n+1},\bm u_j^{*,n+1})=0$
allows this test function to preserve the mr-SAV cancellation.

We call $\bm w_j^{n+1}$ the \emph{filtered velocity}: it is obtained from
the sequence $(\bm u_j^n)_{n\ge0}$ by the linear time filter
$\bm u^{n+1}-2\eta_j\bm u^n+\eta_j\bm u^{n-1}$, whose coefficients do not
depend on $n$. Testing the momentum equation with $\bm w_j^{n+1}$ instead
of $\bm u_j^{n+1}$ is a multiplier technique in the spirit of Nevanlinna
and Odeh \cite{NevanlinnaOdeh1981}, which has been used for energy
estimates of BDF methods, including linearly implicit and
implicit--explicit variants with extrapolated terms
\cite{AkrivisLubich2015}; here the multiplier is matched to the
explicitly treated viscosity, so that the split viscous terms combine
into $\bar\nu\|\nabla\bm w_j^{n+1}\|^2$. The energy estimate therefore
controls $\bm w_j$; the original velocity is recovered by inverting the
filter, which is done below.

For $-1/3<\eta<1$, define
\begin{equation}
\label{eq:full-range-matrices}
 G_\eta:=\frac14
 \begin{pmatrix}5-3\eta&-2\\-2&1+\eta\end{pmatrix},
 \qquad
 B_\eta:=\begin{pmatrix}2\eta&-\eta\\1&0\end{pmatrix}.
\end{equation}
The matrix $G_\eta$ is positive definite, since
$\det G_\eta=(1-\eta)(1+3\eta)/16>0$ and $(G_\eta)_{11}>0$, and, as
direct expansion shows, for every Hilbert-space
sequence,
\begin{equation}
\label{eq:full-range-G-identity}
\begin{aligned}
 &k(\delta_t u^{n+1},u^{n+1}-\eta(2u^n-u^{n-1}))\\
 &\quad=\|(u^{n+1},u^n)\|_{G_\eta}^2
       -\|(u^n,u^{n-1})\|_{G_\eta}^2
       +\frac{1+3\eta}{4}\|u^{n+1}-2u^n+u^{n-1}\|^2.
\end{aligned}
\end{equation}
The matrix $B_\eta$ describes the inverse of the filter: by
\eqref{eq:full-range-multiplier}, $\bm u_j^{n+1}=2\eta_j\bm u_j^n
-\eta_j\bm u_j^{n-1}+\bm w_j^{n+1}$, that is,
$V_j^{n+1}=B_{\eta_j}V_j^n+(\bm w_j^{n+1},\bm0)^T$, a recursive
(infinite impulse response) filter driven by $\bm w_j$. Such a filter is
stable precisely when the roots of its characteristic polynomial, here
$z^2-2\eta z+\eta$, lie in the open unit disk
\cite{OppenheimSchafer2010}, which by the Schur--Cohn conditions~\cite{Miller1971}
$|\eta|<1$, $1-\eta>0$, $1+3\eta>0$ is the case exactly for
$-\frac13<\eta<1$. Consequently, the following series converges, and its sum
satisfies the discrete Lyapunov equation
\cite{LancasterTismenetsky1985}:
\begin{equation}
\label{eq:full-range-P}
 P_\eta:=\sum_{\ell=0}^\infty(B_\eta^T)^\ell B_\eta^\ell,
 \qquad
 P_\eta-B_\eta^TP_\eta B_\eta=I,\qquad P_\eta\ge I.
\end{equation}
The sum $P_\eta$ can be computed in closed form, which is useful for
evaluating the constants below; direct substitution in
\eqref{eq:full-range-P} verifies the formula
\begin{equation}
\label{eq:full-range-P-explicit}
\begin{gathered}
 \Delta_\eta:=(1-\eta)^2(1+3\eta),\qquad
 p_\eta:=\frac{2(1+\eta)}{\Delta_\eta},\qquad
 s_\eta:=-\frac{4\eta^2}{\Delta_\eta},\\
 P_\eta=\begin{pmatrix}p_\eta&s_\eta\\
                         s_\eta&1+\eta^2p_\eta\end{pmatrix}.
\end{gathered}
\end{equation}

For $P_j:=P_{\eta_j}$, define the filter energy and the modified basic
energy by
\begin{equation}
\label{eq:full-range-basic-energies}
 \mathcal F_j^n:=\|V_j^n\|_{P_j}^2,
 \qquad
 \mathcal E_j^n:=\|V_j^n\|_{G_{\eta_j}}^2+|Q_j^n|_G^2.
\end{equation}
All matrix norms use the product-space inner product, as in
\eqref{def:G-matrix}. Introduce constants independent of $k$:
\begin{equation}
\label{eq:full-range-constants}
\begin{gathered}
 \Lambda:=\max_j\lambda_{\max}(P_j),\qquad
 p_*:=\max_j(P_j)_{11},\qquad \theta:=1-\Lambda^{-1},\\
 \sigma:=\frac{1-\theta}{1+\theta}>0,\qquad
 D:=\frac{2p_*}{1-\theta},\qquad
 \varepsilon:=\frac{\bar\nu c_0}{4D}>0,\\
 \widetilde C_F:=\frac{M_F^2}{2\bar\nu}+\frac\gamma2,\qquad
 C_G:=\max_j\lambda_{\max}(G_{\eta_j}),\qquad
 C_q:=c_l^{-1},\\
 c_G:=\min\{c_u^{-1},\min_j\lambda_{\min}(G_{\eta_j})\}>0.
\end{gathered}
\end{equation}
Here $\frac12\le\theta<1$, since $\Lambda\ge(P_{\eta_j})_{11}=p_{\eta_j}\ge2$,
as $p_\eta-2=2\eta^2(5-3\eta)/\Delta_\eta$ by \eqref{eq:full-range-P-explicit}; in particular $0<\sigma\le\frac13$. Define the modified energy and
its decay rate by
\begin{equation}
\label{eq:full-range-energy}
 \widetilde{\mathcal H}_j^n:=\mathcal E_j^n+k\varepsilon\mathcal F_j^n,
 \qquad
 \widetilde\alpha_k:=\min\left\{
 \frac{k\varepsilon\sigma}{C_G+k\varepsilon},
 \frac{k\gamma}{2C_q}\right\}>0.
\end{equation}

The main result bounds this modified energy, and hence the
solutions, uniformly in time for every $k>0$.

\begin{theorem}
\label{thm-long-time-stability}
\label{thm:asymmetric-stability}
Assume \eqref{assumption:viscosity-fluctuation} and
\eqref{assumption:discrete-force}, with constant sample viscosities and
starting velocities in $\bm V$ and scalar starting values in $\mathbb R$.
For every $k>0$, the solutions of the scheme
\eqref{scheme:sav-bdf2-ensemble} satisfy
\begin{align}
 \sup_{n\ge1}\widetilde{\mathcal H}_j^n
 &\le\max\left\{\widetilde{\mathcal H}_j^1+k\widetilde C_F,
                \frac{2k\widetilde C_F}{\widetilde\alpha_k}\right\}
 =:S_j(k),
 \label{eq:asymmetric-uniform}\\
 \limsup_{n\to\infty}\widetilde{\mathcal H}_j^n
 &\le\frac{2k\widetilde C_F}{\widetilde\alpha_k}.
 \label{eq:asymmetric-limsup}
\end{align}
Consequently, with the coercivity constant $c_G$ from
\eqref{eq:full-range-constants}, the samplewise and ensemble physical
bounds are
\begin{align}
 \sup_{n\ge1}(\|\bm u_j^n\|^2+|q_j^n|^2)&\le S_j(k)/c_G,
 \label{eq:main-uniform-physical}\\
 \sup_{n\ge1}\frac1J\sum_{j=1}^J
 (\|\bm u_j^n\|^2+|q_j^n|^2)&\le\frac1{c_G}\max_j S_j(k).
 \label{eq:main-ensemble-bound}
\end{align}
\end{theorem}

\begin{proof}
Test the momentum equation with $k\bm w_j^{n+1}$ and the scalar equation with
$kq_j^{n+1}$. The pressure term vanishes since
$\bm w_j^{n+1}\in\bm V$. Since
$b(\bm u_j^{*,n+1},\bm u_j^{*,n+1},\bm u_j^{*,n+1})=0$, the trilinear term
of the momentum equation equals
$kq_j^{n+1}b(\bm u_j^{*,n+1},\bm u_j^{*,n+1},\bm u_j^{n+1})$ and cancels
the one from the scalar equation.
Using \eqref{eq:full-range-G-identity}, the scalar $G$-identity
\eqref{eq:bdf2-G-identity-preliminary},
and Young's inequality for the forcing and scalar source yields
\begin{equation}
\label{eq:full-range-basic-step}
\begin{aligned}
 &\mathcal E_j^{n+1}-\mathcal E_j^n
 +\frac{1+3\eta_j}{4}
       \|\bm u_j^{n+1}-2\bm u_j^n+\bm u_j^{n-1}\|^2\\
 &\quad+\frac14|q_j^{n+1}-2q_j^n+q_j^{n-1}|^2
 +\frac{k\bar\nu}{2}\|\nabla\bm w_j^{n+1}\|^2
 +\frac{k\gamma}{2}|q_j^{n+1}|^2
 \le k\widetilde C_F.
\end{aligned}
\end{equation}
To obtain dissipation of the original velocity, we use the stable filter
recursion $V_j^{n+1}=B_{\eta_j}V_j^n+(\bm w_j^{n+1},\bm0)^T$.
By \eqref{eq:full-range-P} and $\|V\|_{P_j}^2\le\Lambda\|V\|^2$, we have
\[\|B_{\eta_j}V\|_{P_j}^2=\|V\|_{P_j}^2-\|V\|^2
\le\theta\|V\|_{P_j}^2.\]
Since only the first component of
$(\bm w_j^{n+1},\bm0)^T$ is nonzero, its squared $P_j$-norm is
$(P_j)_{11}\|\bm w_j^{n+1}\|^2\le p_*\|\bm w_j^{n+1}\|^2$.
Taking the $P_j$-norm of the filter recursion, the triangle inequality
and the above bound with $V=V_j^n$ give
\[
 \sqrt{\mathcal F_j^{n+1}}
 =\|V_j^{n+1}\|_{P_j}
 \le\|B_{\eta_j}V_j^n\|_{P_j}
       +\|(\bm w_j^{n+1},\bm0)^T\|_{P_j}
 \le\sqrt{\theta\mathcal F_j^n}
       +\sqrt{p_*}\|\bm w_j^{n+1}\|.
\]
By weighted Young's inequality we derive
\[
 \mathcal F_j^{n+1}
 \le\frac{1+\theta}{2}\mathcal F_j^n
      +\frac{p_*(1+\theta)}{1-\theta}\|\bm w_j^{n+1}\|^2.
\]
Dividing by $(1+\theta)/2$ and rearranging, with
$2/(1+\theta)=1+\sigma$ and $\sigma$, $D$ from \eqref{eq:full-range-constants}, yields
\begin{equation}
\label{eq:full-range-filter-dissipation}
 \mathcal F_j^{n+1}-\mathcal F_j^n
       +\sigma\mathcal F_j^{n+1}
 \le D\|\bm w_j^{n+1}\|^2.
\end{equation}
Multiply this by
$k\varepsilon$ and add it to \eqref{eq:full-range-basic-step}.
Poincar\'e's inequality and $\varepsilon D=\bar\nu c_0/4$ absorb the
right-hand side into half of the remaining viscous dissipation. After
dropping the nonnegative second differences, we obtain
\begin{equation}
\label{eq:full-range-dissipative-step}
\begin{aligned}
 \widetilde{\mathcal H}_j^{n+1}-\widetilde{\mathcal H}_j^n
 +\frac{k\bar\nu}{4}\|\nabla\bm w_j^{n+1}\|^2
 +k\varepsilon\sigma\mathcal F_j^{n+1}
 +\frac{k\gamma}{2}|q_j^{n+1}|^2
 \le k\widetilde C_F.
\end{aligned}
\end{equation}

\emph{Two-step contraction.}  Since $P_j\ge I$, we have
$\|V_j^{n+1}\|^2\le\mathcal F_j^{n+1}$, so
$\|V_j^{n+1}\|_{G_{\eta_j}}^2\le C_G\mathcal F_j^{n+1}$. By
\eqref{eq:G-norm-equivalence-preliminary},
$|Q_j^{n+1}|_G^2\le c_l^{-1}|Q_j^{n+1}|^2
=C_q(|q_j^{n+1}|^2+|q_j^n|^2)$. Hence
\begin{equation}
\label{eq:full-range-coercivity}
 \widetilde{\mathcal H}_j^{n+1}
 \le(C_G+k\varepsilon)\mathcal F_j^{n+1}
       +C_q(|q_j^{n+1}|^2+|q_j^n|^2).
\end{equation}
Let $n\ge2$. Adding \eqref{eq:full-range-dissipative-step} at the levels
$n$ and $n-1$, and dropping the gradient terms and
$k\varepsilon\sigma\mathcal F_j^n$, we obtain
\[
 \widetilde{\mathcal H}_j^{n+1}
 +k\varepsilon\sigma\mathcal F_j^{n+1}
 +\frac{k\gamma}{2}(|q_j^{n+1}|^2+|q_j^n|^2)
 \le\widetilde{\mathcal H}_j^{n-1}+2k\widetilde C_F.
\]
By the definition of $\widetilde\alpha_k$ in \eqref{eq:full-range-energy},
$\widetilde\alpha_k(C_G+k\varepsilon)\le k\varepsilon\sigma$ and
$\widetilde\alpha_kC_q\le k\gamma/2$. Thus, by
\eqref{eq:full-range-coercivity}, the dissipative terms on the left are
bounded below by $\widetilde\alpha_k\widetilde{\mathcal H}_j^{n+1}$, and
\begin{equation}
\label{eq:full-range-two-step}
 (1+\widetilde\alpha_k)\widetilde{\mathcal H}_j^{n+1}
 \le\widetilde{\mathcal H}_j^{n-1}+2k\widetilde C_F,
 \qquad n\ge2.
\end{equation}
For the second level, \eqref{eq:full-range-dissipative-step} with $n=1$
gives $\widetilde{\mathcal H}_j^2\le\widetilde{\mathcal H}_j^1+k\widetilde C_F$.

\emph{Iteration.} Set $\rho_k:=1+\widetilde\alpha_k$ and
$L_k:=2k\widetilde C_F/\widetilde\alpha_k$, so that
$\rho_kL_k=L_k+2k\widetilde C_F$. Subtracting this identity from
\eqref{eq:full-range-two-step} gives
\begin{equation}
\label{eq:full-range-shifted}
 \widetilde{\mathcal H}_j^{n+1}-L_k
 \le\rho_k^{-1}\bigl(\widetilde{\mathcal H}_j^{n-1}-L_k\bigr),
 \qquad n\ge2.
\end{equation}
For $n\ge1$, write
$n=m+2i$ with $m\in\{1,2\}$ and $i=\lfloor(n-1)/2\rfloor$. Iterating
\eqref{eq:full-range-shifted} $i$ times and using
$\widetilde{\mathcal H}_j^m\le\widetilde{\mathcal H}_j^1+k\widetilde C_F$
and $L_k\ge0$, we obtain
\begin{equation}
\label{eq:full-range-decay}
 \widetilde{\mathcal H}_j^n
 \le L_k+\rho_k^{-\lfloor(n-1)/2\rfloor}
       \bigl(\widetilde{\mathcal H}_j^1+k\widetilde C_F\bigr),
 \qquad n\ge1.
\end{equation}
Moreover, \eqref{eq:full-range-shifted} shows that
$\widetilde{\mathcal H}_j^{n+1}\le L_k$ if
$\widetilde{\mathcal H}_j^{n-1}\le L_k$, and
$\widetilde{\mathcal H}_j^{n+1}\le\widetilde{\mathcal H}_j^{n-1}$
otherwise. Thus
$\widetilde{\mathcal H}_j^{n+1}\le
\max\{\widetilde{\mathcal H}_j^{n-1},L_k\}$, and induction over each
parity gives
$\widetilde{\mathcal H}_j^n\le
\max\{\widetilde{\mathcal H}_j^1+k\widetilde C_F,L_k\}$ for all $n\ge1$.
This is \eqref{eq:asymmetric-uniform}. Letting $n\to\infty$ in
\eqref{eq:full-range-decay} gives \eqref{eq:asymmetric-limsup}.

\emph{Physical and ensemble bounds.} Since
$\|V_j^n\|_{G_{\eta_j}}^2\ge\lambda_{\min}(G_{\eta_j})\|\bm u_j^n\|^2$ and,
by \eqref{eq:G-norm-equivalence-preliminary},
$|Q_j^n|_G^2\ge c_u^{-1}|q_j^n|^2$, we have
\begin{equation}
\label{eq:full-range-physical-coercivity}
 c_G(\|\bm u_j^n\|^2+|q_j^n|^2)\le\mathcal E_j^n
 \le\widetilde{\mathcal H}_j^n.
\end{equation}

\noindent Together with \eqref{eq:asymmetric-uniform}, this gives
\eqref{eq:main-uniform-physical}. Averaging over $j$ and bounding each
$S_j(k)$ by $\max_jS_j(k)$ gives \eqref{eq:main-ensemble-bound}.
\end{proof}

The bound \eqref{eq:asymmetric-limsup} grows linearly in $k$ for
large $k$. A closer look at the large-step case shows that the solutions
are nevertheless absorbed by a ball whose radius is independent of $k$
and of the initial data.

\begin{corollary}
\label{cor:absorbing-radius}
Under the assumptions of Theorem~\ref{thm-long-time-stability},
let $\widetilde\alpha_1$ be the value of $\widetilde\alpha_k$ in
\eqref{eq:full-range-energy} at $k=1$, and define
\begin{equation}
\label{eq:full-range-radius}
 \widetilde R^2:=\max\left\{
 \frac{2\widetilde C_F}{c_G\widetilde\alpha_1},\quad
 \frac{2\widetilde C_F}{\varepsilon\widetilde\alpha_1}
 +\frac{2\widetilde C_F}{\gamma}
       \left(1+\frac{2}{\widetilde\alpha_1}\right)\right\}.
\end{equation}
Then, for every $k>0$,
\begin{equation}
\label{eq:full-range-physical-bound}
 \limsup_{n\to\infty}
 (\|\bm u_j^n\|^2+|q_j^n|^2)\le\widetilde R^2,
 \qquad
 \limsup_{n\to\infty}\frac1J\sum_{j=1}^J
 (\|\bm u_j^n\|^2+|q_j^n|^2)\le\widetilde R^2.
\end{equation}
Thus $\widetilde R$ is independent of $k$ and of the initial data;
for each $k$, every larger ball is absorbing, uniformly on sets with
bounded $\widetilde{\mathcal H}_j^1$.
\end{corollary}

\begin{proof}

\emph{A radius independent of $k$.} We use $\rho_k$, $L_k$, and the
estimates from the proof of Theorem~\ref{thm-long-time-stability}. The first term
$k\varepsilon\sigma/(C_G+k\varepsilon)$ in the minimum defining
$\widetilde\alpha_k$ is increasing in $k$, and its quotient by $k$
is decreasing. The second term is proportional to $k$. Hence
\[
 \widetilde\alpha_k\ge k\widetilde\alpha_1\quad(0<k\le1),
 \qquad
 \widetilde\alpha_k\ge\widetilde\alpha_1\quad(k\ge1).
\]
For $0<k\le1$, this gives $L_k\le2\widetilde C_F/\widetilde\alpha_1$.
Then \eqref{eq:asymmetric-limsup} and
\eqref{eq:full-range-physical-coercivity} yield
\[
 \limsup_{n\to\infty}(\|\bm u_j^n\|^2+|q_j^n|^2)
 \le\frac{L_k}{c_G}
 \le\frac{2\widetilde C_F}{c_G\widetilde\alpha_1},
\]
which is the first entry in \eqref{eq:full-range-radius}.
For $k\ge1$, we have $L_k\le2k\widetilde C_F/\widetilde\alpha_1$. Since
$\|\bm u_j^n\|^2\le\|V_j^n\|^2\le\mathcal F_j^n$, we have
$k\varepsilon\|\bm u_j^n\|^2\le\widetilde{\mathcal H}_j^n$, and therefore
\[
 \limsup_{n\to\infty}\|\bm u_j^n\|^2
 \le\frac{L_k}{k\varepsilon}
 \le\frac{2\widetilde C_F}{\varepsilon\widetilde\alpha_1}.
\]
For the scalar variable, we discard
$\mathcal E_j^{n+1}$, the second differences, and the gradient term in
\eqref{eq:full-range-basic-step}, and then divide by $k$:
\begin{equation}
\label{eq:full-range-q-bound}
 \frac\gamma2|q_j^{n+1}|^2
 \le\frac{\mathcal E_j^n}{k}+\widetilde C_F
 \le\frac{\widetilde{\mathcal H}_j^n}{k}+\widetilde C_F.
\end{equation}
Hence
\[
 \limsup_{n\to\infty}|q_j^n|^2
 \le\frac2\gamma\Bigl(\frac{L_k}{k}+\widetilde C_F\Bigr)
 \le\frac{2\widetilde C_F}{\gamma}
       \Bigl(1+\frac2{\widetilde\alpha_1}\Bigr).
\]
Adding the two bounds gives the second entry in
\eqref{eq:full-range-radius}. This proves the first bound in
\eqref{eq:full-range-physical-bound} for every $k>0$. The ensemble
bound follows because the $\limsup$ of an average is at most the average
of the $\limsup$s.

\emph{Absorption.} Fix $k>0$, $B>0$, and $R>\widetilde R$, and assume
$\widetilde{\mathcal H}_j^1\le B$. Set
$\tau_n:=\rho_k^{-\lfloor(n-2)/2\rfloor}(B+k\widetilde C_F)$ for
$n\ge2$. For $0<k\le1$, \eqref{eq:full-range-decay} and
\eqref{eq:full-range-physical-coercivity} give
\[
 \|\bm u_j^n\|^2+|q_j^n|^2
 \le\widetilde R^2+\frac{\tau_n}{c_G}.
\]
For $k\ge1$, \eqref{eq:full-range-decay} is applied at level $n$ to
$\|\bm u_j^n\|^2\le\widetilde{\mathcal H}_j^n/(k\varepsilon)$, and at level
$n-1$ in \eqref{eq:full-range-q-bound}. This gives
\[
 \|\bm u_j^n\|^2+|q_j^n|^2
 \le\widetilde R^2
 +\Bigl(\frac1{k\varepsilon}+\frac2{k\gamma}\Bigr)\tau_n.
\]

Since $\rho_k>1$, $\tau_n\to0$ as $n\to\infty$, and $\tau_n$ depends
only on $k$, $B$, and $n$. Hence there is $N=N(k,B,R)$ such that
$\|\bm u_j^n\|^2+|q_j^n|^2\le R^2$ for all $n\ge N$ and all such initial
data. The same $N$ works for every member, and hence for the ensemble
average, provided $B$ bounds $\widetilde{\mathcal H}_j^1$ for every $j$.
\end{proof}

\begin{remark}[Startup and spatial discretization]
\label{rem:asymmetric-startup}
The scheme can be initiated by the mr-SAV backward Euler step, which supplies a bounded
initial energy under \eqref{assumption:viscosity-fluctuation}
(Lemma~\ref{lem:startup-stability}).
Thus $\widetilde{\mathcal H}_j^1$ is bounded in terms of the prescribed data, uniformly for $k$ in bounded intervals.
The proof uses exact solenoidality through the identity
$b(\bm u_j^{*,n+1},\bm u_j^{*,n+1},\bm u_j^{*,n+1})=0$.
It extends to exactly divergence-free conforming finite elements. For a
pair with only weakly divergence-free velocities, such as Taylor--Hood,
the same proof applies if both discrete equations use the
skew-symmetric form
$b_{\rm sk}(\bm w,\bm u,\bm v):=(b(\bm w,\bm u,\bm v)-b(\bm w,\bm v,\bm u))/2$.
\end{remark}

\subsection{\texorpdfstring{Sharpness of the threshold $1/3$}{Sharpness of the threshold 1/3}}
\label{subsec:sharp-threshold}

We show that the positive threshold in
Theorem~\ref{thm-long-time-stability} cannot be increased. Let
$(\lambda,\bm w)$ be an eigenpair of the Stokes operator: $\lambda>0$,
$\bm w\in\bm V\setminus\{\bm0\}$, and for some $\pi\in Q$,
\begin{equation}
\label{eq:stokes-eigenpair}
  (\nabla\bm w,\nabla\bm v)-(\pi,\nabla\cdot\bm v)
  =\lambda(\bm w,\bm v),
  \qquad \bm v\in\bm X.
\end{equation}
To see how such a mode evolves, take $\bm F_j\equiv\bm0$, omit the
convective term, and insert $\bm u_j^m=z_0^m\bm w$ with $z_0\ne0$ into
\eqref{scheme:sav-bdf2-momentum}. Then
$\bm u_j^{*,n+1}=(2z_0^n-z_0^{n-1})\bm w$, and with the pressure
$p_j^{n+1}=\bigl(\bar\nu z_0^{n+1}+\nu_j'(2z_0^n-z_0^{n-1})\bigr)\pi$,
\eqref{eq:stokes-eigenpair} applied to the implicit term
$\bar\nu\bm u_j^{n+1}$ and to the explicit term
$\nu_j'\bm u_j^{*,n+1}$ reduces the momentum equation to
\[
 z_0^{n-1}\Bigl[\frac{3z_0^2-4z_0+1}{2k}
 +\lambda\bigl(\bar\nu z_0^2+\nu_j'(2z_0-1)\bigr)\Bigr](\bm w,\bm v)=0,
 \qquad \bm v\in\bm X.
\]
Multiplying the bracket by $2k$ and setting $s=k\lambda$ shows that this
holds if and only if $P_{k\lambda}(z_0)=0$, where $P_s$ is the
characteristic polynomial of the scheme on this mode,
\begin{equation}
\label{eq:mode-polynomial}
  P_s(z):=(3+2s\bar\nu)z^2+(4s\nu_j'-4)z+(1-2s\nu_j').
\end{equation}

The next lemma locates the roots of $P_s$ and shows that the
threshold is sharp: for $\nu_j'>\bar\nu/3$, exponentially growing
solutions exist for every $k>0$, and for $\nu_j'=\bar\nu/3$, the
absorbing radius cannot in general be chosen independently of $k$.

\begin{lemma}
\label{prop:sharp-threshold}
Fix a member $j$ with $\nu_j>0$.
\begin{enumerate}
\item If $\nu_j'\le\bar\nu/3$, then both roots of $P_s$ lie in the
open unit disk for every $s>0$. If $\nu_j'>\bar\nu/3$, then both roots
lie in the open unit disk if and only if $s<4/(3\nu_j'-\bar\nu)$, and
$P_s$ has a real root $z_0<-1$ for every $s>4/(3\nu_j'-\bar\nu)$.
\item Let $\bm F_j\equiv\bm0$, $\nu_j'>\bar\nu/3$, $k>0$, and $k\lambda>4/(3\nu_j'-\bar\nu)$,
and let $z_0<-1$ be a root of $P_{k\lambda}$ given by~(i). Then
$\bm u_j^n:=z_0^n\bm w$, with a suitable pressure, solves
\eqref{scheme:sav-bdf2-momentum}--\eqref{scheme:sav-bdf2-continuity}
without the convective term. In particular, for
$\nu_j'>\bar\nu/3$ and every $k>0$, the linearization of the scheme at
$(\bm u_j,q_j)=(\bm0,1)$ has exponentially growing solutions.
\item Let $\Omega=\{|x|<R\}$ be a disk, and let
$\bm w_m:=J_1(\sqrt{\lambda_m}\,r)\,\bm e_\theta$ in polar coordinates
$(r,\theta)$, with
$\lambda_m:=(j_{1,m}/R)^2$, where $J_1$ is the Bessel function of the
first kind of order one and $j_{1,m}$ its $m$th positive zero. If
$\bm F_j\equiv\bm0$ and $\nu_j'>\bar\nu/3$, then for every $k>0$ there
are $m$ and a root $z_0<-1$ of $P_{k\lambda_m}$ such that
$(\bm u_j^n,q_j^n):=(z_0^n\bm w_m,1)$, with a suitable pressure, solves
the full scheme \eqref{scheme:sav-bdf2-ensemble} with the starting levels
$(\bm w_m,1)$ and $(z_0\bm w_m,1)$. Thus
$\|\bm u_j^n\|=|z_0|^n\|\bm w_m\|\to\infty$, and no bound of the form of
Theorem~\ref{thm-long-time-stability} holds.
\item Let $\Omega$ be as in~(iii) and $\nu_j'=\bar\nu/3$. Then there
exist forces satisfying \eqref{assumption:discrete-force} with $M_F$
independent of $k$ such that, for every $k>0$, the full scheme
\eqref{scheme:sav-bdf2-ensemble} has a solution with
$\|\bm u_j^n\|=\kappa k$ for all $n\ge0$, where $\kappa>0$ is independent
of $k$ and $n$. In particular, these trajectories cannot be absorbed by
a ball whose radius is independent of $k$.
\end{enumerate}
\end{lemma}

\begin{proof}
(i) Write $P_s(z)=az^2+bz+c$ with $a>0$. By the Schur--Cohn conditions~\cite{Miller1971}, both roots lie in the open unit disk if and only if
$|c|<a$ and $a\pm b+c>0$. Here $a-c=2+2s\nu_j>0$, $a+b+c=2s\nu_j>0$,
and $a-b+c=8+2s(\bar\nu-3\nu_j')$. Since
$a+c=\frac12[(a+b+c)+(a-b+c)]$, both roots lie in the open unit disk if
and only if $a-b+c>0$. If $\nu_j'\le\bar\nu/3$, then $a-b+c\ge8$ for
every $s>0$. If $\nu_j'>\bar\nu/3$, then
$a-b+c>0$ if and only if $s<4/(3\nu_j'-\bar\nu)$. If
$s>4/(3\nu_j'-\bar\nu)$, then $P_s(-1)=a-b+c<0$, while
$P_s(z)\to+\infty$ as $z\to-\infty$; hence $P_s$ has a root $z_0<-1$.

(ii) Fix $n\ge1$, let $\bm u_j^m:=z_0^m\bm w$ for $m=n-1,n,n+1$, and
let $p_j^{n+1}\in Q$ be the pressure chosen before
\eqref{eq:mode-polynomial}. By the computation there, for every
$\bm v\in\bm X$ the left-hand side of \eqref{scheme:sav-bdf2-momentum}
without the convective term equals
$(2k)^{-1}z_0^{n-1}P_{k\lambda}(z_0)(\bm w,\bm v)=0$, and the right-hand
side vanishes since $\bm F_j\equiv\bm0$. Moreover,
\eqref{scheme:sav-bdf2-continuity} holds because $\bm w\in\bm V$. Since
$n\ge1$ is arbitrary, $(\bm u_j^n,p_j^n)$ solves the scheme for the
Stokes equations with the starting levels $\bm w$ and $z_0\bm w$, and
$\|\bm u_j^n\|=|z_0|^n\|\bm w\|\to\infty$.

Since the convective and trilinear terms are at least quadratic in the velocity,
the linearization at $(\bm0,1)$ consists of the scheme for the Stokes
equations and the decoupled recursion
$\delta_tq_j'^{\,n+1}+\gamma q_j'^{\,n+1}=0$ for the perturbation
$q_j'$ of $q_j$; hence $(z_0^n\bm w,0)$ solves it. Finally, since $\Omega$ is bounded, the Stokes operator has a compact
inverse, so its eigenvalues tend to infinity \cite{Temam1997}. Hence,
for every $k>0$, there is an eigenpair with
$k\lambda>4/(3\nu_j'-\bar\nu)$.

(iii) Write $\bm w_m=f_m(r)\,\bm e_\theta$ with
$f_m(r):=J_1(\sqrt{\lambda_m}\,r)$. These fields are smooth,
divergence-free, and vanish at $r=R$. Bessel's equation
$f_m''+r^{-1}f_m'-r^{-2}f_m=-\lambda_mf_m$ gives
$-\Delta\bm w_m=\lambda_m\bm w_m$, so $(\lambda_m,\bm w_m)$ satisfies
\eqref{eq:stokes-eigenpair} with $\pi=0$; see
\cite[Section~1]{Kelliher2010} for the Stokes spectrum of the disk.
Since $\lambda_m\to\infty$, for every $k>0$ some $m$ satisfies
$k\lambda_m>4/(3\nu_j'-\bar\nu)$; let $\bm w:=\bm w_m$ and let $z_0$ be as
in~(ii). Moreover, $(\bm w\cdot\nabla)\bm w=-r^{-1}f_m^2\,\bm e_r=\nabla\Pi$,
where $\Pi(r):=-\int_0^rf_m(r')^2r'^{-1}\,dr'$ minus its mean lies in $Q$.
Set
$z_*:=2z_0^n-z_0^{n-1}$, so that $\bm u_j^{*,n+1}=z_*\bm w$. The convective
term equals $-q_j^{n+1}z_*^2(\Pi,\nabla\cdot\bm v)$ and is absorbed by
taking $p_j^{n+1}:=-z_*^2\Pi$ as the pressure; the pressure of~(ii)
vanishes, since $\pi=0$. Moreover,
$b(\bm u_j^{*,n+1},\bm u_j^{*,n+1},\bm u_j^{n+1})=0$, so that
\eqref{scheme:sav-bdf2-q} reduces to
$\delta_tq_j^{n+1}+\gamma q_j^{n+1}=\gamma$, which $q_j\equiv1$
satisfies. Hence, by~(ii), $(z_0^n\bm w,1)$ solves the full scheme.

(iv) Let $\bm w:=\bm w_m$ for some $m$, fix $F_0>0$, and set
$\bm F_j^{n+1}:=(-1)^{n+1}F_0\bm w$, $c_n:=(-1)^nkF_0/4$, and
$(\bm u_j^n,q_j^n):=(c_n\bm w,1)$ for $n\ge0$. Then
$\|\bm F_j^{n+1}\|_{-1}=F_0\|\bm w\|_{-1}$ for all $n$, so
\eqref{assumption:discrete-force} holds with $M_F=F_0\|\bm w\|_{-1}$, and
$\|\bm u_j^n\|=\kappa k$ with $\kappa:=F_0\|\bm w\|/4$. It remains to
verify the scheme. We have
$\bm u_j^{*,n+1}=(2c_n-c_{n-1})\bm w=-3\bm u_j^{n+1}$, and the split
viscous terms cancel:
$\bar\nu(\nabla\bm u_j^{n+1},\nabla\bm v)
+\nu_j'(\nabla\bm u_j^{*,n+1},\nabla\bm v)
=(\bar\nu-3\nu_j')(\nabla\bm u_j^{n+1},\nabla\bm v)=0$. Moreover,
$\delta_t\bm u_j^{n+1}=(2k)^{-1}(3c_{n+1}-4c_n+c_{n-1})\bm w
=4c_{n+1}k^{-1}\bm w=\bm F_j^{n+1}$. With $\Pi$ from~(iii), the
convective term $9c_{n+1}^2b(\bm w,\bm w,\bm v)
=-9c_{n+1}^2(\Pi,\nabla\cdot\bm v)$ is absorbed by taking
$p_j^{n+1}:=-9c_{n+1}^2\Pi$ as the pressure, and
$b(\bm u_j^{*,n+1},\bm u_j^{*,n+1},\bm u_j^{n+1})=0$, so $q_j\equiv1$
solves \eqref{scheme:sav-bdf2-q}. For example, the nodal values
of $F_0\cos(\pi t)\bm w$ give this forcing when $k$ is an odd integer.
\end{proof}

\begin{remark}[Choice of the shared viscosity]
\label{rem:nubar-choice}
Let $\nu_{\max}:=\max_j\nu_j$, $\nu_{\min}:=\min_j\nu_j$, and
$\nu_{\rm mid}:=(\nu_{\max}+\nu_{\min})/2$.
\begin{enumerate}
\item \emph{Admissibility.} Condition
\eqref{assumption:viscosity-fluctuation} reads
$\bar\nu>\frac34\nu_{\max}$. For example, $\bar\nu=\nu_{\max}$ gives the
conclusions of Theorem~\ref{thm-long-time-stability} and
Corollary~\ref{cor:absorbing-radius} for any finite set
of positive viscosities. The threshold in
Lemma~\ref{prop:sharp-threshold} restricts $\bar\nu$, not the
spread of the samples. The ensemble mean satisfies this condition exactly when
$\nu_{\max}<\frac43\cdot\frac1J\sum_{j=1}^J\nu_j$, for two members
exactly when $\nu_{\max}<2\nu_{\min}$;
Section~\ref{sec:large-steps} shows that when it fails, the scheme
becomes unstable for large steps. A deterministic $\bar\nu$ fixed before
sampling, for instance from the support of the viscosity distribution,
avoids this check, keeps the members independent (see~(iv)), and lets one
factorization serve ensembles of all sizes
(Section~\ref{exp:sec:sampling}).
\item \emph{Cost of the spread.} Since
$p_\eta\ge(1-\eta)^{-2}$ on $(-\frac13,1)$ by
\eqref{eq:full-range-P-explicit}, and $1-\eta_j=\nu_j/\bar\nu$, we have
$\Lambda\ge(\bar\nu/\nu_{\min})^2>\frac9{16}(\nu_{\max}/\nu_{\min})^2$.
Hence the constants in \eqref{eq:full-range-constants} and the
radius \eqref{eq:full-range-radius} deteriorate as
$\nu_{\max}/\nu_{\min}$ grows. A wide spread also weakens the damping of
stiff modes. For the scheme applied to the Stokes equations, the roots
of \eqref{eq:mode-polynomial} tend to those of $z^2-2\eta_jz+\eta_j$ as
$k\lambda\to\infty$, so stiff modes are damped only by the factor
$\sqrt{1-\nu_j/\bar\nu}$ per step if $\nu_j<\bar\nu$, and by
$|\eta_j|+(\eta_j^2-\eta_j)^{1/2}$ if $\nu_j>\bar\nu$. Both factors tend
to $1$ as $\nu_j/\bar\nu$ approaches $0$ or $\frac43$, whereas unsplit
BDF2 damps stiff modes almost completely. This does not affect
Theorem~\ref{thm-long-time-stability}, but it matters for nonsmooth data
and transients.
\item \emph{A practical rule.} The midrange $\nu_{\rm mid}$
minimizes $\max_j|\nu_j'|$ and is admissible exactly when
$\nu_{\max}<2\nu_{\min}$. The choice
$\bar\nu:=\max\{\nu_{\rm mid},\,0.8\,\nu_{\max}\}$ is the midrange if
$\nu_{\max}\le\frac53\nu_{\min}$, and it keeps $\eta_j\ge-\frac14$. For
the samples \eqref{eq:viscosity-regimes}, $\nu_{\max}/\nu_{\min}=1.6$,
and the rule gives the midrange, which equals the mean used in
Section~\ref{subsec:convergence-tests}.
\item \emph{Groups.} For a wide spread, the samples can be partitioned
into groups with ratio at most $r_g$. Each group is advanced with its
own shared viscosity and time-independent coefficient matrix; the
number of matrices is about $\ln(\nu_{\max}/\nu_{\min})/\ln r_g$. Since the
groups do not interact, all results apply groupwise. The group mean is
admissible for any values in the group if $r_g<\frac43$. Centering each
group at its mean cancels the leading common component of the bias that
the splitting adds to noncentered ensemble statistics. In Monte Carlo
sampling, fixing the groups and their centers before sampling, for
example at the conditional mean of the input distribution, keeps the
samples independent. Spatially varying viscosities are discussed in
Remark~\ref{rem:random-viscosity-choice}.
\end{enumerate}
\end{remark}

\begin{remark}[Higher-order BDF methods]
\label{rem:higher-order-bdf}
Replacing BDF2 and the second-order extrapolation in
\eqref{scheme:sav-bdf2-ensemble} by BDF$\ell$ and the extrapolation of
order $\ell$, $3\le\ell\le5$, gives, with suitable starting levels, a
scheme of order $\ell$ with the same linear structure. The admissible
spread, however, shrinks. As in
Lemma~\ref{prop:sharp-threshold}, the characteristic roots of the Stokes
part on an eigenmode tend, as $k\lambda\to\infty$, to those of
$(1-\eta_j)z^\ell+\eta_j(z-1)^\ell$, and a root outside the closed unit
disk yields exponentially growing solutions for every $k>0$. Since
$|z|\le1$ if and only if $\operatorname{Re}(z/(z-1))\le\frac12$,
stability for every $k$ requires
\[
 -\frac1{2^\ell-1}\le\eta_j\le\frac1{1+(2\cos(\pi/\ell))^\ell},
\]
that is, $\nu_j/\bar\nu\in[\frac12,\frac87]$, $[\frac45,\frac{16}{15}]$,
and approximately $[0.917,1.032]$ for $\ell=3,4,5$. Unlike for BDF2,
this bounds the spread of the samples, $\nu_{\max}/\nu_{\min}\le\frac{16}7$,
$\frac43$, and about $1.125$, whatever $\bar\nu$ is, so wider ensembles
must be split into groups as in Remark~\ref{rem:nubar-choice}(iv).
We will investigate in future work whether these conditions are also
sufficient for an analogue of Theorem~\ref{thm-long-time-stability}.
\end{remark}

\subsection{\texorpdfstring{Spatially varying random viscosity}{Spatially varying random viscosity}}
\label{subsec:random-viscosity}

In uncertainty quantification the viscosity is often a random field,
for instance a truncated Karhunen--Lo\`eve expansion. We extend the scheme and its stability analysis to this case: for
general fields under a signed pointwise condition, and on the full range
for proportional fields, which share one spatial profile.

\paragraph{Setting}
Let $\nu_j\in L^\infty(\Omega)$ with $\operatorname*{ess\,inf}_\Omega\nu_j>0$,
$j=1,\ldots,J$. The viscous term is then written either as
$-\nabla\cdot(\nu_j\nabla\bm u_j)$, as in \eqref{eq:parameterized-nse},
or in the stress form $-\nabla\cdot(2\nu_j\bm D(\bm u_j))$ with
$\bm D(\bm u):=\tfrac12(\nabla\bm u+\nabla\bm u^\top)$, which is the
physically appropriate one for variable viscosity. With
\begin{equation}
\label{def:viscous-strain}
  \bm\varepsilon(\bm v):=\nabla\bm v
  \qquad\text{or}\qquad
  \bm\varepsilon(\bm v):=\sqrt2\,\bm D(\bm v),
\end{equation}
its weak form is $(\nu_j\bm\varepsilon(\bm u_j),\bm\varepsilon(\bm v))$
in both cases. Since
$(\nabla\bm u^\top,\nabla\bm v)=(\nabla\cdot\bm u,\nabla\cdot\bm v)$ on
$\bm X$ (two integrations by parts),
\begin{equation}
\label{eq:strain-identity}
  (\bm\varepsilon(\bm u),\bm\varepsilon(\bm v))
  =(\nabla\bm u,\nabla\bm v)+\kappa_\varepsilon(\nabla\cdot\bm u,\nabla\cdot\bm v),
  \qquad \bm u,\bm v\in\bm X,
\end{equation}
with $\kappa_\varepsilon=0$ for the first and $\kappa_\varepsilon=1$ for
the second form. The same identity gives
$\|\nabla\cdot\bm v\|\le\|\nabla\bm v\|$, so that
$\|\nabla\bm v\|\le\|\bm\varepsilon(\bm v)\|\le\sqrt2\,\|\nabla\bm v\|$ on
$\bm X$, and
\begin{equation}
\label{eq:strain-norm-V}
  \|\bm\varepsilon(\bm v)\|=\|\nabla\bm v\|,
  \qquad \bm v\in\bm V .
\end{equation}
In particular, $(\nu_j\bm\varepsilon(\bm u),\bm\varepsilon(\bm v))
=\nu_j(\nabla\bm u,\nabla\bm v)$ for constant $\nu_j$, $\bm u\in\bm V$,
and $\bm v\in\bm X$.

Let $\bar\nu\in L^\infty(\Omega)$ with $\operatorname*{ess\,inf}_\Omega\bar\nu>0$
be a shared field (choices are discussed in
Remark~\ref{rem:random-viscosity-choice}), and set
\[
\begin{gathered}
  \nu_j':=\nu_j-\bar\nu,\qquad
  (\nu_j')^{+}:=\max\{\nu_j',0\},\qquad
  (\nu_j')^{-}:=\max\{-\nu_j',0\},\\
  \rho:=\max_{1\le j\le J}|\nu_j'|,\qquad
  \rho_\infty:=\|\rho\|_{L^\infty(\Omega)};
\end{gathered}
\]
for constant viscosities and constant $\bar\nu$ these reduce to
\eqref{def:viscosity-splitting}. We consider general fields satisfying
the signed condition
\begin{equation}
\label{assumption:spatial-signed-fluctuation}
 d_{\nu,\pm}:=\operatorname*{ess\,inf}_{x\in\Omega}
 \min_{1\le j\le J}
 \Bigl(\bar\nu(x)-3(\nu_j')^{+}(x)-\frac32(\nu_j')^{-}(x)\Bigr)>0,
\end{equation}
which keeps $\nu_j(x)$ within $(\bar\nu(x)/3,\frac43\bar\nu(x))$, that
is, $-2\bar\nu(x)/3<\nu_j'(x)<\bar\nu(x)/3$, with a positive uniform
margin; and proportional fields
\begin{equation}
\label{assumption:proportional-viscosity}
 \nu_j(x)=a_j\chi(x),\qquad \bar\nu(x)=a_0\chi(x),
 \qquad a_0>0,\quad 0<a_j<\frac43a_0,
\end{equation}
with $\chi\in L^\infty(\Omega)$ and
$\chi_*:=\operatorname*{ess\,inf}_\Omega\chi>0$, which admit the full
range $0<\nu_j(x)<\frac43\bar\nu(x)$, that is,
$-\bar\nu(x)<\nu_j'(x)<\bar\nu(x)/3$. Since
$3(\nu_j')^{+}+\frac32(\nu_j')^{-}\le3|\nu_j'|$, \eqref{assumption:spatial-signed-fluctuation}
is implied by the pointwise symmetric condition
$\operatorname*{ess\,inf}_\Omega(\bar\nu-3\rho)>0$; for constant viscosities,
\eqref{assumption:proportional-viscosity} with $\chi\equiv1$ is
\eqref{assumption:viscosity-fluctuation}.

\paragraph{The scheme}
For $n\ge1$, find $(\bm u_j^{n+1},p_j^{n+1},q_j^{n+1})\in\bm X\times Q\times\mathbb R$
such that, for all $(\bm v,r)\in\bm X\times Q$,
\begin{subequations}
\label{scheme:random-viscosity}
\begin{align}
  &(\delta_t\bm u_j^{n+1},\bm v)
  +q_j^{n+1}b(\bm u_j^{*,n+1},\bm u_j^{*,n+1},\bm v)
  +(\bar\nu\,\bm\varepsilon(\bm u_j^{n+1}),\bm\varepsilon(\bm v))
  \notag\\
  &\qquad
  +(\nu_j'\,\bm\varepsilon(\bm u_j^{*,n+1}),\bm\varepsilon(\bm v))
  -(p_j^{n+1},\nabla\cdot\bm v)
  =\langle\bm F_j^{n+1},\bm v\rangle,
  \label{scheme:random-viscosity-momentum}\\
  &(\nabla\cdot\bm u_j^{n+1},r)=0,
  \label{scheme:random-viscosity-continuity}\\
  &\delta_tq_j^{n+1}+\gamma q_j^{n+1}
  -b(\bm u_j^{*,n+1},\bm u_j^{*,n+1},\bm u_j^{n+1})=\gamma,
  \label{scheme:random-viscosity-q}
\end{align}
\end{subequations}

with both starting levels in $\bm V\times\mathbb R$. By
\eqref{eq:strain-identity}, it coincides with
\eqref{scheme:sav-bdf2-ensemble} for constant viscosities and constant
$\bar\nu$. Since $\bar\nu$ depends neither on $j$ nor on $n$, the
implicit velocity--pressure operator, now associated with
\begin{equation}
\label{def:random-brinkman-form}
  a_k(\bm w,\bm v):=\frac{3}{2k}(\bm w,\bm v)
  +(\bar\nu\,\bm\varepsilon(\bm w),\bm\varepsilon(\bm v)),
\end{equation}
is again shared by all members and time levels.

For general fields we use the modified energy
$\widehat{\mathcal H}_j^n:=\mathcal G_j^n+\widehat\Phi_j^n$, with
$\mathcal G_j^n$ from \eqref{eq:basic-G-energy} and the correction term
$\widehat\Phi_j^n$, which compensates the explicit fluctuation term (see
\eqref{eq:spatial-signed-compensation}):
\begin{equation}
\label{eq:spatial-signed-compensator}
\begin{aligned}
\widehat\Phi_j^n:=&\frac{k}{2}\int_\Omega\Bigl[
 (\nu_j')^{+}\left(3|\bm\varepsilon(\bm u_j^n)|^2
                   +|\bm\varepsilon(\bm u_j^{n-1})|^2\right)\\
 &+(\nu_j')^{-}\left(|\bm\varepsilon(\bm u_j^n)|^2
          +|\bm\varepsilon(\bm u_j^n-\bm u_j^{n-1})|^2\right)\Bigr]\,dx.
\end{aligned}
\end{equation}

The next theorem extends unique solvability and the long-time
stability of Theorem~\ref{thm-long-time-stability} to spatially varying
viscosities, for general fields under the signed condition and for
proportional fields.

\begin{theorem}
\label{thm:random-viscosity-stability}
Let $\nu_j,\bar\nu\in L^\infty(\Omega)$ with
$\operatorname*{ess\,inf}_\Omega\bar\nu>0$, consider
\eqref{scheme:random-viscosity} with either viscosity form in
\eqref{def:viscous-strain}, assume \eqref{assumption:discrete-force},
and let both starting levels belong to $\bm V\times\mathbb R$. Then, for
every $k>0$, each step of \eqref{scheme:random-viscosity} is uniquely
solvable, $\bm u_j^{n+1}\in\bm V$, and the decoupled implementation of
Section~\ref{subsec:decoupling-solvability} applies with $a_k$ given by
\eqref{def:random-brinkman-form}; in particular,
$D_j^{n+1}\ge3/(2k)+\gamma$. Assume, moreover, either
\textup{(i)} the signed condition \eqref{assumption:spatial-signed-fluctuation}
for general fields, or
\textup{(ii)} the proportional condition \eqref{assumption:proportional-viscosity}.
Then, for every $k>0$, $\sup_{n\ge1}(\|\bm u_j^n\|^2+|q_j^n|^2)$ is
bounded in terms of $k$, the data, and the starting levels, and
\[
 \limsup_{n\to\infty}(\|\bm u_j^n\|^2+|q_j^n|^2)\le R^2,
\]
where $R$ depends only on $M_F$, $\gamma$, $\Omega$, and the
viscosities, and not on $k$ or the initial data.
\end{theorem}

\begin{proof}
Since $a_k(\bm v,\bm v)\ge\frac3{2k}\|\bm v\|^2+\nu_*\|\nabla\bm v\|^2$
with $\nu_*:=\operatorname*{ess\,inf}_\Omega\bar\nu$, the proof of
Proposition~\ref{prop:stepwise-solvability} applies verbatim.

\emph{Case~\textup{(i)}.} With $c_0$, $c_l$, $c_u$ from \eqref{eq:poincare} and
\eqref{eq:G-norm-equivalence-preliminary}, set
\begin{equation}
\label{eq:spatial-signed-constants}
\begin{gathered}
 \widehat C_F:=\frac{M_F^2}{2d_{\nu,\pm}}+\frac\gamma2,\qquad
 \widehat A_k:=\frac1{c_lc_0}+\frac{3k\rho_\infty}{2},\qquad
 \widehat\alpha_k:=\min\Bigl\{\frac{kd_{\nu,\pm}}{2\widehat A_k},
 \frac{k\gamma c_l}{2}\Bigr\},\\
 \widehat R^2:=\max\Bigl\{\frac{2c_u\widehat C_F}{\widehat\alpha_1},\;
 2\widehat C_F\Bigl(1+\frac2{\widehat\alpha_1}\Bigr)
 \max\Bigl\{\frac1{c_0d_{\nu,\pm}},\frac1\gamma\Bigr\}\Bigr\},
\end{gathered}
\end{equation}
where $\widehat\alpha_1$ is the value of $\widehat\alpha_k$ at $k=1$.

Let $x=\bm\varepsilon(\bm u_j^{n+1})$, $y=\bm\varepsilon(\bm u_j^n)$,
and $z=\bm\varepsilon(\bm u_j^{n-1})$, so that
$\bm\varepsilon(\bm u_j^{*,n+1})=2y-z$. Where $\nu_j'\ge0$, direct
expansion gives, with $c:=(\nu_j')^{+}$ and
$\psi_c^+(x,y):=\frac c2(3|x|^2+|y|^2)$,
$\bar\nu|x|^2+c(2y-z):x=(\bar\nu-3c)|x|^2+\psi_c^+(x,y)
-\psi_c^+(y,z)+c\bigl(|x+y|^2+\frac12|x-z|^2\bigr)$. Where $\nu_j'<0$,
\begin{equation}
\label{eq:spatial-negative-identity}
\begin{aligned}
 \bar\nu|x|^2-c(2y-z):x
 &=\left(\bar\nu-\frac32c\right)|x|^2
       +\psi_c^-(x,y)-\psi_c^-(y,z)+\frac c2|x-y+z|^2,\\
 \psi_c^-(x,y)&:=\frac c2\left(|x|^2+|x-y|^2\right),
\end{aligned}
\end{equation}
with $c:=(\nu_j')^{-}$. Dropping the nonnegative last terms, integrating, and
using $\bar\nu-3(\nu_j')^{+}-\frac32(\nu_j')^{-}\ge d_{\nu,\pm}$ and
$\|\bm\varepsilon(\bm v)\|\ge\|\nabla\bm v\|$ yields
\begin{equation}
\label{eq:spatial-signed-compensation}
\begin{aligned}
 &k(\bar\nu\bm\varepsilon(\bm u_j^{n+1}),\bm\varepsilon(\bm u_j^{n+1}))
 +k(\nu_j'\bm\varepsilon(\bm u_j^{*,n+1}),\bm\varepsilon(\bm u_j^{n+1}))\\
 &\qquad\ge\widehat\Phi_j^{n+1}-\widehat\Phi_j^n
       +kd_{\nu,\pm}\|\nabla\bm u_j^{n+1}\|^2.
\end{aligned}
\end{equation}
Moreover, $\psi_c^\pm(x,y)\le\frac32c|x|^2+c|y|^2$ and $(\nu_j')^{+}+(\nu_j')^{-}=|\nu_j'|\le\rho$, so
\eqref{eq:strain-norm-V} gives
$0\le\widehat\Phi_j^n\le\frac{3k\rho_\infty}{2}
(\|\nabla\bm u_j^n\|^2+\|\nabla\bm u_j^{n-1}\|^2)$ for solenoidal
velocities. Now test
\eqref{scheme:random-viscosity-momentum} with $k\bm u_j^{n+1}$, multiply
\eqref{scheme:random-viscosity-q} by $kq_j^{n+1}$, and add. The pressure
term vanishes, the convective terms cancel, and the $G$-identity
\eqref{eq:bdf2-G-identity-preliminary} for the velocity and the scalar
gives
\[
\begin{aligned}
 &\mathcal G_j^{n+1}-\mathcal G_j^n
 +\frac14\|\bm u_j^{n+1}-2\bm u_j^n+\bm u_j^{n-1}\|^2
 +\frac14|q_j^{n+1}-2q_j^n+q_j^{n-1}|^2
 +k\gamma|q_j^{n+1}|^2\\
 &\quad+k(\bar\nu\bm\varepsilon(\bm u_j^{n+1}),\bm\varepsilon(\bm u_j^{n+1}))
 +k(\nu_j'\bm\varepsilon(\bm u_j^{*,n+1}),\bm\varepsilon(\bm u_j^{n+1}))
 =k\langle\bm F_j^{n+1},\bm u_j^{n+1}\rangle+k\gamma q_j^{n+1}.
\end{aligned}
\]
By \eqref{eq:duality-estimate}, \eqref{assumption:discrete-force}, and
Young's inequality,
$\langle\bm F_j^{n+1},\bm u_j^{n+1}\rangle
\le\frac{d_{\nu,\pm}}{2}\|\nabla\bm u_j^{n+1}\|^2
+\frac{M_F^2}{2d_{\nu,\pm}}$ and
$\gamma q_j^{n+1}\le\frac\gamma2|q_j^{n+1}|^2+\frac\gamma2$. Together
with \eqref{eq:spatial-signed-compensation}, this gives, for every
$n\ge1$,
\begin{equation}
\label{eq:spatial-signed-one-step}
\begin{aligned}
 &\widehat{\mathcal H}_j^{n+1}-\widehat{\mathcal H}_j^n
 +\frac14\|\bm u_j^{n+1}-2\bm u_j^n+\bm u_j^{n-1}\|^2
 +\frac14|q_j^{n+1}-2q_j^n+q_j^{n-1}|^2\\
 &\qquad+\frac{kd_{\nu,\pm}}{2}\|\nabla\bm u_j^{n+1}\|^2
 +\frac{k\gamma}{2}|q_j^{n+1}|^2
 \le k\widehat C_F.
\end{aligned}
\end{equation}

By \eqref{eq:G-norm-equivalence-preliminary}, \eqref{eq:poincare}, and
the bound on $\widehat\Phi_j^{n+1}$,
\[
 \widehat{\mathcal H}_j^{n+1}
 \le\widehat A_k\bigl(\|\nabla\bm u_j^{n+1}\|^2+\|\nabla\bm u_j^n\|^2\bigr)
 +c_l^{-1}\bigl(|q_j^{n+1}|^2+|q_j^n|^2\bigr).
\]
Adding \eqref{eq:spatial-signed-one-step} at the levels $n$ and $n-1$ and
using $\widehat\alpha_k\widehat A_k\le kd_{\nu,\pm}/2$ and
$\widehat\alpha_kc_l^{-1}\le k\gamma/2$ gives
$(1+\widehat\alpha_k)\widehat{\mathcal H}_j^{n+1}
\le\widehat{\mathcal H}_j^{n-1}+2k\widehat C_F$ for $n\ge2$, and
\eqref{eq:spatial-signed-one-step} with $n=1$ gives
$\widehat{\mathcal H}_j^2\le\widehat{\mathcal H}_j^1+k\widehat C_F$. The
iteration in the proof of Theorem~\ref{thm-long-time-stability} then
gives
\[
 \sup_{n\ge1}\widehat{\mathcal H}_j^n
 \le\max\Bigl\{\widehat{\mathcal H}_j^1+k\widehat C_F,
 \frac{2k\widehat C_F}{\widehat\alpha_k}\Bigr\},
 \qquad
 \limsup_{n\to\infty}\widehat{\mathcal H}_j^n
 \le\frac{2k\widehat C_F}{\widehat\alpha_k},
\]
and \eqref{eq:G-norm-equivalence-preliminary} gives
$\|\bm u_j^n\|^2+|q_j^n|^2\le c_u\mathcal G_j^n\le c_u\widehat{\mathcal H}_j^n$,
which proves the uniform bound.
Since $\widehat A_k$ and $k/\widehat A_k$ increase with $k$,
$\widehat\alpha_k\ge k\widehat\alpha_1$ for $0<k\le1$ and
$\widehat\alpha_k\ge\widehat\alpha_1$ for $k\ge1$. For $0<k\le1$, this
gives the first entry of $\widehat R^2$. For $k\ge1$, dropping
$\widehat{\mathcal H}_j^{n+1}\ge0$ in \eqref{eq:spatial-signed-one-step}
gives $\frac{d_{\nu,\pm}}{2}\|\nabla\bm u_j^{n+1}\|^2
+\frac\gamma2|q_j^{n+1}|^2\le\widehat{\mathcal H}_j^n/k+\widehat C_F$,
and \eqref{eq:poincare} with
$\limsup_n\widehat{\mathcal H}_j^n/k\le2\widehat C_F/\widehat\alpha_1$
gives the second entry. Thus the radius bound holds with $R=\widehat R$.

\emph{Case~\textup{(ii)}.} Here $\nu_*=a_0\chi_*$, and
$\nu_j'=-\eta_j\bar\nu$ with $\eta_j:=1-a_j/a_0\in(-\frac13,1)$
independent of $x$. Hence
$\bm w_j^{n+1}=\bm u_j^{n+1}-\eta_j\bm u_j^{*,n+1}\in\bm V$ and
\[
\begin{aligned}
 &(\bar\nu\bm\varepsilon(\bm u_j^{n+1}),\bm\varepsilon(\bm w_j^{n+1}))
 +(\nu_j'\bm\varepsilon(\bm u_j^{*,n+1}),\bm\varepsilon(\bm w_j^{n+1}))\\
 &\qquad=\int_\Omega\bar\nu|\bm\varepsilon(\bm w_j^{n+1})|^2\,dx
 \ge\nu_*\|\nabla\bm w_j^{n+1}\|^2.
\end{aligned}
\]
All temporal identities, the convective cancellation, and the filter
estimate are unchanged, so the proofs of
Theorem~\ref{thm-long-time-stability} and
Corollary~\ref{cor:absorbing-radius} apply with the coercivity
constant $\nu_*$, that is, with $\eta_j$ in the matrices and
energies and with $\bar\nu$ replaced by $\nu_*$ in $\widetilde C_F$ and
$\varepsilon$ in \eqref{eq:full-range-constants}. In particular,
\eqref{eq:main-uniform-physical}--\eqref{eq:main-ensemble-bound} and
\eqref{eq:full-range-physical-bound} hold, with $R=\widetilde R$ from
\eqref{eq:full-range-radius} so modified.
\end{proof}

\begin{remark}[Spatially varying filter]
\label{rem:variable-filter-obstruction}
For general fields, the choice
$\eta_j(x)=-\nu_j'(x)/\bar\nu(x)$ does not support the same proof. For
bounded coefficients, $\eta_j\bm u_j^{*,n+1}$ need not belong to
$\bm X$; even for Lipschitz $\eta_j$, with $\bm s=\bm u_j^{*,n+1}\in\bm V$,
\[
 \nabla\cdot(\bm u_j^{n+1}-\eta_j\bm s)=-\nabla\eta_j\cdot\bm s,
 \qquad
 b(\bm s,\bm s,\eta_j\bm s)=-\frac12\int_\Omega(\bm s\cdot\nabla\eta_j)|\bm s|^2\,dx,
\]
so the pressure and convective cancellations are lost, and
$\bm\varepsilon(\eta_j\bm s)$ contains coefficient gradients. Whether the
full pointwise range holds for general fields remains open.
\end{remark}

\begin{remark}[Choice of the shared field]
\label{rem:random-viscosity-choice}
Let $\nu_{\min}:=\min_j\nu_j$ and $\nu_{\max}:=\max_j\nu_j$ pointwise.
For a shared value $c\in[\nu_{\min}(x),\nu_{\max}(x)]$, the margin in
\eqref{assumption:spatial-signed-fluctuation} is
$\min\{4c-3\nu_{\max}(x),\,(3\nu_{\min}(x)-c)/2\}$. It is maximal at
$c=(2\nu_{\max}+\nu_{\min})/3$, with value $(4\nu_{\min}-\nu_{\max})/3$,
so this choice of $\bar\nu(x)$ reduces
\eqref{assumption:spatial-signed-fluctuation} to
$\operatorname*{ess\,inf}_\Omega(4\nu_{\min}-\nu_{\max})>0$, a pointwise
contrast below $4$; the pointwise midrange gives only
$2\nu_{\min}-\nu_{\max}$, a contrast below $2$. For random fields
$\nu(x,\omega)=\nu_0(x)\bigl(1+Y(x,\omega)\bigr)$ with
$\operatorname*{ess\,inf}\nu_0>0$ and $|Y|\le\delta<\frac35$, the
deterministic choice $\bar\nu=(1+\delta/3)\nu_0$ satisfies
\eqref{assumption:spatial-signed-fluctuation} for every $J$, with
$d_{\nu,\pm}\ge(1-5\delta/3)\operatorname*{ess\,inf}\nu_0$, and keeps
Monte Carlo samples independent (cf.\
Remark~\ref{rem:nubar-choice}(iv)). Conversely, let $\delta>\frac35$ and
let the samples be independent. If, for some $\delta'\in(\frac35,\delta)$
and some $S\subset\Omega$ of positive measure, the events
$\{Y\ge\delta'\text{ a.e.\ on }S\}$ and
$\{Y\le-\delta'\text{ a.e.\ on }S\}$ have positive probability, then,
with probability tending to one as $J\to\infty$, both occur among the
samples, so that $4\nu_{\min}-\nu_{\max}\le(3-5\delta')\nu_0<0$ on $S$
and no shared field satisfies
\eqref{assumption:spatial-signed-fluctuation}. For a given
ensemble, \eqref{assumption:spatial-signed-fluctuation} can be checked a
posteriori; it is a sufficient condition for the energy argument, not an
instability threshold. In case~(ii) of
Theorem~\ref{thm:random-viscosity-stability}, $a_0$ is chosen as
$\bar\nu$ in Remark~\ref{rem:nubar-choice}.
\end{remark}

\section{Numerical experiments}
\label{numerics}
\subsection{Convergence tests}
\label{subsec:convergence-tests}

We test the
mr-SAV--\allowbreak BDF2 ensemble scheme by coupled space--time
refinements with $k=0.1h$ and $k=0.1h^2$ in three viscosity regimes,
and by a time-step refinement on a fixed fine mesh.

\paragraph{Setup}
Let $\Omega=(0,1)^2$, $T=0.5$, and, for each member,
\begin{equation}
\begin{gathered}
    \psi_j=e^{\nu_j}\cos(2t)\,x^2(1-x)^2y^2(1-y)^2,
    \qquad
    p_j=e^{\nu_j}\cos(2t)\left(xy-\tfrac14\right),
    \\
    \bm u_j=\left(\partial_y\psi_j,-\partial_x\psi_j\right)^\top.
\end{gathered}
\label{eq:manufactured-solution-convergence}
\end{equation}
Then $\nabla\cdot\bm u_j=0$, $\bm u_j|_{\partial\Omega}=\bm0$, and
$\int_\Omega p_j\,d\bm x=0$; the force
$\bm F_j:=\partial_t\bm u_j+(\bm u_j\cdot\nabla)\bm u_j-\nu_j\Delta\bm u_j+\nabla p_j$
makes $(\bm u_j,p_j)$ an exact solution, and the exact auxiliary
variable is $q_j\equiv1$. The force is sampled at the new time level,
$\bm F_j^{n+1}:=\bm F_j(t^{n+1})$ for $n\ge0$, so that
$\bm R_{F,j}^{n+1}=\bm0$ in \eqref{eq:error-residuals}; the cell
average \eqref{def:discrete-force-average} would introduce an $O(k)$
consistency error.

We use Taylor--Hood $\mathbb P_2/\mathbb P_1$ elements on uniform
meshes of right triangles with parallel diagonals. Since Taylor--Hood
velocities are only weakly divergence-free,
$b(\bm v,\bm v,\bm v)=-\frac12(\nabla\cdot\bm v,|\bm v|^2)$ need not
vanish. We therefore implement the convective term in the
skew-symmetric form $b_{\rm sk}$ of Remark~\ref{rem:asymmetric-startup}
in both \eqref{scheme:sav-bdf2-momentum} and \eqref{scheme:sav-bdf2-q},
and the computations lie within the setting of
Theorem~\ref{thm-long-time-stability}. We take the nodal interpolant of
$\bm u_j(0)$ as discrete initial velocity, $q_j^0=1$, and
$\gamma=10^3$; the first level is computed by the mr-SAV backward
Euler (mr-SAV--BE) step of Appendix~\ref{subsec:be-startup}, with
the same $k$ as the subsequent BDF2 steps.

With $\bm u_{j,h}^N$ and $p_{j,h}^N$ the computed solution at $t^N=T$,
the pressure being reconstructed by \eqref{eq:affine-decomposition}, we
measure
\begin{align}
    E_{\bm u,0}^{(j)}
    &:=\|\bm u_{j,h}^{N}-\bm u_j(T)\|_{L^2(\Omega)},
    \qquad
    E_{\bm u,1}^{(j)}
    :=\|\bm u_{j,h}^{N}-\bm u_j(T)\|_{H^1(\Omega)},
    \nonumber\\
    E_{p,0}^{(j)}
    &:=\inf_{c\in\mathbb R}
    \|p_{j,h}^{N}-p_j(T)-c\|_{L^2(\Omega)},
    \label{eq:samplewise-convergence-errors}
\end{align}
with the full $H^1$-norm in $E_{\bm u,1}^{(j)}$, and the ensemble
root-mean-square error
\begin{equation}
    \mathcal E_X
    :=
    \Bigl(\frac1J\sum_{j=1}^{J}\bigl(E_X^{(j)}\bigr)^2\Bigr)^{1/2}.
    \label{eq:rms-convergence-error}
\end{equation}
The error integrals use a collapsed Gauss rule exact for polynomials of
degree ten, with the exact solution evaluated at the quadrature points. All errors are absolute; observed orders are
$\log_2(E_h/E_{h/2})$, and analogously for $k$.

\paragraph{Coupled space--time refinement}
We take $J=10$ and
\begin{equation}
    \nu_j^{(m)}
    =
    \bigl(3.0+0.2(j-1)\bigr)\times10^{-m},
    \qquad
    j=1,\ldots,10,
    \qquad
    m=3,4,5,
    \label{eq:viscosity-regimes}
\end{equation}
referred to as Regimes~I, II, and~III. The shared viscosity is the
ensemble mean $\bar\nu^{(m)}=3.9\times10^{-m}$, which equals the
midrange and is the choice of Remark~\ref{rem:nubar-choice}(iii). In all
three regimes,
\begin{equation}
    \frac{\nu_j'}{\bar\nu}\in\Bigl[-\frac3{13},\frac3{13}\Bigr],
    \qquad
    \frac{10}{13}\le\frac{\nu_j}{\bar\nu}\le\frac{16}{13}<\frac43,
    \label{eq:viscosity-fluctuation-condition}
\end{equation}
so \eqref{assumption:viscosity-fluctuation} holds with the same relative
fluctuation, while the absolute viscosity decreases by two orders of
magnitude. We use $h=1/8,\ldots,1/128$ with $k=0.1h$, and
$h=1/4,\ldots,1/64$ with $k=0.1h^2$. The approximation orders of
Taylor--Hood elements and of BDF2 suggest
\begin{equation}
    E_{\bm u,0}^{(j)}
    \lesssim
    h^3+k^2,
    \qquad
    E_{\bm u,1}^{(j)}
    +
    E_{p,0}^{(j)}
    \lesssim
    h^2+k^2.
    \label{eq:fully-discrete-error-structure}
\end{equation}
These scales are heuristic: no fully discrete estimate is proved here,
and Theorem~\ref{thm:finite-time-temporal-error} and
\eqref{eq:error-pressure-estimate} control the nodal velocity
$L^2$-error but only an accumulated velocity $H^1$-error and a
time-integrated pressure error. The velocity $L^2$ scale is
$h^3+0.01h^2$ for $k=0.1h$, so the $k^2$ term eventually gives order
two, and $h^3+0.01h^4$ for $k=0.1h^2$, so the spatial order three should
persist. Tables~\ref{tab:conv-rms-kh} and~\ref{tab:conv-rms-kh2} report
the ensemble RMS errors.

\begin{table}[!htbp]
\centering
\captionsetup{font=footnotesize,labelfont=footnotesize}
\caption{Ensemble RMS errors \eqref{eq:rms-convergence-error} and observed orders for $k=0.1h$, with $J=10$, $T=0.5$, $\gamma=10^3$, and the viscosities \eqref{eq:viscosity-regimes}.}
\label{tab:conv-rms-kh}
\scriptsize
\setlength{\tabcolsep}{4pt}
\renewcommand{\arraystretch}{1.15}
\begin{tabular}{ccccccc}
\toprule
\multicolumn{7}{c}{Regime~I: $\bar\nu=3.9\times10^{-3}$} \\
\midrule
$h$ & $\mathcal E_{\bm u,0}$ & Order & $\mathcal E_{\bm u,1}$ & Order & $\mathcal E_{p,0}$ & Order \\
\midrule
$1/8$ & 2.36994E-05 & -- & 1.38487E-03 & -- & 5.47073E-04 & -- \\
$1/16$ & 3.00217E-06 & 2.98 & 3.54022E-04 & 1.97 & 1.36768E-04 & 2.00 \\
$1/32$ & 4.00996E-07 & 2.90 & 8.91168E-05 & 1.99 & 3.41920E-05 & 2.00 \\
$1/64$ & 6.11014E-08 & 2.71 & 2.23213E-05 & 2.00 & 8.54799E-06 & 2.00 \\
$1/128$ & 1.15903E-08 & 2.40 & 5.58309E-06 & 2.00 & 2.13852E-06 & 2.00 \\
\midrule
\multicolumn{7}{c}{Regime~II: $\bar\nu=3.9\times10^{-4}$} \\
\midrule
$h$ & $\mathcal E_{\bm u,0}$ & Order & $\mathcal E_{\bm u,1}$ & Order & $\mathcal E_{p,0}$ & Order \\
\midrule
$1/8$ & 2.46536E-05 & -- & 1.41975E-03 & -- & 5.45156E-04 & -- \\
$1/16$ & 3.02688E-06 & 3.03 & 3.54480E-04 & 2.00 & 1.36289E-04 & 2.00 \\
$1/32$ & 4.07303E-07 & 2.89 & 8.88481E-05 & 2.00 & 3.40722E-05 & 2.00 \\
$1/64$ & 6.40443E-08 & 2.67 & 2.22452E-05 & 2.00 & 8.51809E-06 & 2.00 \\
$1/128$ & 1.25774E-08 & 2.35 & 5.56374E-06 & 2.00 & 2.13072E-06 & 2.00 \\
\midrule
\multicolumn{7}{c}{Regime~III: $\bar\nu=3.9\times10^{-5}$} \\
\midrule
$h$ & $\mathcal E_{\bm u,0}$ & Order & $\mathcal E_{\bm u,1}$ & Order & $\mathcal E_{p,0}$ & Order \\
\midrule
$1/8$ & 2.56987E-05 & -- & 1.46736E-03 & -- & 5.44966E-04 & -- \\
$1/16$ & 3.18492E-06 & 3.01 & 3.70817E-04 & 1.98 & 1.36241E-04 & 2.00 \\
$1/32$ & 4.27040E-07 & 2.90 & 9.39974E-05 & 1.98 & 3.40603E-05 & 2.00 \\
$1/64$ & 6.49088E-08 & 2.72 & 2.26210E-05 & 2.05 & 8.51506E-06 & 2.00 \\
$1/128$ & 1.26890E-08 & 2.35 & 5.57416E-06 & 2.02 & 2.12432E-06 & 2.00 \\
\bottomrule
\end{tabular}
\end{table}

\begin{table}[!htbp]
\centering
\captionsetup{font=footnotesize,labelfont=footnotesize}
\caption{As Table~\ref{tab:conv-rms-kh}, for $k=0.1h^2$.}
\label{tab:conv-rms-kh2}
\scriptsize
\setlength{\tabcolsep}{4pt}
\renewcommand{\arraystretch}{1.15}
\begin{tabular}{ccccccc}
\toprule
\multicolumn{7}{c}{Regime~I: $\bar\nu=3.9\times10^{-3}$} \\
\midrule
$h$ & $\mathcal E_{\bm u,0}$ & Order & $\mathcal E_{\bm u,1}$ & Order & $\mathcal E_{p,0}$ & Order \\
\midrule
$1/4$ & 1.84796E-04 & -- & 5.24129E-03 & -- & 2.18835E-03 & -- \\
$1/8$ & 2.30795E-05 & 3.00 & 1.38494E-03 & 1.92 & 5.47073E-04 & 2.00 \\
$1/16$ & 2.86786E-06 & 3.01 & 3.54003E-04 & 1.97 & 1.36768E-04 & 2.00 \\
$1/32$ & 3.59006E-07 & 3.00 & 8.91097E-05 & 1.99 & 3.41921E-05 & 2.00 \\
$1/64$ & 4.49177E-08 & 3.00 & 2.23194E-05 & 2.00 & 8.54804E-06 & 2.00 \\
\midrule
\multicolumn{7}{c}{Regime~II: $\bar\nu=3.9\times10^{-4}$} \\
\midrule
$h$ & $\mathcal E_{\bm u,0}$ & Order & $\mathcal E_{\bm u,1}$ & Order & $\mathcal E_{p,0}$ & Order \\
\midrule
$1/4$ & 1.95713E-04 & -- & 5.48191E-03 & -- & 2.18069E-03 & -- \\
$1/8$ & 2.40406E-05 & 3.03 & 1.41989E-03 & 1.95 & 5.45156E-04 & 2.00 \\
$1/16$ & 2.87867E-06 & 3.06 & 3.54460E-04 & 2.00 & 1.36289E-04 & 2.00 \\
$1/32$ & 3.57768E-07 & 3.01 & 8.88394E-05 & 2.00 & 3.40722E-05 & 2.00 \\
$1/64$ & 4.47536E-08 & 3.00 & 2.22428E-05 & 2.00 & 8.51808E-06 & 2.00 \\
\midrule
\multicolumn{7}{c}{Regime~III: $\bar\nu=3.9\times10^{-5}$} \\
\midrule
$h$ & $\mathcal E_{\bm u,0}$ & Order & $\mathcal E_{\bm u,1}$ & Order & $\mathcal E_{p,0}$ & Order \\
\midrule
$1/4$ & 1.98513E-04 & -- & 5.54368E-03 & -- & 2.17992E-03 & -- \\
$1/8$ & 2.51097E-05 & 2.98 & 1.46747E-03 & 1.92 & 5.44966E-04 & 2.00 \\
$1/16$ & 3.04295E-06 & 3.04 & 3.70797E-04 & 1.98 & 1.36241E-04 & 2.00 \\
$1/32$ & 3.79218E-07 & 3.00 & 9.39890E-05 & 1.98 & 3.40603E-05 & 2.00 \\
$1/64$ & 4.55140E-08 & 3.06 & 2.26186E-05 & 2.05 & 8.51506E-06 & 2.00 \\
\bottomrule
\end{tabular}
\end{table}

The three regimes behave alike, with no deterioration as the viscosity
decreases. The velocity $H^1$- and pressure errors converge with order
two. For $k=0.1h^2$ the velocity $L^2$-errors converge with order
three on all meshes; for $k=0.1h$ the order decreases from three to
about $2.4$ at the finest level, where the $O(k^2)$ term in
\eqref{eq:fully-discrete-error-structure} becomes visible. The errors
of the endpoint members $j=1,10$ differ from the RMS errors by less than
$1.7\%$, and their orders by at most $0.02$.

The pressure errors for the two time-step laws differ by less than a
relative $6\times10^{-6}$. For $1/64\le h\le1/8$, they agree to within a
relative $4\times10^{-6}$ with the RMS over $j$ of
$e^{\nu_j}\cos(1)\,h^2/\sqrt{240}$, the $L^2$-norm modulo constants of
the $\mathbb P_1$ nodal interpolation error of $p_j(T)$ on these meshes.
On these meshes, the computed pressure differs from the nodal
interpolant of $p_j(T)$, modulo constants, by at most $9\times10^{-3}$
times $E_{p,0}^{(j)}$, for all members. The pressure error is thus
dominated by interpolation; only for $h=1/128$ and $k=0.1h$ does a
further component appear, of up to $0.15\,E_{p,0}^{(j)}$, which changes
$E_{p,0}^{(j)}$ by at most $0.2\%$. The tables do not isolate the
temporal pressure error.

\paragraph{Fixed-mesh temporal refinement}
For Regime~III we fix $h=1/128$ and halve $k$ from $1/10$ to $1/160$.
This study prescribes the nodal interpolants of the exact solution at
$t^0$ and $t^1$, so the startup step is not involved. Since the errors
remain fully discrete, the temporal rate is visible only while the
temporal error exceeds the fixed-mesh spatial error; we therefore use
the velocity $L^2$-error.

\begin{table}[!htbp]
\centering
\captionsetup{
    font=footnotesize,
    labelfont=footnotesize,
    position=top,
    skip=6pt
}
\caption{Independent temporal convergence in the velocity $L^2$-norm
for Regime III, with
$\nu_j=(3.0+0.2(j-1))\times10^{-5}$, $j=1,\ldots,10$. The spatial
mesh is fixed at $h=1/128$, with $J=10$, $T=0.5$, and
$\gamma=10^3$.
The quantities $E_{\bm u,0}^{(1)}$ and $E_{\bm u,0}^{(10)}$ are the
absolute errors for the endpoint realizations, while
$\mathcal E_{\bm u,0}$ is the ensemble root-mean-square error.}
\label{tab:independent-temporal-convergence}
\scriptsize
\setlength{\tabcolsep}{5pt}
\renewcommand{\arraystretch}{1.2}
\begin{tabular}{ccccccc}
\hline
$k$
& $E_{\bm u,0}^{(1)}$ & Order
& $E_{\bm u,0}^{(10)}$ & Order
& $\mathcal E_{\bm u,0}$ & Order \\
\hline
$1/10$ & 3.52468E-05 & -- & 3.51512E-05 & -- & 3.51990E-05 & -- \\
$1/20$ & 1.03234E-05 & 1.77 & 1.02945E-05 & 1.77 & 1.03090E-05 & 1.77 \\
$1/40$ & 2.77988E-06 & 1.89 & 2.77205E-06 & 1.89 & 2.77597E-06 & 1.89 \\
$1/80$ & 7.20268E-07 & 1.95 & 7.18237E-07 & 1.95 & 7.19252E-07 & 1.95 \\
$1/160$ & 1.83287E-07 & 1.97 & 1.82770E-07 & 1.97 & 1.83028E-07 & 1.97 \\
\hline
\end{tabular}
\end{table}

The rates in Table~\ref{tab:independent-temporal-convergence} approach
two, with the same rates for both endpoint members and the RMS error.
Extrapolating the third-order behavior in Table~\ref{tab:conv-rms-kh2}
from $h=1/64$ gives a spatial velocity $L^2$-error of about
$6\times10^{-9}$ at $h=1/128$, far below the errors in
Table~\ref{tab:independent-temporal-convergence}, which are thus
dominated by the temporal error. The velocity $H^1$- and pressure errors
flatten earlier (data not tabulated), consistent with spatial
saturation, and are not used to infer temporal rates.

For $h=1/128$ and $k=0.1h$, the terminal auxiliary variables satisfy
$|q_j^N-1|\le2^{-52}\approx2.2\times10^{-16}$, the double-precision
machine epsilon, for all members and in all three regimes. This is not
a stringent test.
Since $\bm u_j=e^{\nu_j}\cos(2t)\,\bm w$ with $\bm w\in\bm V$
independent of $t$, the trilinear term of \eqref{scheme:sav-bdf2-q}
vanishes on the exact nodal values, which therefore satisfy
\eqref{scheme:sav-bdf2-q} with $q_j\equiv1$ without residual; hence
$q_j^n-1$ is driven only by the velocity discretization error, through
the exact recursion
\[
 q_j^{n+1}-1=
 \frac{4(q_j^n-1)-(q_j^{n-1}-1)
 +2k\,b(\bm u_j^{*,n+1},\bm u_j^{*,n+1},\bm u_j^{n+1})}
 {3+2k\gamma},
\]
in which the trilinear term has weight
$(3/(2k)+\gamma)^{-1}\le\gamma^{-1}=10^{-3}$. A stringent test requires a
non-separable solution and a moderate $\gamma$.

In summary, the coupled refinements show the expected fully discrete
rates in all three regimes, including the spatial velocity $L^2$-order
three for $k=0.1h^2$, and the fixed-mesh refinement shows second-order
temporal accuracy of the velocity in $L^2$. The test is
mild: the members differ only through the factor $e^{\nu_j}$, by less
than $0.2\%$ (as $|\nu_i-\nu_l|\le1.8\times10^{-3}$), and the viscous
terms ($\nu_j\le4.8\times10^{-3}$) are small, so the explicit treatment
of the fluctuation is not challenged. Robustness and accuracy for
under-resolved high-Reynolds-number flows, the pressure estimate
\eqref{eq:error-pressure-estimate}, and second-order convergence of
$q_j$ are not tested here; large time steps and the sharpness of
\eqref{assumption:viscosity-fluctuation} are tested in
Section~\ref{sec:large-steps}.

\subsection{Square-cavity flow with random initial data: statistics under joint refinement}
\label{exp:sec:joint}

We consider unforced flow on $\Omega=(0,1)^2$ up to $T=1$, with homogeneous
Dirichlet velocity conditions on all four walls. We fix the shared
viscosity $\bar\nu=1/3200$ before sampling and give each realization a
viscosity that is constant in space and time, with independent and
identically distributed (i.i.d.) $Z_j$,
\begin{equation}
 \label{exp:eq:joint-viscosity}
 \nu_j=\bar\nu(1+0.25Z_j),
 \qquad Z_j\overset{\mathrm{iid}}\sim\mathcal U[-1,1].
\end{equation}
Thus $\nu_j/\bar\nu\in[0.75,1.25]$ and $|\nu_j'|/\bar\nu\le\frac14$,
three quarters of the upper threshold $\frac13$, so
\eqref{assumption:viscosity-fluctuation} holds a priori, with margin
$\frac1{12}$ from $\frac43$. Since $\bar\nu$ does not depend on the
samples, the members are independent, and at each level all of them
share one operator.

The initial velocity is prescribed by
\begin{align}
 \psi_j(x,y)&=64x^2(1-x)^2y^2(1-y)^2
 \left[1+0.1\sum_{(a,b)\in\mathcal K}
 Y_{j,ab}\sin(a\pi x)\sin(b\pi y)\right],
 \label{exp:eq:streamfunction}\\
 \boldsymbol u_j^0&=(\partial_y\psi_j,-\partial_x\psi_j),
 \label{exp:eq:initial-velocity}
\end{align}
where $\mathcal K=\{(1,1),(2,1),(1,2),(2,2)\}$ and all $Y_{j,ab}$ are
independent $\mathcal U[-1,1]$ variables, independent also of the viscosities.
These initial data are smooth, divergence-free, and zero on $\partial\Omega$;
the discrete initial velocities are their mixed $L^2$ projections onto the
discretely divergence-free velocity space. The largest computed
ratio $|\nu_j/\bar\nu-1|$ is $0.2491$.

\emph{Experiment design.}
We use the mr-SAV--BDF2 ensemble scheme with $\mathbb P_2/\mathbb P_1$
Taylor--Hood elements, $\gamma=10^3$, and zero grad-div coefficient.
The convective form is the skew-symmetric one of
Section~\ref{subsec:convergence-tests}. The first level
$(\bm u_j^1,p_j^1,q_j^1)$ is computed from $(\bm u_j^0,q_j^0)$, with
$q_j^0=1$, by the mr-SAV--BE step of
Appendix~\ref{subsec:be-startup} with the same $k$, and the boundary
data are imposed at each new time level.
The four nested meshes in Table~\ref{exp:tab:refinement} are formed by
splitting each cell of a uniform $N\times N$ square grid along a consistent
diagonal, giving $2N^2$ triangles. The mesh spacing and time step are halved
and the ensemble size doubled between consecutive levels.
A pool of 1024 independent input vectors is generated once; each level
uses its first $J$ members, preserving the random parameters of shared
indices. Each level is solved with its own initial projection, and
all levels use the same $\bar\nu$. Thus random inputs are nested, whereas discrete trajectories
depend on the level. For accuracy, the time steps in
Table~\ref{exp:tab:refinement} are chosen to resolve the convective time
scale.

\begin{table}[!htbp]
 \centering
\captionsetup{font=footnotesize,labelfont=footnotesize}
 \caption{Joint refinement for the homogeneous-boundary experiment, $T=1$.}
 \label{exp:tab:refinement}
 \scriptsize
\setlength{\tabcolsep}{3.5pt}
\renewcommand{\arraystretch}{1.2}
 \begin{tabular}{ccccc}
 \toprule
 $N$ & Triangles & Time steps $N_t$ & $k$ & Samples $J$\\
 \midrule
 8 & 128 & 200 & $1/200$ & 64\\
 16 & 512 & 400 & $1/400$ & 128\\
 32 & 2048 & 800 & $1/800$ & 256\\
 64 & 8192 & 1600 & $1/1600$ & 512\\
 \bottomrule
 \end{tabular}
\end{table}

\paragraph{Moment-based cross-level diagnostics}
As in the moment diagnostics of \cite{Bansal2021}, we compute the
sample mean and componentwise sample variance of the terminal
velocities $\boldsymbol u_{i,j}$ of the members $j=1,\ldots,J_i$ at
refinement level $i$:
\begin{align}
 \mathbb E_{J_i}[\boldsymbol u_i]
 &:=\frac1{J_i}\sum_{j=1}^{J_i}\boldsymbol u_{i,j},
 \label{exp:eq:moment-mean}\\
 \operatorname{Var}_{J_i}[\boldsymbol u_i]
 &:=\frac{J_i}{J_i-1}
 \left(\mathbb E_{J_i}[\boldsymbol u_i^2]
       -\bigl(\mathbb E_{J_i}[\boldsymbol u_i]\bigr)^2\right),
 \label{exp:eq:moment-variance}
\end{align}
where squares of velocity vectors are taken componentwise. The
postprocessing evaluates the algebraically equivalent centered sum with
denominator $J_i-1$. Since $\bar\nu$ is fixed before sampling, the
members are i.i.d., and \eqref{exp:eq:moment-variance} is an unbiased
estimator of the componentwise variance of the computed terminal
velocity (cf.\ Remark~\ref{rem:nubar-choice}(iv)).

After the statistics on adjacent levels are represented on the common
finer mesh, we define the $L^2(\Omega)$ Cauchy differences at $T=1$ by
\begin{align}
 \varepsilon_{\mathrm{mean}}^{(i)}
 &:=\left\|\mathbb E_{J_i}[\boldsymbol u_i]
     -\mathbb E_{J_{i+1}}[\boldsymbol u_{i+1}]\right\|_{L^2(\Omega)},
 \label{exp:eq:moment-mean-cauchy}\\
 \varepsilon_{\mathrm{var}}^{(i)}
 &:=\left\|\operatorname{Var}_{J_i}[\boldsymbol u_i]
     -\operatorname{Var}_{J_{i+1}}[\boldsymbol u_{i+1}]\right\|_{L^2(\Omega)}.
 \label{exp:eq:moment-var-cauchy}
\end{align}
Both velocity components enter these norms. The $\mathbb P_2$ mean and
$\mathbb P_4$ variance fields are evaluated on the finer nested triangulation; their
squared differences are integrated there exactly, with a Gauss rule of
degree eight.
Because the mesh width, time step, and ensemble size vary simultaneously,
the differences are composite cross-level diagnostics, not isolated
spatial, temporal, or sampling errors.

The one-point and two-point empirical velocity distributions at
points $\bm x,\bm y\in\Omega$ are
\begin{equation}
 \label{exp:eq:empirical-laws}
 \mu_{i,\bm x}^{(1)}=\frac1{J_i}\sum_{j=1}^{J_i}
    \delta_{\boldsymbol u_{i,j}(\bm x)},
 \qquad
 \mu_{i,\bm x,\bm y}^{(2)}=\frac1{J_i}\sum_{j=1}^{J_i}
    \delta_{(\boldsymbol u_{i,j}(\bm x),\boldsymbol u_{i,j}(\bm y))}.
\end{equation}
Both velocities in a two-point atom belong to the same realization.
With $W_1$ denoting the $1$-Wasserstein distance with Euclidean cost,
we measure
\begin{align}
 \mathcal W_1^{(1)}(i)&=\int_{\Omega}
 W_1\bigl(\mu_{i,\bm x}^{(1)},\mu_{i+1,\bm x}^{(1)}\bigr)\,\mathrm d\bm x,
 \label{exp:eq:spatial-w-one}\\
 \mathcal W_1^{(2)}(i)&=\int_{\Omega^2}
 W_1\bigl(\mu_{i,\bm x,\bm y}^{(2)},\mu_{i+1,\bm x,\bm y}^{(2)}\bigr)\,\mathrm d\bm x\,\mathrm d\bm y.
 \label{exp:eq:spatial-w-two}
\end{align}
The state spaces are $\mathbb R^2$ and $\mathbb R^4$, respectively.
Replacing $\boldsymbol u$ by $|\boldsymbol u|$ gives the speed diagnostics.
The integrals are approximated with 2048 scrambled Sobol points in
$\Omega$ or point pairs in $\Omega^2$, at which both levels are evaluated.
Each transport problem uses all samples; since $J_{i+1}=2J_i$, splitting
each coarse-level atom into two atoms of half weight makes the exact
transport problem a minimum-cost assignment problem.

\begin{figure}[!ht]
 \centering
 \includegraphics[width=0.75\textwidth]{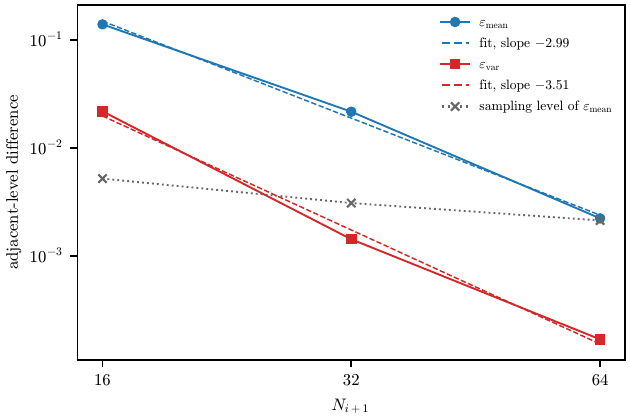}
 \captionsetup{font=footnotesize,labelfont=footnotesize}
 \caption{Adjacent-level $L^2(\Omega)$ differences of the ensemble mean and
 componentwise sample variance at $T=1$, computed between the refinement
 levels in Table~\ref{exp:tab:refinement}. The variance uses $J_i-1$
 normalization; fits use all three comparisons. The dotted
 line is the sampling level of $\varepsilon_{\mathrm{mean}}$ discussed in the
 text.}
 \label{exp:fig:moments}
\end{figure}

\emph{Results.}
Table~\ref{exp:tab:joint-errors} reports the adjacent-level
differences, and Figure~\ref{exp:fig:moments} shows the mean and
variance differences. Both decrease rapidly: from $8\to16$ to $16\to32$
by factors of $6.5$ and $15$, and from $16\to32$ to $32\to64$ by factors
of $9.7$ and $8.5$. Least-squares fits $\varepsilon\approx C N_{i+1}^{-p}$
over the three comparisons give exponents $2.99$ and $3.51$, respectively
(Table~\ref{exp:tab:joint-orders}). These values are not asymptotic
rates. Since consecutive levels share their first $J_i$ inputs, the mean
difference contains the sampling term
$\mathbb E_{J_i}[\bm u_{i+1}]-\mathbb E_{J_{i+1}}[\bm u_{i+1}]$, whose
expected squared $L^2(\Omega)$ norm is $(2J_i)^{-1}$ times the integral
over $\Omega$ of the sum of the componentwise variances of
$\bm u_{i+1}$, so that its size decreases only like $J_i^{-1/2}$. For the finest
comparison, its square root is $2.1\times10^{-3}$, close to
$\varepsilon_{\mathrm{mean}}=2.2\times10^{-3}$ (dotted line in
Figure~\ref{exp:fig:moments}): the discretization differences have
fallen below the sampling level, and on finer levels the differences
are expected to decrease only like $J^{-1/2}\propto N^{-1/2}$.

\begin{table}[!ht]
 \centering
\captionsetup{font=footnotesize,labelfont=footnotesize}
 \caption{Adjacent-level statistical differences at $T=1$.
 The distribution columns refer to velocity vectors.}
 \label{exp:tab:joint-errors}
 \scriptsize
\setlength{\tabcolsep}{3.5pt}
\renewcommand{\arraystretch}{1.2}
 \begin{tabular}{ccccc}
 \toprule
 $N_i\to N_{i+1}$ & $\varepsilon_{\mathrm{mean}}$ & $\varepsilon_{\mathrm{var}}$
 & $\mathcal W_1^{(1)}$ & $\mathcal W_1^{(2)}$\\
 \midrule
 $8\to16$ & 1.3973E-01 & 2.1833E-02 & 1.4258E-01 & 2.2422E-01\\
 $16\to32$ & 2.1655E-02 & 1.4217E-03 & 2.3296E-02 & 3.9885E-02\\
 $32\to64$ & 2.2211E-03 & 1.6784E-04 & 3.8144E-03 & 8.1710E-03\\
 \bottomrule
 \end{tabular}
\end{table}
\begin{table}[!ht]
 \centering
\captionsetup{font=footnotesize,labelfont=footnotesize}
 \caption{Observed exponents along the joint refinement sequence. The fits use $N_{i+1}=16,32,64$ (all three points) and
 $N_{i+1}=32,64$ (last two points).}
 \label{exp:tab:joint-orders}
 \scriptsize
\setlength{\tabcolsep}{3.5pt}
\renewcommand{\arraystretch}{1.2}
 \begin{tabular}{ccc}
 \toprule
 Diagnostic & All three points & Last two points\\
 \midrule
 Mean & 2.9876 & 3.2854\\
 Componentwise variance & 3.5116 & 3.0825\\
 One-point velocity law & 2.6121 & 2.6105\\
 Two-point velocity law & 2.3891 & 2.2872\\
 One-point speed law & 2.7121 & 2.8433\\
 Two-point speed law & 2.4997 & 2.4834\\
 \bottomrule
 \end{tabular}
\end{table}
\FloatBarrier

\paragraph{Wasserstein distributional diagnostics}
To assess cross-level consistency beyond the first two sample moments,
we compare the one-point and two-point empirical distributions of the
velocity $\boldsymbol u$ and its magnitude $|\boldsymbol u|$ using the
spatially averaged $1$-Wasserstein distances defined above.
\begin{figure}[htbp!]
 \centering
 \includegraphics[width=0.48\textwidth]{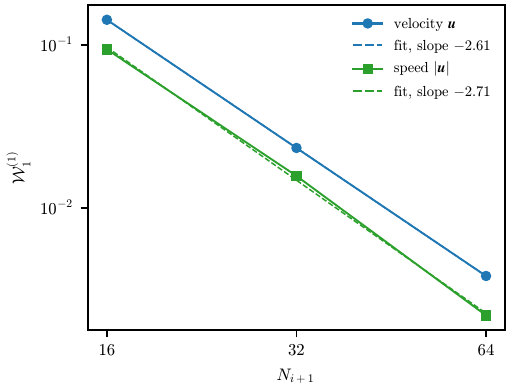}
 \hfill
 \includegraphics[width=0.48\textwidth]{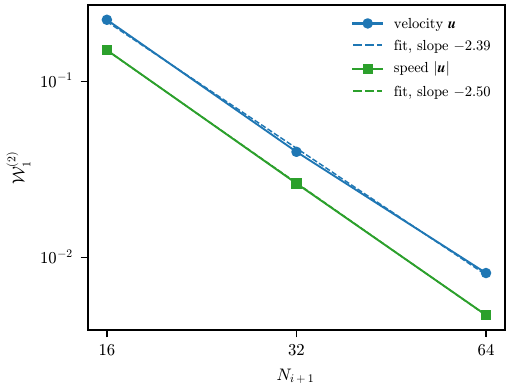}
 \captionsetup{font=footnotesize,labelfont=footnotesize}
 \caption{Adjacent-level spatially averaged $1$-Wasserstein distances at
 $T=1$ for one-point (left) and two-point (right) empirical distributions
 of the velocity $\boldsymbol u$ and its magnitude $|\boldsymbol u|$.
 Fits use all three comparisons.}
 \label{exp:fig:wasserstein}
\end{figure}
\FloatBarrier
The one-point and two-point distances decrease similarly
(Figure~\ref{exp:fig:wasserstein}). The fitted exponents over the three
comparisons (Table~\ref{exp:tab:joint-orders}) are $2.61$ and $2.71$
for one-point velocity and speed, and $2.39$ and $2.50$ for their
two-point counterparts; the empirical laws contain the same sampling
component, which limits the decrease on finer levels.
These are composite rates along the joint refinement path; they do not
identify separate spatial, temporal, or sampling orders.
\FloatBarrier

\subsection{Lid-driven cavity flow with random viscosity: ensemble-size sensitivity}
\label{exp:sec:sampling}

We next consider lid-driven flow on $\Omega=(0,1)^2$ with zero initial velocity,
zero body force, and final time $T=0.5$. The side and bottom walls are
stationary, and the top-wall velocity is
\begin{equation}
 \label{exp:eq:lid}
 \boldsymbol g(x,t)=
 \left(30x^2(1-x)^2[1-\exp(-(t/0.1)^2)],0\right).
\end{equation}
This smooth profile vanishes, with its first derivatives, at the
upper corners and at $t=0$, and is thus compatible with the stationary
sidewalls and zero initial data. Randomness enters only through the
spatially constant viscosity
\begin{equation}
 \label{exp:eq:sampling-viscosity}
 \nu_j=\nu_0(1+\sigma\sqrt{\lambda_0}Y_j),
 \qquad Y_j\overset{\mathrm{iid}}\sim\mathcal U[-\sqrt3,\sqrt3],
\end{equation}
where $\nu_0=10^{-2}$, $\sigma=0.15$, $\lambda_0=\sqrt{\pi L_c}/2$, and
$L_c=0.25$. The sampling study below was computed with the arithmetic mean $\bar\nu_J$ of each group of $J$ members
as shared viscosity. For this choice, unlike in
Section~\ref{exp:sec:joint}, the support $\nu_j\in\nu_0[1-a,1+a]$,
$a=\sigma\sqrt{3\lambda_0}\approx0.173$, does not guarantee
\eqref{assumption:viscosity-fluctuation} a priori; all computed
sampling groups satisfy $\max_j|\nu_j-\bar\nu_J|/\bar\nu_J<1/3$,
with a maximum of $0.197935$, and hence
\eqref{assumption:viscosity-fluctuation}.
The finite element pair, startup step, and SAV
parameters are as in Section~\ref{exp:sec:joint}. We fix a locally refined mesh of $87744$
triangles and $k=1.25\times10^{-3}$ ($N_t=400$).
Since the lid velocity is nonzero, this experiment tests the
implementation with nonhomogeneous boundary data, outside the
homogeneous setting of the analysis.

\emph{Experiment design.}
Following the ensemble Monte Carlo perspective of Luo and
Wang~\cite{LuoWang2018}, we compare $J\in\{10,20,40,80,160\}$ using
$R=10$ independent replications. Within each replication, the ensembles
use prefixes of one pool of 160 i.i.d. inputs. Each group recomputes its
arithmetic mean viscosity and is solved anew. An independent pool of
$J_{\mathrm{ref}}=5000$ inputs provides a single reference ensemble on the
same space--time grid, shared by all $R$ replications. Because the discrete operator
depends on the sampled mean viscosity, shared inputs need not yield
identical trajectories across ensemble sizes.

\emph{Error measures.}
Let $\boldsymbol m_{J,r}^n$ be the sample mean velocity at $t^n$ of
the $J$-member ensemble in replication $r$ and
$\boldsymbol m_{\mathrm{ref}}^n$ the reference mean. Setting
$\boldsymbol e_{J,r}^n=\boldsymbol m_{J,r}^n-\boldsymbol m_{\mathrm{ref}}^n$,
we report
\begin{align}
 \mathcal E_{L^2}(J)&=\max_{1\le n\le N_t}
 \left(\frac1R\sum_{r=1}^R
 \|\boldsymbol e_{J,r}^n\|_{L^2(\Omega)}^2\right)^{1/2},
 \label{exp:eq:error-l2}\\
 \mathcal E_{H^1}(J)&=
 \left(\frac{k}{R}\sum_{r=1}^R\sum_{n=1}^{N_t}
 \|\nabla\boldsymbol e_{J,r}^n\|_{L^2(\Omega)}^2\right)^{1/2}.
 \label{exp:eq:error-h1}
\end{align}
The root mean square in~\eqref{exp:eq:error-l2} is taken before the temporal
maximum;~\eqref{exp:eq:error-h1} uses the spatial $H^1$ seminorm.
All norms are evaluated by finite element integration.

Paired checks with two independent groups of 32 inputs give spatial
$L^2$ and gradient orders $3.0834$ and $2.0474$ on meshes of 5484,
21936, and 87744 triangles, and temporal orders $2.0956$ and $2.0160$
at 100, 200, and 400 steps on the finest mesh. Richardson estimates,
with rates capped at three and two in space and two in time, are each
below $5\%$ of the corresponding $J=160$ sampling error; the largest
ratio is $3.7734\%$. These are empirical bias estimates.

\begin{figure}[htbp!]
	\centering
	\includegraphics[width=0.75\textwidth]
	{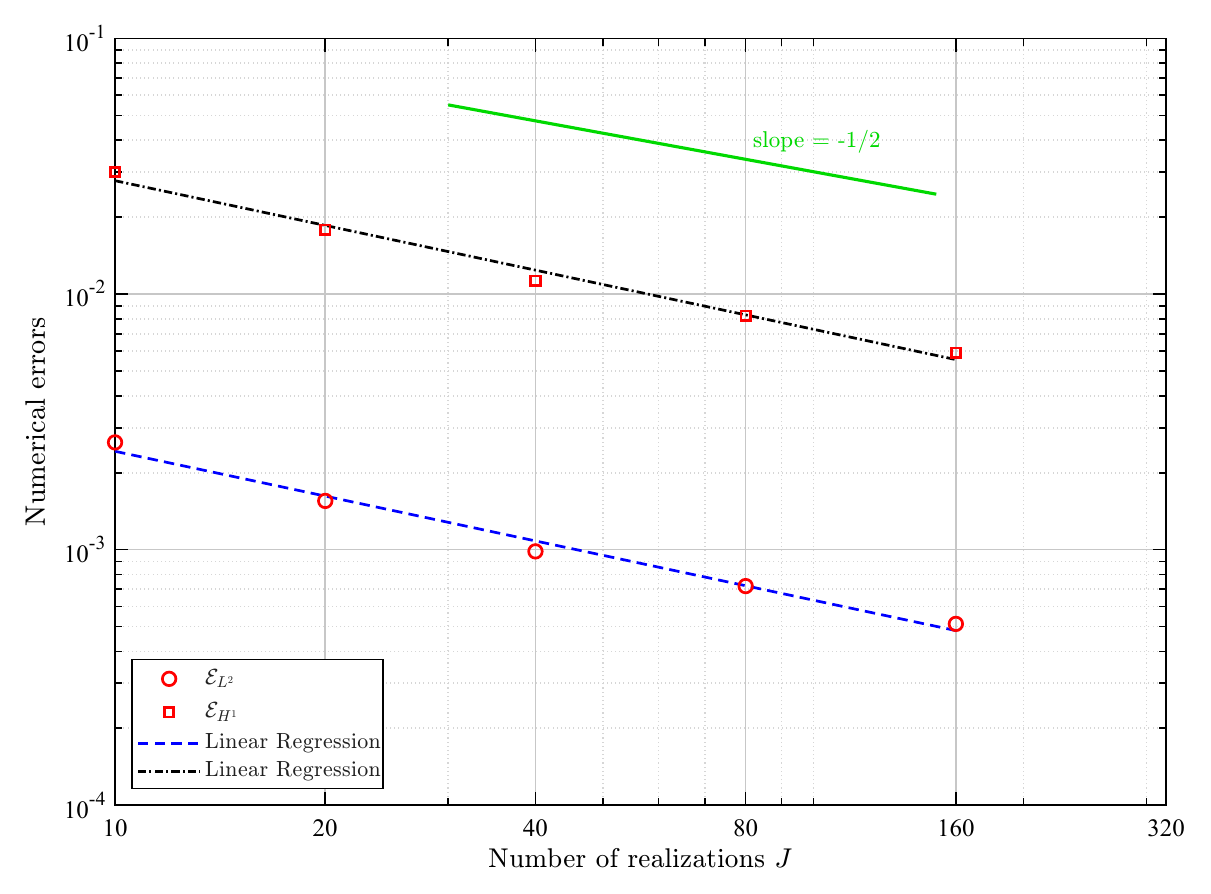}
	\captionsetup{
		font=footnotesize,
		labelfont=footnotesize
	}
	\caption{Ensemble-size sensitivity:
		$\mathcal E_{L^2}(J)$ and
		$\mathcal E_{H^1}(J)$ versus $J$.
		Fits use all five ensemble sizes;
		the green $J^{-1/2}$ line is a visual guide only.}
	\label{exp:fig:sampling}
\end{figure}

\emph{Results.}
Table~\ref{exp:tab:sampling-errors} gives the errors and signed adjacent
log--log slopes
\[
 s_X(J_i)=\frac{\log\mathcal E_X(J_i)-\log\mathcal E_X(J_{i-1})}
 {\log J_i-\log J_{i-1}},\qquad X\in\{L^2,H^1\}.
\]
Both errors decrease with $J$. Least-squares log--log slopes over all
five ensemble sizes are $-0.5826$ and $-0.5814$ for $\mathcal E_{L^2}$ and $\mathcal E_{H^1}$,
respectively (Figure~\ref{exp:fig:sampling}). From 10000 bootstrap
resamples of the ten replication blocks, the 95\% intervals for the
positive exponents are $[0.4795,0.6611]$ and $[0.4754,0.6632]$.
Both contain $1/2$, compatible with $J^{-1/2}$ scaling over this range.
These intervals condition on the fixed, finite reference; the
$J$-dependent shared operator precludes an exact i.i.d. variance identity.

\begin{table}[htbp]
 \centering
\captionsetup{font=footnotesize,labelfont=footnotesize}
 \caption{Ensemble-size sensitivity for random-viscosity lid-driven cavity
 flow: RMS deviations~\eqref{exp:eq:error-l2} and~\eqref{exp:eq:error-h1}
 from one independent $J_{\mathrm{ref}}=5000$ reference shared across
 replications, with adjacent log--log slopes.
 All runs use $T=0.5$, $k=1.25\times10^{-3}$, and $R=10$.}
 \label{exp:tab:sampling-errors}
 \scriptsize
\setlength{\tabcolsep}{3.5pt}
\renewcommand{\arraystretch}{1.2}
 \begin{tabular}{lccccc}
 \toprule
 $J$ & 10 & 20 & 40 & 80 & 160 \\ \midrule
 $\mathcal E_{L^2}(J)$ &
   2.6285E-03 & 1.5514E-03 & 9.8451E-04 & 7.2018E-04 & 5.1230E-04 \\
 slope & -- & $-0.76$ & $-0.66$ & $-0.45$ & $-0.49$ \\
 \midrule
 $\mathcal E_{H^1}(J)$ &
   3.0052E-02 & 1.7851E-02 & 1.1287E-02 & 8.2350E-03 & 5.8987E-03 \\
 slope & -- & $-0.75$ & $-0.66$ & $-0.45$ & $-0.48$ \\
 \bottomrule
 \end{tabular}
\end{table}

\FloatBarrier

\paragraph{Computational cost}
To quantify the gain from the shared operator, we time the lid-driven
configuration of this section with the skew-symmetric form and a
deterministic shared viscosity $\bar\nu=\nu_0$, fixed before sampling.
Then $\nu_j/\bar\nu\in[1-a,1+a]$, so
\eqref{assumption:viscosity-fluctuation} holds a priori, and a single
factorization serves every ensemble size. Moreover, the members are then
independent, so the ensemble with $J$ members is part of the ensemble
with $2J$ members, and an ensemble-size study needs only one computation
with the largest ensemble. We use a uniform $64\times64$ mesh with
$8192$ triangles, so that each Brinkman solve has $36{,}482$ unknowns,
$k=1.25\times10^{-3}$, $T=0.5$, and the solver PARDISO with one thread,
and compare four implementations. In (E), the ensemble scheme solves the
$2J$ right-hand sides of a step in one call; (E$_{\rm seq}$) uses the
same scheme and factorization but solves the members one after another.
(M$_{\rm SAV}$) computes member by member with the mr-SAV--BDF2 scheme and
each member's own viscosity treated implicitly, which needs two numerical
factorizations per member, one for the startup step; the symbolic
analysis is reused.
(M$_{\rm LI}$) computes member by member with the standard linearly
implicit BDF2 scheme, with the convective term
$b_{\rm sk}(\bm u_j^{*,n+1},\bm u_j^{n+1},\bm v)$, the member's own
viscosity, and no auxiliary variable, which needs a numerical
refactorization at every step. Table~\ref{exp:tab:cost} reports the CPU
times. For $J=160$, (E) takes $750$\,s, whereas (M$_{\rm SAV}$) and
(M$_{\rm LI}$) take about $2.0$ and $11$ times as long; these ratios hardly
depend on $J$. The member-by-member computations
(M$_{\rm SAV}$) cost about as much as (E$_{\rm seq}$), since the
factorizations of a member are cheap compared with its $400$ steps of
solves.
The gain of (E) over (M$_{\rm SAV}$) thus comes from solving all
right-hand sides of a step in one call, which the shared operator makes
possible; the two factorizations of (E) take only $0.9$\,s. The
gain over (M$_{\rm LI}$) comes from avoiding a factorization at every
step. Since the right-hand sides of a step are independent, the batched
solve can also be parallelized, which we do not pursue here. As a check, the kinetic energies at $T$ computed by (E)
differ from those of (M$_{\rm SAV}$) and (M$_{\rm LI}$) by relative
amounts of at most $4\times10^{-8}$ and $6\times10^{-7}$, respectively, and
$\max_j|q_j^N-1|\le9\times10^{-11}$ in (E).

\begin{table}[htbp]
 \centering
 \captionsetup{font=footnotesize,labelfont=footnotesize}
 \caption{CPU times in seconds for the lid-driven cavity with
 $\bar\nu=\nu_0$, a uniform $64\times64$ mesh, $k=1.25\times10^{-3}$, and
 $T=0.5$ ($400$ steps). All times include the factorizations and
 symbolic analyses. Values in parentheses are extrapolated linearly in
 $J$ from $10$ members for (M$_{\rm SAV}$), $2$ members for
 (M$_{\rm LI}$), and $J\le40$ for (E$_{\rm seq}$). The last two columns
 give the ratios to (E).}
 \label{exp:tab:cost}
 \scriptsize
 \setlength{\tabcolsep}{4pt}
 \renewcommand{\arraystretch}{1.15}
 \begin{tabular}{rrrrrrr}
 \toprule
 $J$ & (E) & (E$_{\rm seq}$) & (M$_{\rm SAV}$) & (M$_{\rm LI}$) & M$_{\rm SAV}$/E & M$_{\rm LI}$/E \\
 \midrule
10 & 47.7 & 97.1 & 94.9 & (521) & 2.0 & 11 \\
20 & 87.5 & 198.9 & (189) & (1043) & 2.2 & 12 \\
40 & 165.2 & 381.2 & (377) & (2085) & 2.3 & 13 \\
80 & 358.2 & (775) & (754) & (4171) & 2.1 & 12 \\
160 & 749.9 & (1549) & (1507) & (8342) & 2.0 & 11 \\
 \bottomrule
 \end{tabular}
\end{table}

\subsection{Kolmogorov flow}
\label{sec:Kolmogorovflow}
We consider the two-dimensional Navier--Stokes equations with Kolmogorov
forcing and periodic boundary conditions on $\Omega=(0,2\pi)^2$. For
each ensemble member, the forcing is $\bm F_j=(\nu_jm^3\cos(my),0)^T$
with $m=2$. For every $\nu_j$, the laminar Kolmogorov flow
$\bm u=(m\cos(my),0)^T$, with constant pressure, is a steady solution,
common to all members. We take $J=3$ with viscosities
$(\nu_1,\nu_2,\nu_3)=(0.0095,0.0100,0.0105)$ and the shared viscosity
$\bar\nu=0.01$, so that $\nu_j/\bar\nu\in[0.95,1.05]$ and
$Re=1/\bar\nu=100$, based on the amplitude $m$ and the length scale
$1/m$ of the laminar flow.

For this periodic problem we use a Fourier--Galerkin discretization,
implemented with FFTs, with $N^2$ grid points and the wave vectors
$\bm\kappa$ with $|\kappa_1|,|\kappa_2|\le K:=\lfloor N/3\rfloor$; here
$N=64$. The discrete velocities are exactly divergence-free with zero
mean, and the convective term is computed without aliasing error, so
\eqref{eq:skew-general} holds exactly; the arguments of
Section~\ref{sec:long-time-stability} then carry over, with periodic
instead of no-slip conditions. The shared operator is diagonal in
Fourier space, so a time step costs a few FFTs per member. All members
are initialized with a common velocity field obtained after a spin-up of
$T_s=100$, with $q_j^0=1$. The spin-up starts from
$(\partial_y\psi_s,-\partial_x\psi_s)^T$,
$\psi_s(x,y)=\sin(2y)+10^{-3}\sin(2x)\sin(2y)$, and uses $\nu_s=0.01$,
$\bm F_s=(8\nu_s\cos(2y),0)^T$, and the first-order semi-implicit Euler
scheme (implicit viscosity, explicit convection, no auxiliary variable)
with time step $0.002$. The simulation clock is reset after the spin-up,
and time integration begins with the mr-SAV--BE step of
Appendix~\ref{subsec:be-startup}, followed by the mr-SAV--BDF2 scheme,
both with $k=0.002$.

\subsubsection{Long-time stability}
\label{subsubsec:kolmogorov-stability}
\begin{figure}[tbp]
  \centering
  \includegraphics[width=\linewidth]{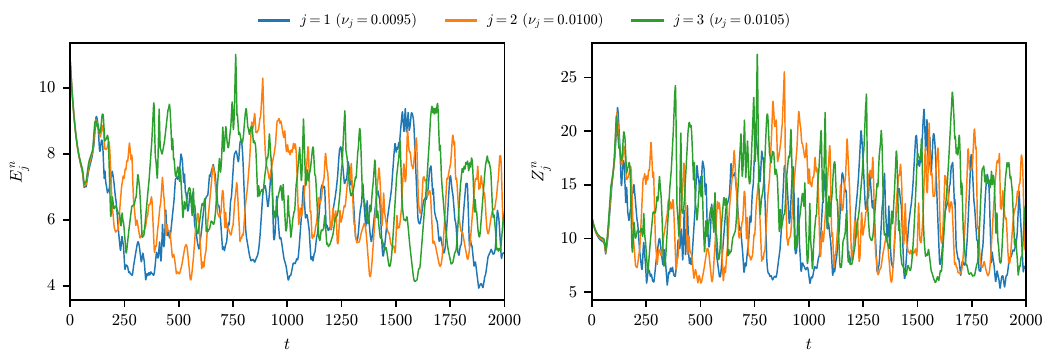}
  \captionsetup{font=footnotesize,labelfont=footnotesize}
  \caption{Fourier--Galerkin computation of Kolmogorov flow with
    $\gamma=10^3$, $N=64$, and $k=0.002$: kinetic energy $E_j^n$ (left)
    and enstrophy $Z_j^n$ (right) of the three members. Time is
    measured from the end of the common spin-up.}
  \label{fig:kolmogorov-energy-enstrophy}
\end{figure}

We examine the long-time behavior of the ensemble with $\gamma=10^3$ up
to $T=2000$ through the kinetic energy and enstrophy of each member,
\[
  E_j^n=\frac12\|\bm u_j^n\|_{L^2(\Omega)}^2,
  \qquad
  Z_j^n=\frac12\|\partial_xu_{j,2}^n-\partial_yu_{j,1}^n\|_{L^2(\Omega)}^2,
\]
where $u_{j,1}^n$ and $u_{j,2}^n$ are the components of $\bm u_j^n$.
Figure~\ref{fig:kolmogorov-energy-enstrophy} shows their evolution.
The flow is chaotic, and for all three members both quantities
fluctuate within bounded ranges, $E_j^n\in[3.9,11.0]$ and $Z_j^n\in[5.4,27.2]$ for $100\le t\le2000$, with no sustained
growth.

\subsubsection{Mean-reverting property}
\label{subsubsec:kolmogorov-mean-reversion}
\begin{figure}[tbp]
  \centering
  \includegraphics[width=\linewidth]{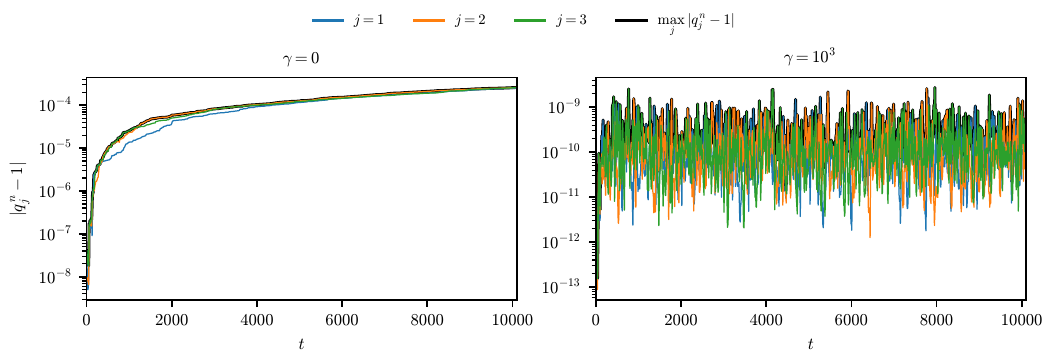}
  \captionsetup{font=footnotesize,labelfont=footnotesize}
  \caption{Deviations $|q_j^n-1|$ without mean reversion, $\gamma=0$
    (left), and with mean reversion, $\gamma=10^3$ (right), for $J=3$,
    $N=64$, $k=0.002$, and $T=10100$. Each curve shows maxima over
    consecutive time windows of length $10$; thick black curves denote
    $\max_j|q_j^n-1|$.}
  \label{fig:kolmogorov-mean-reversion}
\end{figure}

We next compare $\gamma=0$ (no relaxation; the analysis assumes
$\gamma>0$) and $\gamma=10^3$, with the same discretization and initial
data, $k=0.002$, and $T=10100$.
Figure~\ref{fig:kolmogorov-mean-reversion} shows the deviations
$|q_j^n-1|$ and their ensemble envelope $\max_j|q_j^n-1|$. Without mean reversion, the deviations accumulate: the envelope grows
about linearly in time, from $2.6\times10^{-5}$ at $t=1000$ to
$2.6\times10^{-4}$ at $T=10100$, since nothing drives $q_j$ back to one.
With $\gamma=10^3$, the deviations stay below $3\times10^{-9}$
throughout, without growth, because the relaxation damps each deviation
on the time scale $1/\gamma$. Mean reversion thus removes the drift of
the auxiliary variable at no extra cost. The drift without mean
reversion is slow here because $b(\bm v,\bm v,\bm v)=0$ holds exactly; for
discretizations without this property, the source term of
\eqref{scheme:sav-bdf2-q} contains $b(\bm u_j^{*,n+1},\bm u_j^{*,n+1},\bm u_j^{*,n+1})$, which does not vanish
as $k\to0$, and the drift can be much faster.

\FloatBarrier

\subsection{Large time steps and sharpness of the viscosity condition}
\label{sec:large-steps}
We now test Theorem~\ref{thm-long-time-stability} and
Lemma~\ref{prop:sharp-threshold} directly on the Kolmogorov flow of
Section~\ref{sec:Kolmogorovflow}, with a wider spread of viscosities,
$(\nu_1,\nu_2,\nu_3)=(0.005,0.0075,0.0125)$, the forces
$\bm F_j=(\nu_jm^3\cos(my),0)^T$ with $m=2$, and $\gamma=10^3$. The
initial velocity is the field $(\partial_y\psi_s,-\partial_x\psi_s)^T$
of Section~\ref{sec:Kolmogorovflow}, without spin-up, plus a
divergence-free random perturbation of energy $10^{-16}$ in all resolved
modes, and the first level is computed by the mr-SAV--BE step. We use
the Fourier--Galerkin discretization of Section~\ref{sec:Kolmogorovflow};
its discrete Stokes eigenvalues are $|\bm\kappa|^2$, with largest value
$\lambda_{\max}=2K^2$.

We compare two shared viscosities: (a) $\bar\nu=0.0095$, for which
$\nu_3/\bar\nu\approx1.32$, so that
\eqref{assumption:viscosity-fluctuation} holds close to the threshold;
and (b) the ensemble mean $\bar\nu\approx0.0083$, for which
$\nu_3/\bar\nu=1.5$ violates \eqref{assumption:viscosity-fluctuation},
as does the midrange (cf.\ Remark~\ref{rem:nubar-choice}(iii)). In case
(b), Lemma~\ref{prop:sharp-threshold}(i), applied to each Fourier
mode, predicts that the Stokes part
of the scheme is stable exactly for $k<k_c$, where
\[
 k_c:=\frac{4}{(3\nu_3'-\bar\nu)\lambda_{\max}}.
\]

\begin{figure}[tbp]
  \centering
  \includegraphics[width=\linewidth]{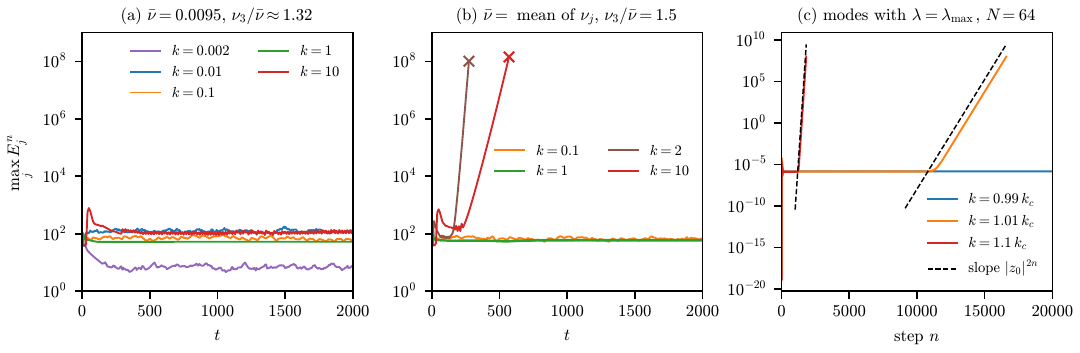}
  \captionsetup{font=footnotesize,labelfont=footnotesize}
  \caption{Fourier--Galerkin computations of Kolmogorov flow with
  $(\nu_1,\nu_2,\nu_3)=(0.005,0.0075,0.0125)$. (a), (b): $\max_jE_j^n$ for
  $N=64$, with $\bar\nu=0.0095$ (a) and $\bar\nu$ equal to the ensemble
  mean (b); crosses mark blow-up ($\max_jE_j^n>10^8$). (c): energy of the
  modes of member $3$ with $\lambda=\lambda_{\max}$ for $N=64$ and
  $\bar\nu$ equal to the ensemble mean, near the predicted critical step
  $k_c$; dashed lines have slope $|z_0|^{2n}$.}
  \label{fig:large-steps}
\end{figure}

\paragraph{Long-time behavior for large time steps}
Figure~\ref{fig:large-steps}(a,b) and Table~\ref{tab:large-steps} show
$\max_jE_j^n$ for $N=64$ and $T=2000$. In case (a), the energies remain
bounded for all time steps from $k=0.002$ to $k=10$, as
Theorem~\ref{thm-long-time-stability} predicts; the largest value, at
$k=10$, occurs in the initial transient. For $k=0.002$, all $q_j$ stay
within $2.2\times10^{-7}$ of one, and the energies are of the size
observed in Section~\ref{sec:Kolmogorovflow}. For larger steps, which do
not resolve the convective time scale of the smallest resolved scales,
$q_j$ drops well below one, first for the member with the smallest
viscosity ($j=1$ at $k=0.01$) and then for all members; this damps the
computed convection, and the energies settle at a higher level. The
solutions stay bounded but are not accurate: unconditional stability
does not imply accuracy for large steps. In case (b), $k_c=1.088$ for
$N=64$. The computations with $k=0.1$ and $k=1$ remain bounded, whereas
those with $k=2$ and $k=10$ blow up, at $t=274$ and $t=570$,
respectively.

\paragraph{Sharpness of the threshold}
Table~\ref{tab:threshold} reports case (b) for $N=64,128,256$ and
$k=0.99k_c$, $1.01k_c$, and $1.1k_c$, with at most $2\times10^4$ steps.
At $k=0.99k_c$ all computations remain bounded; at $1.01k_c$ and
$1.1k_c$ they blow up. In the exponential phase before blow-up, the
modes with $\lambda=\lambda_{\max}$ grow by the factor $|z_0|$ per step,
where $z_0<-1$ is the root of $P_{k\lambda_{\max}}$ for member $3$ in
Lemma~\ref{prop:sharp-threshold}(i)
(Figure~\ref{fig:large-steps}(c)); the fitted factors agree with $|z_0|$
to within $2\times10^{-9}$. Thus, in this test, the convective term does
not shift the onset of instability predicted for the Stokes part. Since
$k_c\propto N^{-2}$, the stable time steps shrink to zero under spatial
refinement, in accordance with Lemma~\ref{prop:sharp-threshold}(ii):
beyond the threshold, no time step is stable uniformly in the spatial
resolution.

\begin{table}[tbp]
  \centering
  \captionsetup{font=footnotesize,labelfont=footnotesize}
  \caption{Case (a), $\bar\nu=0.0095$, $N=64$, $T=2000$: maximal
  energy, mean energy over $t\in[1000,2000]$, and smallest auxiliary
  variable over $t\in[1000,2000]$.}
  \label{tab:large-steps}
  \scriptsize
  \setlength{\tabcolsep}{6pt}
  \renewcommand{\arraystretch}{1.15}
  \begin{tabular}{cccc}
    \toprule
    $k$ & $\max_{n,j}E_j^n$ & mean of $\max_jE_j^n$ & $\min_{n,j}q_j^n$ \\
    \midrule
    $0.002$ & 39.5 & 6.8 & $1-1.0\times10^{-9}$ \\
    $0.01$ & 178.4 & 123.9 & 0.161 \\
    $0.1$ & 92.7 & 68.2 & 0.042 \\
    $1$ & 103.8 & 52.7 & 0.0081 \\
    $10$ & 785.6 & 112.2 & 0.0004 \\
    \bottomrule
  \end{tabular}
\end{table}

\begin{table}[tbp]
  \centering
  \captionsetup{font=footnotesize,labelfont=footnotesize}
  \caption{Case (b), $\bar\nu$ equal to the ensemble mean
  ($\nu_3/\bar\nu=1.5$): predicted critical step $k_c$, maximal energy
  over $2\times10^4$ steps at $k=0.99k_c$, step $n_{\rm b}$ at which
  $\max_jE_j^n$ exceeds $10^8$, and the fitted growth factor per step of
  the modes of member $3$ with $\lambda=\lambda_{\max}$, to be compared
  with the root modulus $|z_0|$ of Lemma~\ref{prop:sharp-threshold}(i).}
  \label{tab:threshold}
  \scriptsize
  \setlength{\tabcolsep}{3.5pt}
  \renewcommand{\arraystretch}{1.15}
  \begin{tabular}{cccccccc}
    \toprule
    & & & $0.99k_c$ & \multicolumn{2}{c}{$1.01k_c$, $|z_0|=1.003051$} &
    \multicolumn{2}{c}{$1.1k_c$, $|z_0|=1.028384$} \\
    \cmidrule(lr){4-4}\cmidrule(lr){5-6}\cmidrule(lr){7-8}
    $N$ & $\lambda_{\max}$ & $k_c$ & $\max_{n,j}E_j^n$ & $n_{\rm b}$ & fitted &
    $n_{\rm b}$ & fitted \\
    \midrule
    64 & 882 & 1.0884 & 137.5 & 16594 & 1.003051 & 1822 & 1.028384 \\
    128 & 3528 & 0.2721 & 50.9 & 16663 & 1.003051 & 1883 & 1.028384 \\
    256 & 14450 & 0.0664 & 45.0 & 16878 & 1.003051 & 1886 & 1.028384 \\
    \bottomrule
  \end{tabular}
\end{table}
\FloatBarrier

\section{Concluding remarks}
\label{sec:conclusions}
The mr-SAV--BDF2 ensemble method combines three properties. All members
share one time-independent shifted Stokes operator, with two linear
solves per member and step. For constant sample viscosities with
$0<\nu_j<\frac43\bar\nu$ and uniformly bounded forcing, the scheme is
dissipative uniformly in time for every $k>0$, with ensemble bounds and
an absorbing radius independent of the time step and the initial data
(Theorem~\ref{thm-long-time-stability} and
Corollary~\ref{cor:absorbing-radius}). For a given $\bar\nu$, the upper
bound is sharp (Lemma~\ref{prop:sharp-threshold}), and it restricts the
choice of $\bar\nu$, not the ensemble (Remark~\ref{rem:nubar-choice}).
The analysis requires a convective form with $b(\bm v,\bm v,\bm v)=0$,
which holds for the time-discrete scheme and for exactly
divergence-free or skew-symmetric spatial discretizations. The unconditional long-time stability
extends to spatially varying random viscosities under a signed pointwise
condition that keeps $\nu_j$ within $(\bar\nu/3,\frac43\bar\nu)$,
bounded away from the endpoints uniformly in space, and on the full
range for proportional spatial profiles
(Theorem~\ref{thm:random-viscosity-stability}). In two space dimensions,
for constant viscosities and sufficiently regular solutions, the
velocity, the auxiliary variable, and the discrete time primitive of the
pressure are second-order accurate in time on bounded intervals
(Appendix~\ref{app:error-analysis}).

Numerical experiments confirm the long-time stability theory. Within
the setting of the analysis, a Fourier--Galerkin discretization gives
bounded solutions for all tested time steps, up to $k=10$, when the
viscosity condition holds, although accuracy is lost once the step no
longer resolves the convective time scale; when the condition is
violated, the solutions blow up as soon as the time step exceeds the
predicted critical value, which tends to zero under spatial refinement
(Section~\ref{sec:large-steps}). Under the viscosity condition, the
time step can thus be chosen for accuracy alone: no member can blow up
over a long horizon, and in our tests a step that is too large for a
member showed itself through $q_j$ falling below one rather than through
blow-up. The manufactured-solution and cavity
experiments show that the method is accurate, member by member and in
its ensemble statistics, also for viscosity deviations of a quarter of
the shared viscosity
(Sections~\ref{subsec:convergence-tests}--\ref{exp:sec:sampling}).
Computing all members with the shared
operator, with all right-hand sides of a step solved together, is about
twice as fast as computing them one by one with the same scheme, and
about eleven times as fast as with a standard linearly implicit scheme
(Section~\ref{exp:sec:sampling}). Kolmogorov flow shows that mean
reversion suppresses the drift of the auxiliary variable and keeps
long-time computations robust (Section~\ref{sec:Kolmogorovflow}).

\FloatBarrier
\appendix
\section{\texorpdfstring{Backward Euler startup and finite-time error analysis}{Backward Euler startup and finite-time error analysis}}
\label{app:error-analysis}

Appendix~\ref{subsec:be-startup} treats the backward Euler startup step.
From Appendix~\ref{subsec:error-setting} on, we
prove second-order temporal convergence of the mr-SAV--\allowbreak BDF2
ensemble scheme in two space dimensions on $[0,T]$ under the full-range
condition \eqref{assumption:viscosity-fluctuation}, without a time-step
restriction. The nonlinear terms are controlled by the modified basic
energy of Theorem~\ref{thm-long-time-stability}, augmented by a filter
energy for the gradients. After estimating the consistency,
extrapolation, and startup errors, we derive a filtered error-energy
recursion for the velocity and auxiliary variable, and recover a
second-order bound for the time-integrated pressure from the inf--sup
condition. Spatially varying viscosities are discussed in
Remark~\ref{rem:spatial-error}. Two-level vectors and $G$-norms
are those of \eqref{eq:stability-two-level-vectors} and
\eqref{def:G-matrix}.

\subsection{\texorpdfstring{The backward Euler startup step}{The backward Euler startup step}}
\label{subsec:be-startup}

The stability results of Section~\ref{sec:long-time-stability} hold for
arbitrary starting levels. In the error analysis, the first level is
computed from $(\bm u_j^0,q_j^0)$, $q_j^0=1$, by the mr-SAV backward
Euler (mr-SAV--BE) step: find
$(\bm u_j^1,p_j^1,q_j^1)\in\bm X\times Q\times\mathbb R$ such that
\begin{subequations}
\label{scheme:fsav-be-startup}
\begin{align}
  &\Bigl(\frac{\bm u_j^1-\bm u_j^0}{k},\bm v\Bigr)
  +q_j^1\,b(\bm u_j^0,\bm u_j^0,\bm v)
  +\bar\nu(\nabla\bm u_j^1,\nabla\bm v)
  \nonumber\\
  &\qquad
  +\nu_j'(\nabla\bm u_j^0,\nabla\bm v)
  -(p_j^1,\nabla\cdot\bm v)
  =\langle\bm F_j^1,\bm v\rangle,
  \label{scheme:fsav-be-startup-momentum}
  \\
  &(\nabla\cdot\bm u_j^1,r)=0,
  \label{scheme:fsav-be-startup-continuity}
  \\
  &\frac{q_j^1-q_j^0}{k}+\gamma q_j^1
  -b(\bm u_j^0,\bm u_j^0,\bm u_j^1)=\gamma,
  \label{scheme:fsav-be-startup-q}
\end{align}
\end{subequations}
for all $(\bm v,r)\in\bm X\times Q$. That is, the BDF2 differences are
replaced by first-order differences, the extrapolated velocity is replaced
by $\bm u_j^0$, and all new-time-level quantities are evaluated at $t^1$.
With $1/k$ in place of $3/(2k)$ in \eqref{def:brinkman-bilinear-form}
and \eqref{def:scalar-denominator}, the proof of
Proposition~\ref{prop:stepwise-solvability} gives unique solvability, a
scalar denominator of at least $1/k+\gamma$, and $\bm u_j^1\in\bm V$.

The next lemma bounds the first level by the data under the full-range
condition, as used in Remark~\ref{rem:asymmetric-startup} and in the
error analysis. For $m=0,1$ set
\begin{equation}
\label{eq:BE-energy-definition}
  \mathcal E_{j,\mathrm{BE}}^m
  :=\|\bm u_j^m\|^2+|q_j^m|^2.
\end{equation}

\begin{lemma}
\label{lem:startup-stability}
Let $\bm u_j^0\in\bm V$, $q_j^0\in\mathbb R$, and $k>0$, and assume
\eqref{assumption:viscosity-fluctuation} and
\eqref{assumption:discrete-force}. Then $\rho<\bar\nu$, so that
$\beta:=\bar\nu-\rho/2>\bar\nu/2$, and for $j=1,\ldots,J$ the solution of
\eqref{scheme:fsav-be-startup} satisfies
\begin{equation}
\label{eq:asymmetric-startup}
  \mathcal E_{j,\mathrm{BE}}^1
  +k\beta\|\nabla\bm u_j^1\|^2
  +k\gamma|q_j^1|^2
  \le\mathcal E_{j,\mathrm{BE}}^0
  +k\rho\|\nabla\bm u_j^0\|^2
  +k\Bigl(\frac{M_F^2}{\beta}+\gamma\Bigr).
\end{equation}
\end{lemma}

\begin{proof}
By \eqref{assumption:viscosity-fluctuation}, $|\nu_j'|<\bar\nu$ for
every $j$, hence $\rho<\bar\nu$. Test
\eqref{scheme:fsav-be-startup-momentum} with $\bm v=k\bm u_j^1$,
multiply \eqref{scheme:fsav-be-startup-q} by $kq_j^1$, and add. The
pressure term vanishes by \eqref{scheme:fsav-be-startup-continuity}, and
the convective terms $\pm kq_j^1b(\bm u_j^0,\bm u_j^0,\bm u_j^1)$ cancel.
The identity $(a-b,a)=\frac12(\|a\|^2-\|b\|^2+\|a-b\|^2)$ yields
\begin{equation}
\label{eq:proof-startup-exact-balance}
\begin{aligned}
  &\frac12
   \left(
     \mathcal E_{j,\mathrm{BE}}^1
     -\mathcal E_{j,\mathrm{BE}}^0
   \right)
  +\frac12\|\bm u_j^1-\bm u_j^0\|^2
  +\frac12|q_j^1-q_j^0|^2
  +k\bar\nu\|\nabla\bm u_j^1\|^2
  \\&\qquad+k\gamma|q_j^1|^2=
  -k\nu_j'(\nabla\bm u_j^0,\nabla\bm u_j^1)
  +k\langle\bm F_j^1,\bm u_j^1\rangle
  +k\gamma q_j^1.
\end{aligned}
\end{equation}
By $|\nu_j'|\le\rho$, \eqref{eq:duality-estimate},
\eqref{assumption:discrete-force}, and Young's inequality, the right-hand
side is at most
$\frac{k\rho}{2}(\|\nabla\bm u_j^0\|^2+\|\nabla\bm u_j^1\|^2)
+\frac{k\beta}{2}\|\nabla\bm u_j^1\|^2+\frac{kM_F^2}{2\beta}
+\frac{k\gamma}{2}|q_j^1|^2+\frac{k\gamma}{2}$. Multiplying by $2$,
dropping the squared increments, and using $2\bar\nu-\rho-\beta=\beta$
proves \eqref{eq:asymmetric-startup}.
\end{proof}

Since $\beta>\bar\nu/2$, \eqref{eq:asymmetric-startup} bounds
$\|\bm u_j^1\|^2+|q_j^1|^2$ and $k\|\nabla\bm u_j^1\|^2$ in terms of the
data, uniformly for $k$ in bounded intervals. With the constants of
\eqref{eq:full-range-constants}, the modified energy
\eqref{eq:full-range-energy} at the first level then satisfies
$\widetilde{\mathcal H}_j^1\le\max\{C_G+k\varepsilon\Lambda,C_q\}
(\mathcal E_{j,\mathrm{BE}}^1+\mathcal E_{j,\mathrm{BE}}^0)$, as claimed in
Remark~\ref{rem:asymmetric-startup}.

\subsection{\texorpdfstring{Setting and main result}{Setting and main result}}
\label{subsec:error-setting}

Assume \eqref{assumption:viscosity-fluctuation}, and use $\eta_j$,
$G_{\eta_j}$, $P_j$, $D$, $\sigma$, and $c_G$ from
\eqref{eq:full-range-multiplier}--\eqref{eq:full-range-constants}.
All constants below are uniform in $j$ for a fixed ensemble; for
families of ensembles, the ratios $\eta_j$ must stay in a common compact
subinterval of $(-1/3,1)$, and the remaining data, regularity, and
consistency bounds must be uniform. Fix $T>0$ and $N\ge2$ with
$Nk\le T$, so that $0<k\le T/2$. For the exact solution,
$b(\bm u_j,\bm u_j,\bm u_j)=0$, so
\eqref{eq:continuous-fsav-system2} with $q_j(0)=1$ gives $q_j\equiv1$,
and the auxiliary-variable error is measured against $1$. Set
\begin{gather}
 \bm U_j^n:=\bm u_j(t^n),\qquad
 \bm U_j^{*,n+1}:=2\bm U_j^n-\bm U_j^{n-1},\qquad
 e_{\bm u,j}^n:=\bm u_j^n-\bm U_j^n,\qquad
 e_{q,j}^n:=q_j^n-1,
 \label{eq:error-exact-levels}\\
 e_{\bm u,j}^{*,n+1}:=2e_{\bm u,j}^n-e_{\bm u,j}^{n-1},\qquad
 E_{V,j}^n:=(e_{\bm u,j}^n,e_{\bm u,j}^{n-1})^T,\qquad
 E_{Q,j}^n:=(e_{q,j}^n,e_{q,j}^{n-1})^T,
 \label{eq:error-notation}\\
 e_{\bm u,j}^0=\bm0,\qquad e_{q,j}^0=0
 \qquad\text{(exact initial data)}.
 \label{eq:error-initial-values}
\end{gather}
The BDF2, extrapolation, and force residuals are
\begin{equation}
\begin{gathered}
 \bm R_{t,j}^{n+1}:=\partial_t\bm u_j(t^{n+1})-\delta_t\bm U_j^{n+1},\qquad
 \bm R_{\mathrm{ext},j}^{n+1}:=\bm U_j^{n+1}-\bm U_j^{*,n+1},\\
 \bm R_{F,j}^{n+1}:=\bm F_j^{n+1}-\bm F_j(t^{n+1}),
\end{gathered}
 \label{eq:error-residuals}
\end{equation}
the first two for $n\ge1$, the last for $n\ge0$. Endpoint sampling,
$\bm F_j^{n+1}=\bm F_j(t^{n+1})$, gives $\bm R_{F,j}^{n+1}=\bm0$. In
general, we assume the second-order condition
\eqref{eq:error-force-consistency} below on the BDF2 steps and, for the
accuracy of the first step (Lemma~\ref{lem:error-startup}),
$\|\bm R_{F,j}^1\|\le Ck$. The cell average
\eqref{def:discrete-force-average} does not, in general, satisfy
\eqref{eq:error-force-consistency}.

Set
\begin{equation}
 \kappa:=\frac{\bar\nu}{4D}>0,\qquad
 c_{\eta,j}:=\frac{1+3\eta_j}{4}>0,\qquad
 \bm z_j^{n+1}:=e_{\bm u,j}^{n+1}-\eta_j e_{\bm u,j}^{*,n+1},
 \label{eq:error-filter-parameters}
\end{equation}
and, for $n\ge1$, define the gradient filter energy, the error energy,
and the corresponding energy of the numerical solution by
\begin{align}
 \mathcal F_{\nabla e,j}^n&:=\|\nabla E_{V,j}^n\|_{P_j}^2,\qquad
 \mathcal H_{e,j}^n:=\|E_{V,j}^n\|_{G_{\eta_j}}^2
                   +|E_{Q,j}^n|_G^2+k\kappa\mathcal F_{\nabla e,j}^n,
 \label{eq:error-compensated-energy}\\
 \mathcal A_j^n&:=\mathcal E_j^n+k\kappa\|\nabla V_j^n\|_{P_j}^2,
 \label{eq:error-numerical-gradient-energy}
\end{align}
with the gradient acting on both components. The filter estimate
\eqref{eq:full-range-filter-dissipation} holds in any Hilbert space;
applied to the gradients in
$E_{V,j}^{n+1}=B_{\eta_j}E_{V,j}^n+(\bm z_j^{n+1},\bm0)^T$, it gives
\begin{equation}
 \mathcal F_{\nabla e,j}^{n+1}-\mathcal F_{\nabla e,j}^n
       +\sigma\mathcal F_{\nabla e,j}^{n+1}
 \le D\|\nabla\bm z_j^{n+1}\|^2.
 \label{eq:error-gradient-filter}
\end{equation}
Since $P_j\ge I$, the definition of $c_G$ gives
\begin{equation}
 \mathcal H_{e,j}^n\ge
 c_G\bigl(\|e_{\bm u,j}^n\|^2+\|e_{\bm u,j}^{n-1}\|^2
           +|e_{q,j}^n|^2+|e_{q,j}^{n-1}|^2\bigr)
 +k\kappa\bigl(\|\nabla e_{\bm u,j}^n\|^2
                  +\|\nabla e_{\bm u,j}^{n-1}\|^2\bigr).
 \label{eq:error-filter-coercivity}
\end{equation}

We assume, uniformly in $j$,
\begin{equation}
    \bm u_j\in L^\infty(0,T;H^2(\Omega)^2\cap\bm V),
    \quad
    \partial_{tt}\bm u_j\in L^\infty(0,T;H^1(\Omega)^2),
    \quad
    \partial_{ttt}\bm u_j\in L^2(0,T;\bm X').
\label{eq:error-main-regularity}
\end{equation}
In particular,
$\bm u_j\in C^1([0,T];H^1(\Omega)^2)\cap C_w([0,T];H^2(\Omega)^2)$, so
the nodal values of $\bm u_j$ and $\partial_t\bm u_j$ are defined. We
also assume $\bm F_j\in C([0,T];\bm X')$, so that the exact weak
momentum identity holds at every node. Regularity up to $t=0$ is an
assumption on the exact solution; for general data it requires
compatibility conditions on the initial data and forcing
\cite{HeywoodRannacher1982}. Finally, we assume
\begin{equation}
    k\sum_{n=1}^{N-1}
       \|\bm R_{F,j}^{n+1}\|_{-1}^2
    \le C_{F,T}k^4.
    \label{eq:error-force-consistency}
\end{equation}

Throughout, $C$ denotes a generic positive constant independent of $k$,
$n$, and $j$. The constant $C_T$ may in addition depend on $T$, the
domain, $M_F$, $\gamma^{\pm1}$, $\bar\nu^{\pm1}$, $\rho$, $c_G^{-1}$, $D$,
$\sigma^{-1}$, $\|\bm u_j^0\|_{H^1}$, the bounds in
\eqref{eq:error-main-regularity}, $C_{F,T}$, and the constants $C_0$ and
$C_{\mathrm{init}}$ of Theorem~\ref{thm:finite-time-temporal-error}, but
not on $k$, $N$, or $j$. No condition on $k$ other than $Nk\le T$ is
required. The constants deteriorate as a viscosity ratio approaches an
endpoint of \eqref{assumption:viscosity-fluctuation}, and can grow
exponentially in $T$ through $\exp(k\sum_n\Gamma_j^{n+1})$, with
$\Gamma_j^{n+1}$ from \eqref{eq:error-Gronwall-coefficient}.

The main result of this appendix bounds the velocity and
auxiliary-variable errors in terms of the first-level error and $k^4$;
with a suitable initialization, this gives second-order convergence.

\begin{theorem}
\label{thm:finite-time-temporal-error}
Let $d=2$, and let the first level $(\bm u_j^1,q_j^1)$ be computed by
\eqref{scheme:fsav-be-startup} or, more generally, satisfy
$\bm u_j^1\in\bm V$ and $\mathcal A_j^1\le C_0$ with
$C_0$ independent of $k$ and $j$.
Let $\bm F_j\in C([0,T];\bm X')$, and use the exact initial data
\eqref{eq:error-initial-values}.
Assume
\eqref{assumption:viscosity-fluctuation},
\eqref{assumption:discrete-force},
\eqref{eq:error-main-regularity}, and
\eqref{eq:error-force-consistency}. Then
\begin{align}
&\max_{1\le m\le N}\mathcal H_{e,j}^m
 +c_{\eta,j}\sum_{n=1}^{N-1}
  \|e_{\bm u,j}^{n+1}-2e_{\bm u,j}^{n}
    +e_{\bm u,j}^{n-1}\|^2
\nonumber\\
&\quad
 +\frac14\sum_{n=1}^{N-1}
  |e_{q,j}^{n+1}-2e_{q,j}^{n}+e_{q,j}^{n-1}|^2
 +\frac{k\bar\nu}{4}\sum_{n=1}^{N-1}
  \|\nabla\bm z_j^{n+1}\|^2
 +k\kappa\sigma\sum_{n=1}^{N-1}\mathcal F_{\nabla e,j}^{n+1}
\nonumber\\
&\quad
 +\frac{k\gamma}{2}\sum_{n=1}^{N-1}|e_{q,j}^{n+1}|^2
 \le C_T\bigl(\mathcal H_{e,j}^1+k^4\bigr).
\label{eq:error-main-estimate}
\end{align}
If the initialization satisfies
$\mathcal H_{e,j}^1\le C_{\mathrm{init}}k^4$, which holds for
\eqref{scheme:fsav-be-startup} when, in addition,
$\partial_t\bm u_j\in L^\infty(0,T;H^2(\Omega)^2)$ and
$\|\bm R_{F,j}^1\|\le Ck$ (Lemma~\ref{lem:error-startup}), then
\begin{equation}
\begin{aligned}
    &\max_{0\le n\le N}
      \bigl(\|e_{\bm u,j}^n\|^2+|e_{q,j}^n|^2\bigr)
      +k\sum_{n=1}^{N-1}
       \|\nabla e_{\bm u,j}^{n+1}\|^2
      \le C_Tk^4.
\end{aligned}
\label{eq:error-physical-estimate}
\end{equation}
Consequently,
\begin{equation}
    \max_{0\le n\le N}
    \bigl(\|e_{\bm u,j}^n\|+|e_{q,j}^n|\bigr)
    \le C_Tk^2.
    \label{eq:error-second-order-conclusion}
\end{equation}
\end{theorem}

The proof, given at the end of
Appendix~\ref{subsec:error-velocity-sav}, combines the residual and
startup estimates of
Appendix~\ref{subsec:error-consistency-startup}, the filtered
error-energy inequality of Lemma~\ref{lem:error-energy-recursion}, and
a discrete Gronwall argument.

\subsection{Consistency and the first time step}
\label{subsec:error-consistency-startup}

We first bound the residuals of the BDF2 difference quotient and
of the second-order extrapolation.

\begin{lemma}
\label{lem:error-consistency-residuals}
If $\partial_{ttt}\bm u_j\in L^2(0,T;\bm X')$, then%
\begin{equation}
    k\sum_{n=1}^{N-1}
    \|\bm R_{t,j}^{n+1}\|_{-1}^2
    \le
    Ck^4
    \|\partial_{ttt}\bm u_j\|_{L^2(0,T;\bm X')}^2.
    \label{eq:error-BDF2-residual-bound}
\end{equation}
If
$\partial_{tt}\bm u_j\in L^2(0,T;H^1(\Omega)^2)$, then
\begin{equation}
    k\sum_{n=1}^{N-1}
    \|\bm R_{\mathrm{ext},j}^{n+1}\|_{H^1}^2
    \le
    Ck^4
    \|\partial_{tt}\bm u_j\|_{L^2(0,T;H^1)}^2.
    \label{eq:error-extrapolation-residual-sum}
\end{equation}
If, in addition,
$\partial_{tt}\bm u_j\in L^\infty(0,T;H^1(\Omega)^2)$, then
\begin{equation}
    \max_{1\le n\le N-1}
    \|\bm R_{\mathrm{ext},j}^{n+1}\|_{H^1}^2
    \le
    Ck^4
    \|\partial_{tt}\bm u_j\|_{L^\infty(0,T;H^1)}^2.
    \label{eq:error-extrapolation-residual-max}
\end{equation}
\end{lemma}

\begin{proof}

The extrapolation is exact for linear and the BDF2 formula for
quadratic polynomials in time, so Taylor's formula with integral remainder
gives
\[
 \bm R_{\mathrm{ext},j}^{n+1}
 =\int_{t^{n-1}}^{t^{n+1}}L_n(s)\,\partial_{tt}\bm u_j(s)\,ds,
 \qquad
 \bm R_{t,j}^{n+1}
 =\int_{t^{n-1}}^{t^{n+1}}K_n(s)\,\partial_{ttt}\bm u_j(s)\,ds,
\]
with $L_n(s)=k-|s-t^n|$ and
\[
    K_n(s)=
    \begin{cases}
      \dfrac{(s-t^{n-1})^2}{4k}, & t^{n-1}\le s\le t^n,\\[1.2ex]
      \dfrac{(t^{n+1}-s)\bigl(4k-3(t^{n+1}-s)\bigr)}{4k}, & t^n\le s\le t^{n+1}.
    \end{cases}
\]
Thus $0\le L_n\le k$, $0\le K_n\le k/3$,
$\int_{t^{n-1}}^{t^{n+1}}L_n^2\,ds=2k^3/3$, and
$\int_{t^{n-1}}^{t^{n+1}}K_n^2\,ds=k^3/12$. The Cauchy--Schwarz
inequality for Bochner integrals yields
\[
 \|\bm R_{\mathrm{ext},j}^{n+1}\|_{H^1}^2
 \le\tfrac{2k^3}{3}\int_{t^{n-1}}^{t^{n+1}}
      \|\partial_{tt}\bm u_j\|_{H^1}^2\,ds,
 \quad
 \|\bm R_{t,j}^{n+1}\|_{-1}^2
 \le\tfrac{k^3}{12}\int_{t^{n-1}}^{t^{n+1}}
      \|\partial_{ttt}\bm u_j\|_{-1}^2\,ds.
\]

We multiply by $k$ and sum over $n$; since every point lies in at
most two of the intervals $[t^{n-1},t^{n+1}]$, this gives
\eqref{eq:error-BDF2-residual-bound} and
\eqref{eq:error-extrapolation-residual-sum}. Bounding the first integral
by $2k\|\partial_{tt}\bm u_j\|_{L^\infty(0,T;H^1)}^2$ gives
\eqref{eq:error-extrapolation-residual-max}.
\end{proof}

The stability of the mr-SAV--BE initialization
\eqref{scheme:fsav-be-startup} does not by itself provide the accuracy
needed to start the BDF2 error recursion; the next lemma shows that its
error is $O(k^2)$ in $L^2$. With $\bm R_{F,j}^1$ from
\eqref{eq:error-residuals}, set
\begin{equation}
    \bm\xi_j^1:=\bm U_j^1-\bm U_j^0,
    \qquad
    \bm R_{\mathrm{BE},j}^1:=\partial_t\bm u_j(t^1)-\frac{\bm\xi_j^1}{k},
\label{eq:error-startup-notation}
\end{equation}
and define the startup residuals
\begin{equation}
\begin{aligned}
    \bm{\mathscr R}_{u,j}^1
    &:=
    \bm R_{\mathrm{BE},j}^1
    +\bm R_{F,j}^1
    -\nu_j'\Delta\bm\xi_j^1
    +(\bm\xi_j^1\cdot\nabla)\bm U_j^1
    +(\bm U_j^0\cdot\nabla)\bm\xi_j^1,
    \\
    \mathscr R_{q,j}^1
    &:=
    b(\bm U_j^0,\bm U_j^0,\bm\xi_j^1).
\end{aligned}
\label{eq:error-startup-residuals}
\end{equation}
Since the initial errors vanish, the first-level energy is
\begin{equation}
 \mathcal H_{e,j}^1
 =\frac{5-3\eta_j}{4}\|e_{\bm u,j}^1\|^2
  +\frac54|e_{q,j}^1|^2
  +k\kappa(P_j)_{11}\|\nabla e_{\bm u,j}^1\|^2.
 \label{eq:error-startup-compensated-energy}
\end{equation}

The next lemma shows that the mr-SAV--BE step produces a
first-level error energy of order $k^4$.

\begin{lemma}
\label{lem:error-startup}
Let $\bm u_j^0=\bm U_j^0$ and $q_j^0=1$, and assume, in addition to
\eqref{eq:error-main-regularity},
$\partial_t\bm u_j\in L^\infty(0,T;H^2(\Omega)^2)$ and
$\|\bm R_{F,j}^1\|\le Ck$ uniformly in $j$ (trivially true for endpoint
sampling; for the cell average \eqref{def:discrete-force-average}, true
if $\partial_t\bm F_j\in L^\infty(0,T;L^2(\Omega)^2)$). Then, uniformly
in $j$,
\begin{equation}
    \|\bm{\mathscr R}_{u,j}^1\|
    +
    |\mathscr R_{q,j}^1|
    \le Ck,
    \label{eq:error-startup-residual-assumption}
\end{equation}
and the mr-SAV--BE initialization \eqref{scheme:fsav-be-startup}
satisfies
\begin{equation}
    \|e_{\bm u,j}^1\|^2
    +|e_{q,j}^1|^2
    +k\bar\nu\|\nabla e_{\bm u,j}^1\|^2
    +k\gamma|e_{q,j}^1|^2
    \le Ck^4.
\label{eq:error-startup-basic-bound}
\end{equation}
Consequently,
\begin{equation}
    \mathcal H_{e,j}^1
    \le Ck^4.
    \label{eq:error-startup-H-bound}
\end{equation}
\end{lemma}

\begin{proof}

\emph{Residuals.} By assumption,
$\|\bm\xi_j^1\|_{H^2}\le\int_0^k\|\partial_t\bm u_j\|_{H^2}\,ds\le Ck$,
and
$\bm R_{\mathrm{BE},j}^1=k^{-1}\int_0^k\int_s^k
\partial_{tt}\bm u_j(r)\,dr\,ds$, so that
$\|\bm R_{\mathrm{BE},j}^1\|
\le\frac k2\|\partial_{tt}\bm u_j\|_{L^\infty(0,T;L^2)}$. The
two-dimensional embedding $H^2(\Omega)\hookrightarrow L^\infty(\Omega)$
gives
\[
 \|(\bm\xi_j^1\cdot\nabla)\bm U_j^1\|
 +\|(\bm U_j^0\cdot\nabla)\bm\xi_j^1\|
 +|\mathscr R_{q,j}^1|
 \le C\bigl(\|\bm U_j^1\|_{H^2}+\|\bm U_j^0\|_{H^2}
            +\|\bm U_j^0\|_{H^2}^2\bigr)\|\bm\xi_j^1\|_{H^2}
 \le Ck.
\]
Together with $\|\Delta\bm\xi_j^1\|\le\sqrt2\,\|\bm\xi_j^1\|_{H^2}$ and
$\|\bm R_{F,j}^1\|\le Ck$, this proves
\eqref{eq:error-startup-residual-assumption}.

\emph{Error equations.} Since $\bm u_j^0=\bm U_j^0$ and
$\bm u_j^1,\bm U_j^1\in\bm V$, we have $e_{\bm u,j}^1\in\bm V$. Write the
exact weak momentum equation at $t^1$ with
$\bm\xi_j^1/k=\partial_t\bm u_j(t^1)-\bm R_{\mathrm{BE},j}^1$ and
$\nu_j\nabla\bm U_j^1=\bar\nu\nabla\bm U_j^1+\nu_j'\nabla\bm U_j^0
+\nu_j'\nabla\bm\xi_j^1$, and subtract it from
\eqref{scheme:fsav-be-startup-momentum}; insert $q_j^0=1$,
$q_j^1=1+e_{q,j}^1$, and $\bm u_j^1=\bm U_j^1+e_{\bm u,j}^1$ in
\eqref{scheme:fsav-be-startup-q}. For every $\bm v\in\bm V$ (the
pressure term vanishes),
\begin{align}
    \Bigl(\frac{e_{\bm u,j}^1}{k},\bm v\Bigr)
    +\bar\nu(\nabla e_{\bm u,j}^1,\nabla\bm v)
    &=b(\bm U_j^1,\bm U_j^1,\bm v)
      -(1+e_{q,j}^1)b(\bm U_j^0,\bm U_j^0,\bm v)
    \notag\\
    &\quad+(\bm R_{\mathrm{BE},j}^1+\bm R_{F,j}^1,\bm v)
      +\nu_j'(\nabla\bm\xi_j^1,\nabla\bm v),
    \label{eq:error-startup-momentum-error}\\
    \frac{e_{q,j}^1}{k}+\gamma e_{q,j}^1
    &=b(\bm U_j^0,\bm U_j^0,\bm U_j^1+e_{\bm u,j}^1).
    \label{eq:error-startup-scalar-error}
\end{align}

\emph{Energy identity.} Take $\bm v=ke_{\bm u,j}^1$ in
\eqref{eq:error-startup-momentum-error}, multiply
\eqref{eq:error-startup-scalar-error} by $ke_{q,j}^1$, and add; the
terms $\mp ke_{q,j}^1b(\bm U_j^0,\bm U_j^0,e_{\bm u,j}^1)$ cancel.
Moreover, by trilinearity,
\[
 b(\bm U_j^1,\bm U_j^1,e_{\bm u,j}^1)-b(\bm U_j^0,\bm U_j^0,e_{\bm u,j}^1)
 =b(\bm\xi_j^1,\bm U_j^1,e_{\bm u,j}^1)+b(\bm U_j^0,\bm\xi_j^1,e_{\bm u,j}^1);
\]
since $b(\bm U_j^0,\bm U_j^0,\bm U_j^0)=0$,
$b(\bm U_j^0,\bm U_j^0,\bm U_j^1)=\mathscr R_{q,j}^1$; and, since
$\nu_j'$ is constant and $\bm\xi_j^1\in H^2(\Omega)^2\cap\bm X$,
$\nu_j'(\nabla\bm\xi_j^1,\nabla e_{\bm u,j}^1)
=(-\nu_j'\Delta\bm\xi_j^1,e_{\bm u,j}^1)$. Hence
\begin{equation}
    \|e_{\bm u,j}^1\|^2+|e_{q,j}^1|^2
    +k\bar\nu\|\nabla e_{\bm u,j}^1\|^2+k\gamma|e_{q,j}^1|^2
    =k(\bm{\mathscr R}_{u,j}^1,e_{\bm u,j}^1)
     +ke_{q,j}^1\mathscr R_{q,j}^1.
    \label{eq:error-startup-energy-identity}
\end{equation}
Young's inequality bounds the right-hand side by
$\frac14\|e_{\bm u,j}^1\|^2+\frac14|e_{q,j}^1|^2
+k^2(\|\bm{\mathscr R}_{u,j}^1\|^2+|\mathscr R_{q,j}^1|^2)$, and
\eqref{eq:error-startup-residual-assumption} gives
\eqref{eq:error-startup-basic-bound}.

Since $(5-3\eta_j)/4<3/2$ and, by \eqref{eq:full-range-constants},
$\kappa(P_j)_{11}\le\kappa p_*=(1-\theta)\bar\nu/8<\bar\nu$,
\eqref{eq:error-startup-compensated-energy} then gives
$\mathcal H_{e,j}^1\le\frac32(\|e_{\bm u,j}^1\|^2+|e_{q,j}^1|^2)
+k\bar\nu\|\nabla e_{\bm u,j}^1\|^2\le Ck^4$, which is
\eqref{eq:error-startup-H-bound}.
\end{proof}

\subsection{Velocity and auxiliary-variable errors}
\label{subsec:error-velocity-sav}

The proof uses the filtered test of
Theorem~\ref{thm-long-time-stability}. Filtering produces one additional
nonlinear consistency term, which is controlled by the accumulated
fourth power of the $L^4$ norm of the numerical extrapolation, available
from the two-dimensional energy bounds.

The next lemma states the error equations and decomposes the
nonlinear error term that the filtered test produces.

\begin{lemma}
\label{lem:error-equations}
For $n\ge1$, the errors satisfy, for every $\bm v\in\bm V$,
\begin{align}
 (\delta_te_{\bm u,j}^{n+1},\bm v)
 &+\bar\nu(\nabla e_{\bm u,j}^{n+1},\nabla\bm v)
       +\nu_j'(\nabla e_{\bm u,j}^{*,n+1},\nabla\bm v)
 \nonumber\\
 &=b(\bm U_j^{n+1},\bm U_j^{n+1},\bm v)
   -q_j^{n+1}b(\bm u_j^{*,n+1},\bm u_j^{*,n+1},\bm v)
 \nonumber\\
 &\quad+\langle\bm R_{t,j}^{n+1}+\bm R_{F,j}^{n+1},\bm v\rangle
       +\nu_j'(\nabla\bm R_{\mathrm{ext},j}^{n+1},\nabla\bm v),
 \label{eq:error-weak-momentum}\\
 (\nabla\cdot e_{\bm u,j}^{n+1},r)&=0,\qquad r\in Q,
 \label{eq:error-divergence}\\
 \delta_te_{q,j}^{n+1}+\gamma e_{q,j}^{n+1}
 &=b(\bm u_j^{*,n+1},\bm u_j^{*,n+1},\bm u_j^{n+1}).
 \label{eq:error-weak-scalar}
\end{align}
Set
\begin{equation}
 \bm\zeta_j^{n+1}:=\bm U_j^{n+1}-\bm u_j^{*,n+1}
       =\bm R_{\mathrm{ext},j}^{n+1}-e_{\bm u,j}^{*,n+1}\in\bm V.
 \label{eq:error-zeta-definition}
\end{equation}
Testing momentum with $\bm z_j^{n+1}$ and the scalar equation with
$e_{q,j}^{n+1}$ gives the combined nonlinear contribution

\[
\begin{aligned}
 \mathscr N_j^{n+1}:={}&b(\bm U_j^{n+1},\bm U_j^{n+1},\bm z_j^{n+1})
 -q_j^{n+1}b(\bm u_j^{*,n+1},\bm u_j^{*,n+1},\bm z_j^{n+1})\\
 &+e_{q,j}^{n+1}b(\bm u_j^{*,n+1},\bm u_j^{*,n+1},\bm u_j^{n+1}),
\end{aligned}
\]
and
\begin{align}
 \mathscr N_j^{n+1}={}&
 b(\bm\zeta_j^{n+1},\bm U_j^{n+1},\bm z_j^{n+1})
 +b(\bm u_j^{*,n+1},\bm\zeta_j^{n+1},\bm z_j^{n+1})
 \nonumber\\
 &-(1-\eta_j)e_{q,j}^{n+1}
        b(\bm u_j^{*,n+1},\bm\zeta_j^{n+1},\bm U_j^{n+1})
 \nonumber\\
 &+\eta_j e_{q,j}^{n+1}
        b(\bm u_j^{*,n+1},\bm u_j^{*,n+1},
                                      \bm R_{\mathrm{ext},j}^{n+1}).
 \label{eq:error-nonlinear-decomposition}
\end{align}
\end{lemma}

\begin{proof}

Subtracting the exact weak momentum identity at $t^{n+1}$, written with
\eqref{eq:error-residuals} and $\nu_j=\bar\nu+\nu_j'$, from
\eqref{scheme:sav-bdf2-momentum} gives \eqref{eq:error-weak-momentum};
\eqref{eq:error-divergence} is immediate, and since the exact scalar
variable is $1$, \eqref{scheme:sav-bdf2-q} gives
\eqref{eq:error-weak-scalar}. All exact and numerical velocities,
extrapolations, and their differences lie in $\bm V$, so that
\begin{equation}
 b(\bm w,\bm a,\bm c)=-b(\bm w,\bm c,\bm a),\qquad
 b(\bm w,\bm a,\bm a)=0,\qquad \bm w\in\bm V.
 \label{eq:error-skew-identities}
\end{equation}

In this proof and the next, abbreviate
$\bm U=\bm U_j^{n+1}$, $\bm S=\bm U_j^{*,n+1}$, $\bm s=\bm u_j^{*,n+1}$,
$\bm v=\bm s-\bm S=e_{\bm u,j}^{*,n+1}$ (no longer a test function),
$\bm R=\bm U-\bm S=\bm R_{\mathrm{ext},j}^{n+1}$,
$\bm e=e_{\bm u,j}^{n+1}$, $\bm z=\bm e-\eta\bm v$,
$\bm\zeta=\bm U-\bm s=\bm R-\bm v$, $a=e_{q,j}^{n+1}$, and
$\eta=\eta_j$. With $q_j^{n+1}=1+a$ and $\bm u_j^{n+1}=\bm U+\bm e$,
\[
 \mathscr N=b(\bm U,\bm U,\bm z)-b(\bm s,\bm s,\bm z)
            +a\,b(\bm s,\bm s,\bm U+\eta\bm v),
\]
and $b(\bm U,\bm U,\bm z)-b(\bm s,\bm s,\bm z)
=b(\bm\zeta,\bm U,\bm z)+b(\bm s,\bm\zeta,\bm z)$. Since
$\bm U+\eta\bm v=(1-\eta)\bm U+\eta\bm R+\eta\bm s$,
$b(\bm s,\bm s,\bm s)=0$, and
$b(\bm s,\bm s,\bm U)=-b(\bm s,\bm\zeta,\bm U)$ by
\eqref{eq:error-skew-identities}, we obtain
$b(\bm s,\bm s,\bm U+\eta\bm v)
=-(1-\eta)b(\bm s,\bm\zeta,\bm U)+\eta\,b(\bm s,\bm s,\bm R)$, which
proves \eqref{eq:error-nonlinear-decomposition}.
\end{proof}

Testing the error equations with the filtered error
$\bm z_j^{n+1}$ gives the following one-step inequality for the error
energy.

\begin{lemma}
\label{lem:error-energy-recursion}
Under the assumptions of Theorem~\ref{thm:finite-time-temporal-error},
\begin{align}
 &\mathcal H_{e,j}^{n+1}-\mathcal H_{e,j}^n
       +\mathscr D_{e,j}^{n+1}
 \le k\Gamma_j^{n+1}\mathcal H_{e,j}^n+k\mathcal B_j^{n+1},
 \label{eq:error-fundamental-recursion}
\end{align}
where
\begin{align}
 \mathscr D_{e,j}^{n+1}:={}&
 c_{\eta,j}\|e_{\bm u,j}^{n+1}-2e_{\bm u,j}^n+e_{\bm u,j}^{n-1}\|^2
 +\frac14|e_{q,j}^{n+1}-2e_{q,j}^n+e_{q,j}^{n-1}|^2
 \nonumber\\
 &+\frac{k\bar\nu}{4}\|\nabla\bm z_j^{n+1}\|^2
       +k\kappa\sigma\mathcal F_{\nabla e,j}^{n+1}
       +\frac{k\gamma}{2}|e_{q,j}^{n+1}|^2,
 \label{eq:error-dissipation-functional}\\
 \Gamma_j^{n+1}:=C\Big[&\bar\nu^{-1}
   \bigl(\|\bm U_j^{n+1}\|_{H^2}^2
       +\|\bm U_j^{*,n+1}\|_{H^2}^2
       +\|\nabla e_{\bm u,j}^{*,n+1}\|^2\bigr)
 \nonumber\\
 &+\gamma^{-1}\|\bm U_j^{n+1}\|_{H^2}^2
                         \|\bm u_j^{*,n+1}\|_{L^4}^2\Big],
 \label{eq:error-Gronwall-coefficient}
\end{align}
and the nonnegative consistency collection is
\begin{align}
 \mathcal B_j^{n+1}:=C\Big[&\bar\nu^{-1}
      \bigl(\|\bm R_{t,j}^{n+1}\|_{-1}^2
                         +\|\bm R_{F,j}^{n+1}\|_{-1}^2\bigr)
 \nonumber\\
 &+\bar\nu^{-1}\|\bm U_j^{n+1}\|_{H^2}^2
                         \|\bm R_{\mathrm{ext},j}^{n+1}\|^2
 \nonumber\\
 &+\bar\nu^{-1}\bigl(\|\nabla\bm u_j^{*,n+1}\|^2+\rho^2\bigr)
                         \|\nabla\bm R_{\mathrm{ext},j}^{n+1}\|^2
 \nonumber\\
 &+\gamma^{-1}\|\bm U_j^{n+1}\|_{H^2}^2
       \|\bm u_j^{*,n+1}\|_{L^4}^2
                         \|\bm R_{\mathrm{ext},j}^{n+1}\|^2
 \nonumber\\
 &+\gamma^{-1}\|\bm u_j^{*,n+1}\|_{L^4}^4
                         \|\nabla\bm R_{\mathrm{ext},j}^{n+1}\|^2\Big].
 \label{eq:error-consistency-collection}
\end{align}
\end{lemma}

\begin{proof}

Test \eqref{eq:error-weak-momentum} with $k\bm z_j^{n+1}$ and
\eqref{eq:error-weak-scalar} with $ke_{q,j}^{n+1}$. Because
$\nu_j'=-\eta_j\bar\nu$, the diffusion terms combine to
$k\bar\nu\|\nabla\bm z_j^{n+1}\|^2$, and
\eqref{eq:full-range-G-identity} together with
\eqref{eq:bdf2-G-identity-preliminary} gives
\begin{align}
 &\mathcal E_{e,j}^{n+1}-\mathcal E_{e,j}^n
  +c_{\eta,j}\|e_{\bm u,j}^{n+1}-2e_{\bm u,j}^n+e_{\bm u,j}^{n-1}\|^2
 \nonumber\\
 &\quad+\frac14|e_{q,j}^{n+1}-2e_{q,j}^n+e_{q,j}^{n-1}|^2
       +k\bar\nu\|\nabla\bm z_j^{n+1}\|^2
       +k\gamma|e_{q,j}^{n+1}|^2
 \nonumber\\
 &=k\mathscr N_j^{n+1}
       +k\langle\bm R_{t,j}^{n+1}+\bm R_{F,j}^{n+1},\bm z_j^{n+1}\rangle
       +k\nu_j'(\nabla\bm R_{\mathrm{ext},j}^{n+1},
                                      \nabla\bm z_j^{n+1}),
 \label{eq:error-exact-balance}
\end{align}
where $\mathcal E_{e,j}^n:=\|E_{V,j}^n\|_{G_{\eta_j}}^2+|E_{Q,j}^n|_G^2$.
With the abbreviations of the previous proof, the Cauchy--Schwarz
inequality and \eqref{eq:error-filter-coercivity} give
$\|\bm v\|^2=\|2e_{\bm u,j}^n-e_{\bm u,j}^{n-1}\|^2
\le5(\|e_{\bm u,j}^n\|^2+\|e_{\bm u,j}^{n-1}\|^2)
\le5c_G^{-1}\mathcal H_{e,j}^n$, and $\|\bm\zeta\|\le\|\bm R\|+\|\bm v\|$.
We use
\begin{align}
 |b(\bm\zeta,\bm U,\bm z)|
 &\le C\|\bm U\|_{H^2}(\|\bm R\|+\|\bm v\|)\|\nabla\bm z\|,
 \label{eq:error-nonlinear-one}\\
 |b(\bm s,\bm\zeta,\bm z)|
 &\le C\Big[\|\nabla\bm s\|\|\nabla\bm R\|
        +\|\bm S\|_{H^2}\|\bm v\|
        +\|\bm v\|\|\nabla\bm v\|\Big]\|\nabla\bm z\|,
 \label{eq:error-nonlinear-two}\\
 |b(\bm s,\bm\zeta,\bm U)|
 &\le C\|\bm s\|_{L^4}\|\bm U\|_{H^2}(\|\bm R\|+\|\bm v\|),
 \qquad
 |b(\bm s,\bm s,\bm R)|
 \le\|\bm s\|_{L^4}^2\|\nabla\bm R\|.
 \label{eq:error-nonlinear-three}
\end{align}
For \eqref{eq:error-nonlinear-one}, write
$b(\bm\zeta,\bm U,\bm z)=-b(\bm\zeta,\bm z,\bm U)$ and use
$H^2\hookrightarrow L^\infty$. For \eqref{eq:error-nonlinear-two},
expand $b(\bm s,\bm\zeta,\bm z)=b(\bm s,\bm R,\bm z)
-b(\bm S,\bm v,\bm z)-b(\bm v,\bm v,\bm z)$ and use
\eqref{eq:trilinear-continuity},
$b(\bm S,\bm v,\bm z)=-b(\bm S,\bm z,\bm v)$ with
$H^2\hookrightarrow L^\infty$, and
$b(\bm v,\bm v,\bm z)=-b(\bm v,\bm z,\bm v)$ with Ladyzhenskaya's
inequality $\|\bm v\|_{L^4}^2\le C\|\bm v\|\|\nabla\bm v\|$.
\eqref{eq:error-nonlinear-three} follows from
\eqref{eq:error-skew-identities}, H\"older's inequality, and
$H^2\hookrightarrow W^{1,4}$.

Now apply Young's inequality to each term. The right-hand sides of
\eqref{eq:error-nonlinear-one}--\eqref{eq:error-nonlinear-two}, the two
dual residuals, and $|\nu_j'|\|\nabla\bm R\|\|\nabla\bm z\|$ with
$|\nu_j'|\le\rho$ contain eight products of the form
$X\|\nabla\bm z\|$; each is bounded by
$\frac{\bar\nu}{16}\|\nabla\bm z\|^2+4\bar\nu^{-1}X^2$. Since
$|\eta_j|<1$ and $0<1-\eta_j<4/3$, each of the two scalar terms is
bounded by $\frac\gamma4|a|^2$ plus $C\gamma^{-1}$ times the square of
its remaining factor, e.g.,
$k|\eta_j a\,b(\bm s,\bm s,\bm R)|\le\frac{k\gamma}{4}|a|^2
+k\gamma^{-1}\|\bm s\|_{L^4}^4\|\nabla\bm R\|^2$. Terms containing
$\|\bm v\|^2\le5c_G^{-1}\mathcal H_{e,j}^n$ form
$k\Gamma_j^{n+1}\mathcal H_{e,j}^n$, and the others form
$k\mathcal B_j^{n+1}$:
\begin{align}
 &k|\mathscr N_j^{n+1}|
 +k|\langle\bm R_{t,j}^{n+1}+\bm R_{F,j}^{n+1},\bm z_j^{n+1}\rangle|
 +k|\nu_j'(\nabla\bm R,\nabla\bm z)|
 \nonumber\\
 &\quad\le \frac{k\bar\nu}{2}\|\nabla\bm z\|^2
       +\frac{k\gamma}{2}|a|^2
       +k\Gamma_j^{n+1}\mathcal H_{e,j}^n+k\mathcal B_j^{n+1}.
 \label{eq:error-coercive-balance}
\end{align}
Insert \eqref{eq:error-coercive-balance} into
\eqref{eq:error-exact-balance} and add $k\kappa$ times
\eqref{eq:error-gradient-filter}. Since $k\kappa D=k\bar\nu/4$, half of
the remaining dissipation $\frac{k\bar\nu}{2}\|\nabla\bm z_j^{n+1}\|^2$
absorbs the right-hand side of the filter estimate, which proves
\eqref{eq:error-fundamental-recursion}.
\end{proof}

The next lemma shows that the coefficients of this inequality are
summable in time.

\begin{lemma}
\label{lem:error-coefficient-summability}
Under the assumptions of Theorem~\ref{thm:finite-time-temporal-error},
\begin{equation}
 k\sum_{n=1}^{N-1}\Gamma_j^{n+1}\le C_T,
 \label{eq:error-lambda-summability}
\end{equation}
and
\begin{equation}
 k\sum_{n=1}^{N-1}\mathcal B_j^{n+1}\le C_Tk^4.
 \label{eq:error-B-summability}
\end{equation}
\end{lemma}

\begin{proof}

We first bound the numerical solution, independently of the error
recursion. Adding $k\kappa$ times the gradient version of
\eqref{eq:full-range-filter-dissipation} for $V_j^n$ to
\eqref{eq:full-range-basic-step}, and dropping the second differences,
gives
\begin{equation}
 \mathcal A_j^{n+1}-\mathcal A_j^n
 +\frac{k\bar\nu}{4}\|\nabla\bm w_j^{n+1}\|^2
 +k\kappa\sigma\|\nabla V_j^{n+1}\|_{P_j}^2
 +\frac{k\gamma}{2}|q_j^{n+1}|^2
 \le k\widetilde C_F.
 \label{eq:finite-time-stability-error-use}
\end{equation}
For the backward Euler startup, the hypothesis $\mathcal A_j^1\le C_0$
holds: \eqref{eq:asymmetric-startup} bounds
$\|\bm u_j^1\|^2+|q_j^1|^2$ and $k\|\nabla\bm u_j^1\|^2$ uniformly for
$0<k\le T/2$, and $\|\bm u_j^0\|_{H^1}$, $G_{\eta_j}$, and $P_j$ are
bounded. Summing
\eqref{eq:finite-time-stability-error-use} yields
$\mathcal A_j^m+k\kappa\sigma\sum_{n=1}^{m-1}\|\nabla V_j^{n+1}\|_{P_j}^2
\le\mathcal A_j^1+T\widetilde C_F$ for $1\le m\le N$. Since $P_j\ge I$,
the levels $0,1$ of $k\sum_n\|\nabla\bm u_j^n\|^2$ are bounded by
$\mathcal A_j^1/\kappa$ and the levels $2,\dots,N$ by
$(\mathcal A_j^1+T\widetilde C_F)/(\kappa\sigma)$. With
\eqref{eq:full-range-physical-coercivity}, this gives
\begin{equation}
 \max_{0\le n\le N}(\|\bm u_j^n\|^2+|q_j^n|^2)
             +k\sum_{n=0}^{N}\|\nabla\bm u_j^n\|^2\le C_T.
 \label{eq:error-use-stability}
\end{equation}
By $\|2\bm a-\bm b\|^2\le5(\|\bm a\|^2+\|\bm b\|^2)$, Ladyzhenskaya's
inequality, and the discrete Cauchy--Schwarz inequality, we obtain
\begin{gather}
 \max_{1\le n<N}\|\bm u_j^{*,n+1}\|^2
       +k\sum_{n=1}^{N-1}\|\nabla\bm u_j^{*,n+1}\|^2\le C_T,
 \label{eq:error-extrapolated-stability}\\
 k\sum_{n=1}^{N-1}\|\bm u_j^{*,n+1}\|_{L^4}^2\le C_T,
 \label{eq:error-extrapolated-L4-summability}\\
 k\sum_{n=1}^{N-1}\|\bm u_j^{*,n+1}\|_{L^4}^4
 \le C\max_{1\le n<N}\|\bm u_j^{*,n+1}\|^2\,
              k\sum_{n=1}^{N-1}\|\nabla\bm u_j^{*,n+1}\|^2
 \le C_T.
 \label{eq:error-extrapolated-L4-fourth}
\end{gather}
Since $e_{\bm u,j}^{*,n+1}=\bm u_j^{*,n+1}-\bm U_j^{*,n+1}$, these bounds
and \eqref{eq:error-main-regularity} give
$k\sum_{n=1}^{N-1}\|\nabla e_{\bm u,j}^{*,n+1}\|^2\le C_T$ and
$\max_{1\le n<N}(\|\bm U_j^{n+1}\|_{H^2}+\|\bm U_j^{*,n+1}\|_{H^2})\le C$.
Together with \eqref{eq:error-extrapolated-L4-summability}, this proves
\eqref{eq:error-lambda-summability}, without using the error estimate
being proved.

By Lemma~\ref{lem:error-consistency-residuals},
\eqref{eq:error-force-consistency}, and the uniform $H^2$ bound of
$\bm U_j^{n+1}$, the terms of $\mathcal B_j^{n+1}$ without
solution-dependent weights sum to $C_Tk^4$. In each remaining term, we
bound $\|\bm R_{\mathrm{ext},j}^{n+1}\|_{H^1}^2$ by its maximum
\eqref{eq:error-extrapolation-residual-max} and sum the weight by
\eqref{eq:error-extrapolated-stability}--\eqref{eq:error-extrapolated-L4-fourth};
e.g.,
\[
 k\sum_{n=1}^{N-1}\|\bm u_j^{*,n+1}\|_{L^4}^4
       \|\nabla\bm R_{\mathrm{ext},j}^{n+1}\|^2
 \le\max_{1\le n<N}\|\bm R_{\mathrm{ext},j}^{n+1}\|_{H^1}^2\,
 k\sum_{n=1}^{N-1}\|\bm u_j^{*,n+1}\|_{L^4}^4\le C_Tk^4.
\]
This proves \eqref{eq:error-B-summability}.
\end{proof}

\begin{proof}[Proof of Theorem~\ref{thm:finite-time-temporal-error}]

Dropping $\mathscr D_{e,j}^{n+1}$ in
\eqref{eq:error-fundamental-recursion} gives
$\mathcal H_{e,j}^{n+1}\le(1+k\Gamma_j^{n+1})\mathcal H_{e,j}^n
+k\mathcal B_j^{n+1}$. Iterating and using $1+x\le e^x$, we obtain
$\mathcal H_{e,j}^{m+1}\le\exp\bigl(k\sum_{n=1}^m\Gamma_j^{n+1}\bigr)
\bigl(\mathcal H_{e,j}^1+k\sum_{n=1}^m\mathcal B_j^{n+1}\bigr)$ for
$1\le m<N$, and Lemma~\ref{lem:error-coefficient-summability} gives
\begin{equation}
 \max_{1\le n\le N}\mathcal H_{e,j}^n
       \le C_T(\mathcal H_{e,j}^1+k^4).
 \label{eq:error-compensated-energy-bound}
\end{equation}
Summing \eqref{eq:error-fundamental-recursion} over $n$ and using
$\mathcal H_{e,j}^N\ge0$ gives
\[
 \sum_{n=1}^{N-1}\mathscr D_{e,j}^{n+1}
 \le\mathcal H_{e,j}^1
   +\max_{1\le n\le N}\mathcal H_{e,j}^n\,k\sum_{n=1}^{N-1}\Gamma_j^{n+1}
   +k\sum_{n=1}^{N-1}\mathcal B_j^{n+1}
 \le C_T(\mathcal H_{e,j}^1+k^4),
\]
which proves \eqref{eq:error-main-estimate}. By \eqref{eq:error-filter-coercivity},
$e_{\bm u,j}^0=\bm0$, $e_{q,j}^0=0$, and
$\mathcal F_{\nabla e,j}^{n+1}\ge\|\nabla e_{\bm u,j}^{n+1}\|^2$, the
left-hand side of \eqref{eq:error-physical-estimate} is at most
$C_T(\mathcal H_{e,j}^1+k^4)$.

With $\mathcal H_{e,j}^1\le C_{\mathrm{init}}k^4$, this proves
\eqref{eq:error-physical-estimate} and
\eqref{eq:error-second-order-conclusion}. All recursion coefficients
multiply the previous-level energy, so no smallness condition on $k$ is
used in the Gronwall step.
\end{proof}

\subsection{Discrete time derivative and pressure error}
\label{subsec:error-time-derivative-pressure}

We bound the BDF2 difference quotient of the velocity error in
$\|\ell\|_{\bm V'}:=\sup_{\bm0\ne\bm v\in\bm V}
|\langle\ell,\bm v\rangle|/\|\nabla\bm v\|$, and then the pressure
error. The pressure estimate uses intermediate bounds from the proof,
not the $\bm V'$-bound itself (Remark~\ref{rem:Vprime-vs-Xprime}).

\begin{lemma}
\label{lem:error-discrete-time-derivative}
Under the assumptions of Theorem~\ref{thm:finite-time-temporal-error},
suppose that
\begin{equation}
    \mathcal H_{e,j}^1\le C_{\mathrm{init}}k^4,
    \label{eq:error-derivative-startup-assumption}
\end{equation}
which holds by \eqref{eq:error-startup-H-bound} under the startup
hypotheses of Lemma~\ref{lem:error-startup}. Then
\begin{equation}
    k\sum_{n=1}^{N-1}\|\delta_te_{\bm u,j}^{n+1}\|_{\bm V'}^2
    \le C_Tk^4.
    \label{eq:error-time-derivative-Vprime}
\end{equation}
\end{lemma}

\begin{proof}
By \eqref{eq:error-derivative-startup-assumption},
Theorem~\ref{thm:finite-time-temporal-error} gives
\eqref{eq:error-physical-estimate}. For $\bm v\in\bm X$, define
\begin{align}
    \langle\mathscr M_j^{n+1},\bm v\rangle
    &:=b(\bm U_j^{n+1},\bm U_j^{n+1},\bm v)
       -q_j^{n+1}b(\bm u_j^{*,n+1},\bm u_j^{*,n+1},\bm v)
    \notag\\
    &=b(\bm\zeta_j^{n+1},\bm U_j^{n+1},\bm v)
      +b(\bm u_j^{*,n+1},\bm\zeta_j^{n+1},\bm v)
      -e_{q,j}^{n+1}b(\bm u_j^{*,n+1},\bm u_j^{*,n+1},\bm v),
    \label{eq:error-momentum-nonlinear-functional}
\end{align}
with $\bm\zeta_j^{n+1}$ from \eqref{eq:error-zeta-definition}. Since the
first arguments $\bm\zeta_j^{n+1}$ and $\bm u_j^{*,n+1}$ belong to
$\bm V$, \eqref{eq:skew-general} moves $\bm v$ into the second slot of
each term. H\"older's inequality and the two-dimensional Ladyzhenskaya
inequality $\|\bm w\|_{L^4}^2\le C\|\bm w\|\|\nabla\bm w\|$ then give
\begin{align}
    \|\mathscr M_j^{n+1}\|_{\bm X'}
    \le C\Big[&
    \bigl(\|\bm\zeta_j^{n+1}\|\|\nabla\bm\zeta_j^{n+1}\|\bigr)^{1/2}
    \bigl(\|\bm U_j^{n+1}\|\|\nabla\bm U_j^{n+1}\|\bigr)^{1/2}
    \notag\\
    &+\bigl(\|\bm\zeta_j^{n+1}\|\|\nabla\bm\zeta_j^{n+1}\|\bigr)^{1/2}
    \bigl(\|\bm u_j^{*,n+1}\|\|\nabla\bm u_j^{*,n+1}\|\bigr)^{1/2}
    \notag\\
    &+|e_{q,j}^{n+1}|\,\|\bm u_j^{*,n+1}\|\|\nabla\bm u_j^{*,n+1}\|
    \Big].
    \label{eq:error-momentum-functional-bound}
\end{align}
Since $e_{\bm u,j}^{*,n+1}=2e_{\bm u,j}^n-e_{\bm u,j}^{n-1}$ and
$e_{\bm u,j}^0=\bm0$,
\begin{equation}
    k\sum_{n=1}^{N-1}\|\nabla e_{\bm u,j}^{*,n+1}\|^2
    \le10k\sum_{n=1}^{N-1}\|\nabla e_{\bm u,j}^{n}\|^2
    \le C_Tk^4,
    \label{eq:error-extrapolated-error-sum}
\end{equation}
where the level $n=1$ is bounded by
\eqref{eq:error-derivative-startup-assumption}, since
$k\|\nabla e_{\bm u,j}^1\|^2\le\kappa^{-1}\mathcal H_{e,j}^1$ by
\eqref{eq:error-filter-coercivity}, and the levels $n\ge2$ by
\eqref{eq:error-physical-estimate}. As
$\bm\zeta_j^{n+1}=\bm R_{\mathrm{ext},j}^{n+1}-e_{\bm u,j}^{*,n+1}$,
this estimate, \eqref{eq:error-physical-estimate},
\eqref{eq:error-extrapolation-residual-sum}, and
\eqref{eq:error-extrapolation-residual-max} give
\begin{equation}
    \max_{1\le n\le N-1}
    \bigl(\|\bm\zeta_j^{n+1}\|^2+|e_{q,j}^{n+1}|^2\bigr)
    +k\sum_{n=1}^{N-1}\|\nabla\bm\zeta_j^{n+1}\|^2
    \le C_Tk^4.
    \label{eq:error-zeta-control}
\end{equation}
Square \eqref{eq:error-momentum-functional-bound}, multiply by $k$, and
sum. Bounding $\|\bm U_j^{n+1}\|_{H^1}$ by the regularity
\eqref{eq:error-main-regularity}, and $\|\bm u_j^{*,n+1}\|$ and
$k\sum_n\|\nabla\bm u_j^{*,n+1}\|^2$ by
\eqref{eq:error-extrapolated-stability}, the Cauchy--Schwarz inequality
in $n$ bounds the result by
$C_T\bigl[\max_n\|\bm\zeta_j^{n+1}\|\,
(k\sum_n\|\nabla\bm\zeta_j^{n+1}\|^2)^{1/2}+\max_n|e_{q,j}^{n+1}|^2\bigr]$.
Hence \eqref{eq:error-zeta-control} yields
\begin{align}
    k\sum_{n=1}^{N-1}\|\mathscr M_j^{n+1}\|_{\bm X'}^2
    &\le C_Tk^4,
    \label{eq:error-momentum-functional-Xprime}\\
    k\sum_{n=1}^{N-1}\|\mathscr M_j^{n+1}\|_{\bm V'}^2
    &\le k\sum_{n=1}^{N-1}\|\mathscr M_j^{n+1}\|_{\bm X'}^2
    \le C_Tk^4.
    \label{eq:error-momentum-functional-sum}
\end{align}
Finally, \eqref{eq:error-weak-momentum} and $|\nu_j'|\le\rho$ give, for
every $\bm v\in\bm V$,
$|\langle\delta_te_{\bm u,j}^{n+1},\bm v\rangle|\le a_{n+1}\|\nabla\bm v\|$
with
\[
    a_{n+1}:=\bar\nu\|\nabla e_{\bm u,j}^{n+1}\|
    +\rho\|\nabla e_{\bm u,j}^{*,n+1}\|
    +\|\mathscr M_j^{n+1}\|_{\bm V'}
    +\|\bm R_{t,j}^{n+1}\|_{-1}
    +\|\bm R_{F,j}^{n+1}\|_{-1}
    +\rho\|\nabla\bm R_{\mathrm{ext},j}^{n+1}\|.
\]
Hence $\|\delta_te_{\bm u,j}^{n+1}\|_{\bm V'}\le a_{n+1}$, and
$k\sum_na_{n+1}^2\le C_Tk^4$ by \eqref{eq:error-physical-estimate},
\eqref{eq:error-extrapolated-error-sum},
\eqref{eq:error-momentum-functional-sum},
\eqref{eq:error-BDF2-residual-bound},
\eqref{eq:error-force-consistency}, and
\eqref{eq:error-extrapolation-residual-sum}.
\end{proof}

To recover the pressure through the inf--sup condition
\eqref{eq:continuous-infsup}, the momentum error equation must be tested
with all $\bm v\in\bm X$. On a no-slip domain the $\bm X'$-norm of an
element of $\bm V$ is not controlled by its $\bm V'$-norm:

\begin{remark}[The norms $\|\cdot\|_{\bm V'}$ and $\|\cdot\|_{\bm X'}$ on $\bm V$]
\label{rem:Vprime-vs-Xprime}
Identify $\bm g\in\bm V$ with $\bm v\mapsto(\bm g,\bm v)$. Clearly
$\|\bm g\|_{\bm V'}\le\|\bm g\|_{\bm X'}$, but no reverse bound
$\|\bm g\|_{\bm X'}\le C\|\bm g\|_{\bm V'}$ holds. Let
$\Omega\subset\mathbb R^2$ be a bounded Lipschitz domain,
$\nabla^\perp\phi:=(\partial_y\phi,-\partial_x\phi)^\top$, and
$\operatorname{curl}\bm v:=\partial_xv_2-\partial_yv_1$. For
$\phi\in H^1(\Omega)$ and $\bm v\in\bm X$, one has
$(\nabla^\perp\phi,\bm v)=(\phi,\operatorname{curl}\bm v)$ and
$\|\operatorname{curl}\bm v\|\le\|\nabla\bm v\|$. Take
$\theta_\varepsilon\in C_c^\infty(\Omega)$ with
$0\le\theta_\varepsilon\le1$ and $\theta_\varepsilon\to1$ a.e., and set
$\bm g_\varepsilon:=\nabla^\perp(x\theta_\varepsilon)\in\bm V$. Since
$\nabla^\perp x=-\nabla y$ is $L^2$-orthogonal to $\bm V$,
$(\bm g_\varepsilon,\bm v)=(x(\theta_\varepsilon-1),\operatorname{curl}\bm v)$
for $\bm v\in\bm V$, so
$\|\bm g_\varepsilon\|_{\bm V'}\le\|x(\theta_\varepsilon-1)\|\to0$. On
the other hand, $(\nabla^\perp\phi,\bm v)=(\nabla\phi,\mathsf R\bm v)$
with the rotation $\mathsf R\bm v:=(-v_2,v_1)^\top$, an isometry of
$\bm X$, so \eqref{eq:continuous-infsup} gives
$\|\bm g_\varepsilon\|_{\bm X'}=\|\nabla(x\theta_\varepsilon)\|_{-1}
\ge\beta_0\|x\theta_\varepsilon\|_{L^2(\Omega)/\mathbb R}
\to\beta_0\|x\|_{L^2(\Omega)/\mathbb R}>0$. Hence
Lemma~\ref{lem:error-discrete-time-derivative} does not bound
$\delta_te_{\bm u,j}^{n+1}$ in $\bm X'$, and does not yield a bound for
$k\sum_n\|e_{p,j}^{n+1}\|^2$.
\end{remark}

We therefore estimate the pressure error in a weaker, time-integrated
norm. As noted in Appendix~\ref{subsec:error-setting}, the exact
velocity satisfies the divergence-free weak form at every node, so the
functional
$\bm v\mapsto\langle\bm F_j(t^m),\bm v\rangle-(\partial_t\bm u_j(t^m),\bm v)
-b(\bm U_j^m,\bm U_j^m,\bm v)-\nu_j(\nabla\bm U_j^m,\nabla\bm v)$
lies in $\bm X'$ and vanishes on $\bm V$. By
\eqref{eq:continuous-infsup}, there is a unique $p_j(t^m)\in Q$ for which
the momentum identity holds in $\bm X'$ at $t^m$. For $n=1,\ldots,N$, set
\begin{equation}
    e_{p,j}^n:=p_j^n-p_j(t^n).
    \label{eq:error-pressure-notation}
\end{equation}

The next corollary bounds the time-integrated pressure error at
second order.

\begin{corollary}
\label{cor:error-pressure}
Assume the hypotheses of Theorem~\ref{thm:finite-time-temporal-error}
and \eqref{eq:error-derivative-startup-assumption}. Then
\begin{equation}
    \max_{1\le n\le N-1}
    \Bigl\|k\sum_{m=2}^{n+1}e_{p,j}^{m}\Bigr\|
    \le C_Tk^2,
    \label{eq:error-pressure-estimate}
\end{equation}
where $C_T$ may also depend on $\beta_0^{-1}$, but not on $k$, $N$, or
$j$. If the additional startup hypotheses of
Lemma~\ref{lem:error-startup} hold, then $\|e_{p,j}^1\|\le Ck$, and
\eqref{eq:error-pressure-estimate} also holds with the summation starting
at $m=1$.
\end{corollary}

\begin{proof}
Subtracting the exact momentum identity at $t^m$ from
\eqref{scheme:sav-bdf2-momentum} gives, for every $\bm v\in\bm X$ and
$m=2,\ldots,N$,
\begin{align}
    (\delta_te_{\bm u,j}^{m},\bm v)
    &+\bar\nu(\nabla e_{\bm u,j}^{m},\nabla\bm v)
    +\nu_j'(\nabla e_{\bm u,j}^{*,m},\nabla\bm v)
    -(e_{p,j}^{m},\nabla\cdot\bm v)
    \notag\\
    &=
    \langle\mathscr M_j^{m}+\bm R_{t,j}^{m}+\bm R_{F,j}^{m},\bm v\rangle
    +\nu_j'(\nabla\bm R_{\mathrm{ext},j}^{m},\nabla\bm v),
    \label{eq:error-pressure-momentum-equation}
\end{align}
with $\mathscr M_j^m$ from \eqref{eq:error-momentum-nonlinear-functional},
which is bounded in $\bm X'$ by \eqref{eq:error-momentum-functional-Xprime}.
Multiply by $k$ and sum over $m=2,\ldots,n+1$. Since
$e_{\bm u,j}^0=\bm0$, the BDF2 differences telescope:
\begin{equation}
    k\sum_{m=2}^{n+1}\delta_te_{\bm u,j}^{m}
    =\tfrac12\bigl(3e_{\bm u,j}^{n+1}-e_{\bm u,j}^{n}-3e_{\bm u,j}^{1}\bigr).
    \label{eq:bdf2-telescoping}
\end{equation}
Hence, for every $\bm v\in\bm X$,
\begin{align*}
    \Bigl(k\sum_{m=2}^{n+1}e_{p,j}^{m},\nabla\cdot\bm v\Bigr)
    ={}&\tfrac12\bigl(3e_{\bm u,j}^{n+1}-e_{\bm u,j}^{n}-3e_{\bm u,j}^{1},\bm v\bigr)
    \\
    &+k\sum_{m=2}^{n+1}
    \Bigl[\bar\nu(\nabla e_{\bm u,j}^{m},\nabla\bm v)
        +\nu_j'(\nabla e_{\bm u,j}^{*,m},\nabla\bm v)\Bigr]
    \\
    &-k\sum_{m=2}^{n+1}
    \Bigl[\langle\mathscr M_j^{m}+\bm R_{t,j}^{m}+\bm R_{F,j}^{m},\bm v\rangle
        +\nu_j'(\nabla\bm R_{\mathrm{ext},j}^{m},\nabla\bm v)\Bigr].
\end{align*}
By \eqref{eq:poincare},
$|(\bm w,\bm v)|\le c_0^{-1/2}\|\bm w\|\|\nabla\bm v\|$. The inf--sup
condition \eqref{eq:continuous-infsup} and $|\nu_j'|\le\rho$ therefore give
\[
    \beta_0\Bigl\|k\sum_{m=2}^{n+1}e_{p,j}^{m}\Bigr\|
    \le\frac{7}{2}c_0^{-1/2}\max_{0\le m\le N}\|e_{\bm u,j}^m\|
    +k\sum_{m=2}^{N}a_m,
\]
with $a_m$ as in the proof of
Lemma~\ref{lem:error-discrete-time-derivative}, but with
$\|\mathscr M_j^m\|_{\bm X'}$ in place of $\|\mathscr M_j^m\|_{\bm V'}$.
The first term is at most $C_Tk^2$ by
\eqref{eq:error-physical-estimate}. By the Cauchy--Schwarz inequality,
$k\sum_ma_m\le T^{1/2}(k\sum_ma_m^2)^{1/2}$, and $k\sum_ma_m^2\le C_Tk^4$
by the same estimates as in that proof, with
\eqref{eq:error-momentum-functional-Xprime} for the nonlinear term. This proves
\eqref{eq:error-pressure-estimate}.

For the first level, subtract the exact momentum identity at $t^1$ from
\eqref{scheme:fsav-be-startup-momentum} and rewrite the result with
$\bm{\mathscr R}_{u,j}^1$ from \eqref{eq:error-startup-residuals}, as in
the proof of Lemma~\ref{lem:error-startup}. For all $\bm v\in\bm X$,
\[
  (e_{p,j}^1,\nabla\cdot\bm v)
  =\Bigl(\frac{e_{\bm u,j}^1}{k}-\bm{\mathscr R}_{u,j}^1,\bm v\Bigr)
  +\bar\nu(\nabla e_{\bm u,j}^1,\nabla\bm v)
  +e_{q,j}^1\,b(\bm U_j^0,\bm U_j^0,\bm v).
\]
By \eqref{eq:continuous-infsup}, \eqref{eq:poincare},
\eqref{eq:trilinear-continuity}, \eqref{eq:error-startup-residual-assumption},
and \eqref{eq:error-startup-basic-bound},
\[
  \beta_0\|e_{p,j}^1\|
  \le c_0^{-1/2}\Bigl(\frac{\|e_{\bm u,j}^1\|}{k}+\|\bm{\mathscr R}_{u,j}^1\|\Bigr)
  +\bar\nu\|\nabla e_{\bm u,j}^1\|
  +C_b|e_{q,j}^1|\,\|\nabla\bm U_j^0\|^2
  \le Ck .
\]
Hence $\|ke_{p,j}^1\|\le Ck^2$, and the triangle inequality gives the
last assertion.
\end{proof}

\begin{remark}[On the pressure estimate]
\label{rem:pressure-norm}
Estimate \eqref{eq:error-pressure-estimate} is a second-order bound for
the discrete time primitive of the pressure error, that is, for the
pressure error in a discrete $W^{-1,\infty}(0,T;L^2(\Omega))$-type
norm. In particular, for
$1\le n_1<n_2\le N$,
\[
  \Bigl\|\frac{k}{t^{n_2}-t^{n_1}}\sum_{m=n_1+1}^{n_2}
  \bigl(p_j^m-p_j(t^m)\bigr)\Bigr\|
  \le\frac{2C_Tk^2}{t^{n_2}-t^{n_1}},
\]
so discrete time averages over intervals of fixed length converge with
second order to the corresponding averages of the nodal values
$p_j(t^m)$. Differences of consecutive partial sums give
$\max_{2\le m\le N}\|e_{p,j}^m\|\le2C_Tk$. A second-order bound in the
discrete $L^2(0,T;L^2(\Omega))$ norm would require an estimate of
$\delta_te_{\bm u,j}^{n+1}$ in $\bm X'$ (or in $L^2(\Omega)^2$), which is
not proved here (Remark~\ref{rem:Vprime-vs-Xprime}). It would also
require $\|e_{p,j}^1\|=O(k^{3/2})$: an $O(k)$ first pressure error
contributes $O(k^{3/2})$ to that norm. With the mr-SAV--BE startup,
$\|e_{p,j}^1\|\le Ck$ cannot be expected to improve in general, because
$p_j^1$ responds without damping to the non-solenoidal part of the
explicitly treated terms in $\bm{\mathscr R}_{u,j}^1=O(k)$; a
second-order starting pressure would remove this restriction. In
\eqref{eq:error-pressure-estimate}, the first level enters only through
$ke_{p,j}^1$, so first-order accuracy of $p_j^1$ suffices there.
\end{remark}

\begin{remark}[Spatially varying viscosity]
\label{rem:spatial-error}
The temporal error analysis extends to the scheme
\eqref{scheme:random-viscosity} in both cases of
Theorem~\ref{thm:random-viscosity-stability}, for either viscosity form
in \eqref{def:viscous-strain} and sufficiently regular solutions of the
corresponding continuous problem, provided
$\nu_j,\bar\nu\in W^{1,\infty}(\Omega)$; this regularity enters only the
consistency of the startup step, where $\nu_j'\Delta\bm\xi_j^1$ in
\eqref{eq:error-startup-residuals} becomes a divergence of
$\nu_j'\bm\varepsilon(\bm\xi_j^1)$. For proportional fields
(case~(ii), \eqref{assumption:proportional-viscosity}), $\eta_j$ is
constant, so the filtered error test applies as above, with $\bar\nu$
replaced by $\operatorname*{ess\,inf}_\Omega\bar\nu$ in the lower bound
of the dissipation, in $\kappa$, and in the constants. Under the signed
condition (case~(i), \eqref{assumption:spatial-signed-fluctuation}), no
filter is used: inequality \eqref{eq:spatial-signed-compensation}
follows from pointwise identities and therefore holds for the error
sequence, so the modified error energy, with the correction term
\eqref{eq:spatial-signed-compensator} evaluated on the errors, satisfies
a recursion of the same form, and the one-step inequality
\eqref{eq:spatial-signed-one-step} supplies the required stability
bounds. We do not pursue the details here.
\end{remark}


\begin{thebibliography}{10}

\bibitem{AkrivisLubich2015}
{\sc G.~Akrivis and C.~Lubich}, {\em Fully implicit, linearly implicit and
  implicit--explicit backward difference formulae for quasi-linear parabolic
  equations}, Numerische Mathematik, 131 (2015), pp.~713--735,
  \url{https://doi.org/10.1007/s00211-015-0702-0}.

\bibitem{Bansal2021}
{\sc P.~Bansal}, {\em Numerical approximation of statistical solutions of the
  incompressible {N}avier--{S}tokes equations}.
\newblock arXiv:2107.06073, 2021.

\bibitem{Dahlquist1976}
{\sc G.~Dahlquist}, {\em Error analysis for a class of methods for stiff
  non-linear initial value problems}, in Numerical Analysis, vol.~506 of
  Lecture Notes in Mathematics, Springer, 1976, pp.~60--72,
  \url{https://doi.org/10.1007/BFb0080115}.

\bibitem{Fiordilino2018}
{\sc J.~A. Fiordilino}, {\em A second order ensemble timestepping algorithm for
  natural convection}, SIAM Journal on Numerical Analysis, 56 (2018),
  pp.~816--837, \url{https://doi.org/10.1137/17M1135104}.

\bibitem{GTWWW2012}
{\sc S.~Gottlieb, F.~Tone, C.~Wang, X.~Wang, and D.~Wirosoetisno}, {\em Long
  time stability of a classical efficient scheme for two-dimensional
  {Navier--Stokes} equations}, SIAM Journal on Numerical Analysis, 50 (2012),
  pp.~126--150, \url{https://doi.org/10.1137/110834901}.

\bibitem{GunzburgerJiangSchneier2017}
{\sc M.~Gunzburger, N.~Jiang, and M.~Schneier}, {\em An ensemble-proper
  orthogonal decomposition method for the nonstationary {Navier--Stokes}
  equations}, SIAM Journal on Numerical Analysis, 55 (2017), pp.~286--304,
  \url{https://doi.org/10.1137/16M1056444}.

\bibitem{GunzburgerJiangWang2019IMA}
{\sc M.~Gunzburger, N.~Jiang, and Z.~Wang}, {\em An efficient algorithm for
  simulating ensembles of parameterized flow problems}, IMA Journal of
  Numerical Analysis, 39 (2019), pp.~1180--1205,
  \url{https://doi.org/10.1093/imanum/dry029}.

\bibitem{GJW-CMAM2019}
{\sc M.~Gunzburger, N.~Jiang, and Z.~Wang}, {\em A second-order time-stepping
  scheme for simulating ensembles of parameterized flow problems},
  Computational Methods in Applied Mathematics, 19 (2019), pp.~681--701,
  \url{https://doi.org/10.1515/cmam-2017-0051}.

\bibitem{HairerWanner1996}
{\sc E.~Hairer and G.~Wanner}, {\em Solving Ordinary Differential Equations II:
  Stiff and Differential-Algebraic Problems}, vol.~14 of Springer Series in
  Computational Mathematics, Springer, Berlin, Heidelberg, 2nd~ed., 1996,
  \url{https://doi.org/10.1007/978-3-642-05221-7}.

\bibitem{HanWang-M3AS}
{\sc D.~Han and X.~Wang}, {\em A highly efficient second-order accurate
  long-time dynamics preserving scheme for some geophysical fluid models}.
\newblock arXiv preprint arXiv:2510.05360, 2025,
  \url{https://arxiv.org/abs/2510.05360}.

\bibitem{HeywoodRannacher1982}
{\sc J.~G. Heywood and R.~Rannacher}, {\em Finite element approximation of the
  nonstationary {Navier--Stokes} problem. {I}. {R}egularity of solutions and
  second-order error estimates for spatial discretization}, SIAM Journal on
  Numerical Analysis, 19 (1982), pp.~275--311,
  \url{https://doi.org/10.1137/0719018}.

\bibitem{Hosder2006UQ}
{\sc S.~Hosder, R.~W. Walters, and R.~Perez}, {\em A non-intrusive polynomial
  chaos method for uncertainty propagation in {CFD} simulations}, in 44th AIAA
  Aerospace Sciences Meeting and Exhibit, American Institute of Aeronautics and
  Astronautics, 2006, \url{https://doi.org/10.2514/6.2006-891}.

\bibitem{HuangShen2021}
{\sc F.~Huang and J.~Shen}, {\em Stability and error analysis of a class of
  high-order {IMEX} schemes for {Navier--Stokes} equations with periodic
  boundary conditions}, SIAM Journal on Numerical Analysis, 59 (2021),
  pp.~2926--2954, \url{https://doi.org/10.1137/21M1404144}.

\bibitem{Jiang2015}
{\sc N.~Jiang}, {\em A higher order ensemble simulation algorithm for fluid
  flows}, Journal of Scientific Computing, 64 (2015), pp.~264--288,
  \url{https://doi.org/10.1007/s10915-014-9932-z}.

\bibitem{JiangLayton2014}
{\sc N.~Jiang and W.~Layton}, {\em An algorithm for fast calculation of flow
  ensembles}, International Journal for Uncertainty Quantification, 4 (2014),
  pp.~273--301.

\bibitem{JiangLayton2015}
{\sc N.~Jiang and W.~Layton}, {\em Numerical analysis of two ensemble eddy
  viscosity numerical regularizations of fluid motion}, Numerical Methods for
  Partial Differential Equations, 31 (2015), pp.~630--651,
  \url{https://doi.org/10.1002/num.21908}.

\bibitem{JiangYang2021CMAME}
{\sc N.~Jiang and H.~Yang}, {\em {SAV} decoupled ensemble algorithms for fast
  computation of {Stokes--Darcy} flow ensembles}, Computer Methods in Applied
  Mechanics and Engineering, 387 (2021), p.~114150,
  \url{https://doi.org/10.1016/j.cma.2021.114150}.

\bibitem{JiangYang2021SISC}
{\sc N.~Jiang and H.~Yang}, {\em Stabilized scalar auxiliary variable ensemble
  algorithms for parameterized flow problems}, SIAM Journal on Scientific
  Computing, 43 (2021), pp.~A2869--A2896,
  \url{https://doi.org/10.1137/20M1364679}.

\bibitem{Kelliher2010}
{\sc J.~P. Kelliher}, {\em Eigenvalues of the {S}tokes operator versus the
  {D}irichlet {L}aplacian in the plane}, Pacific Journal of Mathematics, 244
  (2010), pp.~99--132, \url{https://doi.org/10.2140/pjm.2010.244.99}.

\bibitem{LancasterTismenetsky1985}
{\sc P.~Lancaster and M.~Tismenetsky}, {\em The Theory of Matrices}, Academic
  Press, Orlando, FL, 2nd~ed., 1985.

\bibitem{LiShen2020}
{\sc X.~Li and J.~Shen}, {\em Error analysis of the {SAV-MAC} scheme for the
  {Navier--Stokes} equations}, SIAM Journal on Numerical Analysis, 58 (2020),
  pp.~2465--2491, \url{https://doi.org/10.1137/19M1288267}.

\bibitem{LiShenLiu2022}
{\sc X.~Li, J.~Shen, and Z.~Liu}, {\em New {SAV}-pressure correction methods
  for the {Navier--Stokes} equations: Stability and error analysis},
  Mathematics of Computation, 91 (2022), pp.~141--167,
  \url{https://doi.org/10.1090/mcom/3651}.

\bibitem{LinYangDong2019}
{\sc L.~Lin, Z.~Yang, and S.~Dong}, {\em Numerical approximation of
  incompressible {Navier--Stokes} equations based on an auxiliary energy
  variable}, Journal of Computational Physics, 388 (2019), pp.~1--22,
  \url{https://doi.org/10.1016/j.jcp.2019.03.012}.

\bibitem{LuoWang2018}
{\sc Y.~Luo and Z.~Wang}, {\em An ensemble algorithm for numerical solutions to
  deterministic and random parabolic {PDE}s}, SIAM Journal on Numerical
  Analysis, 56 (2018), pp.~859--876, \url{https://doi.org/10.1137/17M1131489}.

\bibitem{Miller1971}
{\sc J.~J.~H. Miller}, {\em On the location of zeros of certain classes of
  polynomials with applications to numerical analysis}, IMA Journal of Applied
  Mathematics, 8 (1971), pp.~397--406,
  \url{https://doi.org/10.1093/imamat/8.3.397}.

\bibitem{MohebujjamanRebholz2017}
{\sc M.~Mohebujjaman and L.~G. Rebholz}, {\em An efficient algorithm for
  computation of {MHD} flow ensembles}, Computational Methods in Applied
  Mathematics, 17 (2017), pp.~121--137,
  \url{https://doi.org/10.1515/cmam-2016-0033}.

\bibitem{NevanlinnaOdeh1981}
{\sc O.~Nevanlinna and F.~Odeh}, {\em Multiplier techniques for linear
  multistep methods}, Numerical Functional Analysis and Optimization, 3 (1981),
  pp.~377--423, \url{https://doi.org/10.1080/01630568108816097}.

\bibitem{OppenheimSchafer2010}
{\sc A.~V. Oppenheim and R.~W. Schafer}, {\em Discrete-Time Signal Processing},
  Pearson, Upper Saddle River, NJ, 3rd~ed., 2010.

\bibitem{ShenXuYang2018}
{\sc J.~Shen, J.~Xu, and J.~Yang}, {\em The scalar auxiliary variable ({SAV})
  approach for gradient flows}, Journal of Computational Physics, 353 (2018),
  pp.~407--416, \url{https://doi.org/10.1016/j.jcp.2017.10.021}.

\bibitem{ShenXuYang2019}
{\sc J.~Shen, J.~Xu, and J.~Yang}, {\em A new class of efficient and robust
  energy stable schemes for gradient flows}, SIAM Review, 61 (2019),
  pp.~474--506, \url{https://doi.org/10.1137/17M1150153}.

\bibitem{Temam1997}
{\sc R.~Temam}, {\em Infinite-Dimensional Dynamical Systems in Mechanics and
  Physics}, vol.~68 of Applied Mathematical Sciences, Springer, New York,
  2~ed., 1997.

\bibitem{Wang2010a}
{\sc X.~Wang}, {\em Approximation of stationary statistical properties of
  dissipative dynamical systems: time discretization}, Math. Comp., 79 (2010),
  pp.~259--280, \url{https://doi.org/10.1090/S0025-5718-09-02256-X}.

\bibitem{Wang2012}
{\sc X.~Wang}, {\em An efficient second order in time scheme for approximating
  long time statistical properties of the two dimensional {Navier--Stokes}
  equations}, Numerische Mathematik, 121 (2012), pp.~753--779.

\bibitem{Wang2016}
{\sc X.~Wang}, {\em Numerical algorithms for stationary statistical properties
  of dissipative dynamical systems}, Discrete Contin. Dyn. Syst., 36 (2016),
  pp.~4599--4618, \url{https://doi.org/10.3934/dcds.2016.36.4599}.

\bibitem{XiuHesthaven2005}
{\sc D.~Xiu and J.~S. Hesthaven}, {\em High-order collocation methods for
  differential equations with random inputs}, SIAM Journal on Scientific
  Computing, 27 (2005), pp.~1118--1139,
  \url{https://doi.org/10.1137/040615201}.

\end{thebibliography}
\end{document}